\documentclass[12pt,a4paper,reqno]{amsart}
\pdfoutput=1
\usepackage[hidelinks]{hyperref}
\usepackage[top=2.5cm,right=2.2cm,left=2.2cm,bottom=2cm,headsep=20pt]{geometry}
\usepackage{amsfonts,amsthm,amssymb,microtype}
\usepackage[justification=centering,skip=0pt]{caption}
\usepackage{enumitem}
\usepackage[new]{old-arrows}
\usepackage{boondox-cal} %for calligraphics
\usepackage{mathtools}
\mathtoolsset{showonlyrefs} %shows only the labels of referenced equations
\usepackage{tikz}
\tikzset{
	node distance=1.8cm,
	main node/.style={circle,inner sep=2pt},
	freccia/.style={->,shorten >=1pt,shorten <=1pt},
	ciclo/.style={out=130, in=50, loop, distance=2cm, ->},
	line width=0.5pt }
\usetikzlibrary{arrows,cd}
\tikzcdset{every arrow/.append style = -{>[scale=0.85]}}

\makeatletter
\newcommand\myitem[1][]{\item[#1]\def\@currentlabel{#1}} %to put \label on axioms
\makeatother

\newcommand{\cof}{ \longhookrightarrow }
\newcommand{\weq}{ \stackrel{\text{\raisebox{-3pt}{\scriptsize $\sim$}}}{\longrightarrow} }

\makeatletter
\newcommand{\doublewidetilde}[1]{{%
  \mathpalette\double@widetilde{#1}%
}}
\newcommand{\double@widetilde}[2]{%
  \sbox\z@{$\m@th{\hspace{-3pt}#1}\widetilde{#2}$}%
  \ht\z@=.9\ht\z@
  \widetilde{\box\z@}%
}
\makeatother

\allowdisplaybreaks[4]

\newtheorem{thm}{Theorem}[section]
\newtheorem{prop}[thm]{Proposition}
\newtheorem{cor}[thm]{Corollary}
\newtheorem{lemma}[thm]{Lemma}
\theoremstyle{remark}
\newtheorem{rmk}[thm]{Remark}
\newtheorem{example}[thm]{Example}
\theoremstyle{definition}
\newtheorem{df}[thm]{Definition}
\DeclareMathOperator{\Ho}{Ho}

\newcommand{\id}{\mathrm{id}}
\newcommand{\CP}{C\hspace{-1pt}P}
\newcommand{\CPq}{\mathsf{CP}}
\newcommand{\RP}{RP}
\newcommand{\RPq}{\mathsf{RP}}
\newcommand{\WP}{\mathsf{W}\hspace{-1pt}\mathsf{P}}
\newcommand{\tub}[2]{
{\rm Tub\hspace{0.15em}}_{\CPq_{H}^{{#2}}}\hspace{-0.1em}(\CPq_{H}^{#1})}
\newcommand{\ctub}[2]{
{\rm Tub\hspace{0.15em}}_{\CP^{{#2}}}\hspace{-0.1em}(\CP^{#1})}

\newcommand{\CHS}{\ensuremath{\mathcal{C\hspace{-4.2pt}p\hspace{-2.5pt}c\hspace{-2.2pt}t\hspace{-1.7pt}H\hspace{-5.3pt}a\hspace{-3pt}u\hspace{-2.5pt}s}}}
\newcommand{\CQS}{\ensuremath{\mathcal{C\hspace{-4.2pt}p\hspace{-2.5pt}c\hspace{-2.2pt}t\hspace{-1.7pt}Q\hspace{-1.8pt}u\hspace{-3pt}a\hspace{-3pt}n\hspace{-2.8pt}t}}}

\DeclareMathOperator{\ima}{im}
\DeclareMathOperator{\coima}{coim}
\DeclareMathOperator{\coker}{coker}

\newcommand{\C}{\mathbb{C}}
\newcommand{\R}{\mathbb{R}}
\newcommand{\N}{\mathbb{N}}
\newcommand{\Z}{\mathbb{Z}}
\renewcommand{\emptyset}{\varnothing}

\DeclareRobustCommand{\SkipTocEntry}[5]{}

\numberwithin{equation}{section}

\begin{document}
\author{F.~D'Andrea}
\address[F.~D'Andrea]{Universit\`a di Napoli ``Federico II'' and I.N.F.N. Sezione di Napoli, Complesso MSA, Via Cintia, 80126 Napoli, Italy.}
\email{francesco.dandrea@unina.it}
\author{P.\ M.~Hajac}
\address[P.\ M. Hajac]{Instytut Matematyczny, Polska Akademia Nauk, ul.~\'Sniadeckich 8, Warszawa, 00-656 Poland}
\email{pmh@impan.pl}
\author{T.~Maszczyk}
\address[T. Maszczyk]{Instytut Matematyki, Uniwersytet Warszawski, ul.\ Banacha 2, 02-097 Warszawa, Poland}
\email{t.maszczyk@uw.edu.pl}
\author[B.~Zieli\'nski]{B.~Zieli\'nski\vspace{3mm}}
\address[B. Zieli\'nski]{Department of  Computer Science, University of \L{}\'od\'z, Pomorska 149/153 90-236
\L{}\'od\'z, Poland}
\email{bzielinski@uni.lodz.pl}
\title[Quantum CW-complexes]{\vspace*{-2cm}Quantum CW-complexes in\\[3pt] a Waldhausen category for unital C*-algebras}
\date{September 2026}

\begin{abstract}
Using the ring structure of the K-groups of finite CW-complexes, Atiyah and Todd unravelled the K-theory of complex projective spaces~$\CP^n$. Herein, in the
realm of noncommutative topology, we develop a new framework of finite quantum CW-complexes using the language of Waldhausen categories,
which allows us to enrich the class of standard morphisms between unital C*-algebras by adding inverses of *-homomorphisms that are isomorphisms in K-theory.
Our concept of quantum CW-complexes subsumes earlier constructions and enjoys a plethora of examples.
Moreover, the framework allows us to reduce problems concerning the multipushout quantum complex projective space $\CPq^n_H$ 
to the much more approachable setting
of the Vaksman--Soibelman  quantum complex projective space $\CPq^n_q$ enjoying the availability of graph-algebraic methods. In particular, these methods 
permit us to transport the ring structure from $K^0(\CP^n)$ to $K^0(\CPq_q^n)$. Finally, we adapt the formalism of Waldhausen categories 
to determine a natural set of free generators of $K^0(\CPq_H^n)$ from   a natural set of free generators of $K^0(\CPq_q^n)$, and to transport 
the ring structure from $K^0(\CPq_q^n)$ to~$K^0(\CPq_H^n)$.
\end{abstract}

\maketitle

{\parskip=1pt\tableofcontents}

\newpage
\section{Introduction}
\noindent
Our first key result, which is Theorem~\ref{hWald}, 
makes all unital C*-algebras part of our Waldhausen-category approach to K-theory. We modify the standard unpointed-Waldhausen-category setting
by adding a condition allowing us to formally invert weak equivalences, so as to obtain a homotopy category suitable to define quantum 
CW-complexes (Definition~\ref{def:finiteqKCWc}) and weak cellular equivalences (Definition~\ref{wceq}).
Then we claim a general K-theory result (Theorem~\ref{Kthm}) for all finite quantum CW-complexes, and a generalization of Atiyah's K\"unneth
theorem (Theorem~\ref{atiyahkun}) for finite strict quantum CW-complexes.

Recall that the classical CW-complex structure is a standard tool to compute topological K-theory. Better still, 
every topological space is weak homotopy equivalent to a CW-complex, and every continuous map between CW-complexes is homotopic to 
a cellular map. 
Although these results do not have known noncommutative counterparts, we hope that our approach to quantum 
CW-complex structures prompted by concrete examples of noncommutative C*-algebras contributes towards developing Noncommutative Topology.

\subsection{Motivation and strategy}
The multipushout quantum complex projective spaces $\CPq^n_H$  generalized the mirror quantum sphere~\cite{hms06}, 
and were motivated by the theory of  finite free distributive lattices~\cite{hkz12}. Computing their K-groups beyond $n=2$ (see~\cite{Rudnik})
turned out to be a challenging problem  solved in~\cite{hnpsz18}. Then, a natural set of free generators of $K_0(C(\CPq^2_H))$ 
was determined and analyzed in~\cite{d2025milnor}. Herein, not only we complete this work (Corollary~\ref{multmult}), but also go beyond by
showing that $\CPq^n_H$ enjoys a finite weak quantum CW-complex structure (Theorem~\ref{weakCWmult}), and 
that $K_0(C(\CPq^n_H))$ enjoys a natural ring structure (Corollary~\ref{multmult}).
We thus obtain the understanding of the K-theory of $\CPq^n_H$ on par with the understanding of the K-theory of $\CP^n$ obtained by Atiyah and Todd
in \cite{AT60} and by Adams in~\cite{a-jf62}.

Our strategy is to use Vaksman--Soibelman quantum complex projective space $\CPq_q^n$ (see~\cite{Mey95}) as a bridge junction
between $\CP^n$ and~$\CPq_H^n$.
The big technical advantage of $\CPq_q^n$ over $\CPq_H^n$ is that $C(\CPq_q^n)$ is a graph C*-algebra~\cite{HS02}, so its K-theory can be unravelled
using graph-theoretic methods~\cite{adht17,d2026k}. In particular, it allows us transfer the ring structure of $K^0(\CP^n)$ to $K_0(C(\CPq_q^n))$
in Theorem~\ref{thm:varphi}. On the other hand,
to transport K-theory from the Vaksman--Soibelman side to the multipushout side,
we prove the Tubular Neighbourhood Lemma (Lemma~\ref{lemma:3.3}), which is a pivotal technical result. Combining all this, we arrive at the aforementioned
Corollary~\ref{multmult}, which was a motivating concrete goal of this paper, and is a key application of new methods introduced herein.

\subsection{Background and state of the art}
Recall that the category of unital C*-algebras and unital *-homomorphisms is symmetric monoidal with respect to
 the maximal tensor product \cite[p.~76]{BN72}, 
which is also exact. 
Our choice of the maximal tensor product is also motivated by the fact that, for any C*-algebra $A$, the functor $B\mapsto K_0(A\otimes_{\max}B)$ is
half-exact and continuous, which is not true for the minimal tensor product. As customary, by the \emph{category of compact quantum spaces} 
we mean the opposite category of the category of unital C*-algebras and unital *-homomorphisms.
The maximal tensor product gives, by duality, a monoidal structure on the category compact quantum spaces, 
extending the Cartesian product of compact 
Hausdorff spaces. 

Some quantizations of CW-complexes (e.g., in the case of multipushout quantum complex  projective spaces) require fattening 
the image of an embedding of a lower-dimensional skeleton into a higher-dimensional skeleton  to its compact tubular neighbourhood.
In the classical setting, the compact tubular neighbourhood can be deformation-retracted onto the skeleton 
of lower dimension, so the fattening is inconsequential. 
However, in the noncommutative setting of our examples, where we implement fattening using the quantum disc, the situation gets more complicated
because the quantum disc is not contractible \cite[Corollary~2.7]{dhn}, unlike its classical counterpart.
To overcome this difficulty, we introduce a framework of
\emph{cw-Waldhausen categories} that allows a calculus of fractions \cite{gz12} leading to a notion of \emph{weak cofibration}.
Subsequently, we introduce the concept of a \emph{weak CW-complex structure}, which is equally effective as a tool for computing K-theory.
Finally,
notice that our construction generalizes Eilers--Loring--Pedersen's
notion of a noncommutative CW-complex \cite{eil-lor-ped-98} by allowing more general cells, gluing maps, and weak rather strict filtrations by skeleta. 

In the setting of unital C*-algebras, weak cofibrations are constructed by formally 
inverting K-equivalences, i.e.\ morphisms inducing isomorphisms of K-groups.
Observe also that, in the UCT class, the K-equivalence is a strengthening of the notion of a KK-equivalence. Indeed, 
every K-equivalence $A\to B$ between unital C*-algebras induces an isomorphism $K_*(A)\cong K_*(B)$ of graded abelian groups and, 
when $A$ is UCT and $K_*(B)$ is injective,
every such isomorphism comes from a unique invertible element in $KK_*(A,B)$ \cite[Theorem~2.1]{RS87}.
Furthermore,
to better understand the relation between K-equivalences and KK-equivalences, one can employ the Cuntz picture for KK-theory. In this picture, an 
element of $KK(A,B)$ is represented by a quasi-homomorphism between $A$ and $B$~\cite{j-c83}. On the other hand, for a K-equivalence we require 
the implementation of the K-theory isomorphism by a unital *-homomorphism, which refines the classification up to KK-equivalence.

For example, $\CP^n$ is KK-equivalent to the discrete space $P_{n+1}$ consisting of $n+1$ points because
they both have $K^0\cong\Z^{n+1}$ and $K^1\cong 0$, and every commutative C*-algebra is UCT.
On the other hand, there is no K-equivalence from
$\CP^n$ to $P_{n+1}$, nor from
$P_{n+1}$ to $\CP^n$,
 because a K-equivalence of commutative C*-algebras induces a ring isomorphism between the K-theories, 
 but $K^0(\CP^n)\cong\Z[x]/(x^{n+1})$ and $K^0(P_{n+1})\cong\Z^{n+1}$
are not isomorphic as rings.
Next,
a noncommutative example is given by the C*-algebras $\C$ and $M_n(\C)$, $n\geq 2$. Indeed, there exists no unital *-homomorphism 
$M_n(\C)\to\C$, and the unique unital *-homomorphism $\C\to M_n(\C)$ is not a K-equivalence. Therefore, $\C$ and $M_n(\C)$ are KK-equivalent 
but not K-equivalent.

Finally, let us note that group quotients $G/K$ are KK-equivalent to their quantum-group deformation quotient ${\sf G}_q/{\sf K}_q$ by~\cite{nt12},
which includes the Vaksman--Soibelman quantum complex projective spaces as special cases. Also, the multipushout quantum complex 
projective spaces are KK-equivalent to complex projective spaces because their K-groups are isomorphic and their C*-algebras are UCT.
Indeed, they are groupoid C*-algebras of locally compact Hausdorff amenable \'etale groupoids \cite{s-aj}, and by \cite[Theorem~10.7]{t-jl99} such groupoid C*-
algebras are UCT.
Hence,  all three types of projective spaces are KK-equivalent. Better yet, as a main point of this paper, we will show that the Vaksman--Soibelman
 and the multipushout quantum complex projective spaces are not only KK-equivalent but also K-equivalent.

\subsection{Notation and conventions}
Throughout the paper, we will denote by $\mathbb{R}$, $\mathbb{C}$, $\mathbb{Z}$, and $\mathbb{N}=\{0,1,2,\ldots\}$ the sets of real, 
complex, integer and natural numbers, 
respectively. 
We shall use $\otimes$ to stand for the maximal tensor product of C*-algebras. If $A$ is a C*-algebra and $n$ a positive integer, by $A^{\otimes n}$ we mean the 
maximal tensor product of $n$ copies of $A$. By convention, we set $A^{\otimes 0}:=\C$, the empty sum equal to~$0$, and the
empty product equal to~$1$.

Compact Hausdorff spaces will be denoted by ${X},{Y},{Z}$, etc. We will use the standard notation ${S}^n$ for  the $n$-dimensional sphere,
${B}^n$ for the $n$-dimensional closed ball, 
${D}={B}^2$ for the closed disc,
$\CP^n$ for the $n$-dimensional complex projective space, and $\RP^n$ for the $n$-dimensional real projective space.
If ${X}$ is a compact Hausdorff space and $n$ a positive integer,  ${X}^{\times n}$ will stand for
the Cartesian product of $n$ copies of ${X}$.
Finally, we will write the unital C*-algebra of complex-valued continuous functions on $X$ as~$C({X})$. Our notation
for the corresponding $\Z/2\Z$-graded K-theory group will be
$K^*({X})\cong K_*(C({X}))$.
 
We consider the category opposite to the category of unital C*-algebras and unital \linebreak *-homomorphisms as the category of \emph{compact quantum spaces},
 and denote it by~\CQS. 
Its objects will be written using a different font to avoid confusion (e.g., $\mathsf{X},\mathsf{Y},\mathsf{Z}$). We will use the notation
$\mathsf{S}^{2n+1}_q$ for the $(2n+1)$-dimensional Vaksman--Soibelman sphere \cite{VS91},
$\mathsf{S}^{2n+1}_H$ for the \mbox{$(2n+1)$-}dimensional Heegaard quantum sphere \cite{hnpsz18},
 $\mathsf{B}^{2n}_q$ for the $2n$-dimensional Hong--Szyma\'nski ball \cite{HS02}, 
$\mathsf{D}_q=\mathsf{B}^2_q$ for the quantum disc~\cite{c-la69,klimek1993},
$\CPq_q^n$ for  the $n$-dimensional Vaksman--Soibelman quantum complex projective spaces \cite{VS91,Mey95},
and $\CPq^n_H$ for the $n$-dimen\-sional multipushout quantum complex projective spaces \cite{hkz12}.
We will write the unital C*-algebra defining a compact quantum space $\mathsf X$ as~$C(\mathsf{X})$. Using the maximal tensor product, we will
define the Cartesian product of compact quantum spaces via
\begin{equation}\label{eq:monprod}
C(\mathsf X\times \mathsf Y):=C(\mathsf X)\otimes C(\mathsf Y) .
\end{equation}
Much as in the classical setting, we will write   $C(\mathsf X^{\times n}):=C(\mathsf X)^{\otimes n}$. 
In particular, we will use the notation $\mathsf{D}_q^{\times n}$, and refer to it as a \emph{quantum polydisc}.
Finally, we will employ the following notation for the $\Z/2\Z$-graded K-theory group:
$K^*(\mathsf{X}):=K_*(C(\mathsf{X}))$.

All the above-mentioned examples of quantum spheres are described by \mbox{${U}(1)$-C*-algebras}. Here, we call $A$ a \emph{${U}(1)$-C*-algebra} if there is a 
group homomorphism ${U}(1)\to \mathrm{Aut}(A)$ that is continuous with respect to the point-norm topology on $\mathrm{Aut}(A)$. If $A=C(X)$. 
Every such an action is 
tantamount to a continuous action of ${U}(1)$ on ${X}$ (e.g., see~\cite{Wil07}). A multiple tensor product of $U(1)$-C*-algebras
carries many ${U}(1)$-actions. We will specify the action on such a tensor product by putting a dot as a subscript on the factors where ${U}(1)$ acts non-trivially. 
For example, by $A_\bullet\otimes B\otimes C_\bullet$
we mean that the action of $u\in{U}(1)$ on $a\otimes b\otimes c$
is given by $\alpha_u(a)\otimes b\otimes\gamma_u(c)$, where $\alpha$ is the action on $A$ and $\gamma$ the action on~$C$.
The $U(1)$-action on the multiple tensor product of $U(1)$-algebras given by the tensor product of all involved $U(1)$-actions is called \emph{diagonal}. 

To end with, let us recall that there is a ${U}(1)$-equivariant 
isomorphism of ${U}(1)$-C*-algebras \cite[Section~2.3]{hnpsz18}
\begin{equation}\label{eq:gaugetrick}
\begin{split}
\phi:A_\bullet\otimes C({S}^1)_\bullet\otimes B_\bullet &\longrightarrow A\otimes C({S}^1)_\bullet\otimes B \,, \\
a\otimes f\otimes b &\longmapsto a_{(0)}\otimes a_{(1)}fb_{(1)}\otimes b_{(0)} \,, \\
A\otimes C(S^1)\ni a_{(0)}\otimes a_{(1)}&:=\alpha(\cdot)(a)\in C(U(1),A)\,,\\
C(S^1)\otimes B\ni b_{(1)}\otimes b_{(0)}&:=\beta(\cdot)(b)\in C(U(1),B)\,.
\end{split}
\end{equation}
Here $\alpha$ and $\beta$ are $U(1)$-actions on $A$ and $B$, respectively, and  both identifications 
$A\otimes C(S^1)=C(U(1),A)$ and $C(S^1)\otimes B=C(U(1),B)$
with the algebras of norm-continuous functions are the obvious ones.

\section{Cofibration-weakening Waldhausen categories}\label{sec:generaltheory}

\noindent
The customary framework for homotopy theory of CW-complexes uses model category structures.  
The question of the existence of a model category structure for C*-algebras has been raised by Schochet \cite{Schoch12} in the context of homotopy 
equivalences. However, in our approach, we invert also maps inducing isomorphisms in K-theory, which we call \emph{K-equivalences}.
We thus extend the Mayer--Vietoris principle for K-theory to the context of a specific structure, which we call a cofibration-weakening-Waldhausen
structure,
or  \emph{cw-Waldhausen} structure in short. (For standard Waldhausen categories see \cite{waldhausen-78}.)
Then, we apply this framework to the category \CQS\ of compact quantum spaces. Since general noncommutative C*-algebras need not admit 
characters, we are forced to consider an unpointed version of the standard Waldhausen structure \cite{rap-steim-19} with possibly distinct initial and 
terminal objects. In the classical situation going from pointed to unpointed does not change much.

In what follows, we provide a categorical framework of \emph{cw-Waldhausen categories} allowing a weakening of the notion of a CW-complex 
structure where weak equivalences are K-equivalences.

\begin{df}\label{def:Waldcat}\leftmargini=3.5em
An \emph{unpointed Waldhausen category} $\mathcal{C}$ is a category  with an initial object~$\emptyset $, a terminal object $\star$, and two 
distinguished classes of maps, \textit{Cof} of  \emph{cofibrations}, depicted as $\cof$, and \textit{W\hspace{-0.5pt}eq} of \emph{weak equivalences}, 
depicted as $\weq$, such that:
\begin{itemize}
\myitem[(\textit{Iso})] \label{ax:iso}
Isomorphisms are both cofibrations and weak equivalences.

\myitem[(\textit{Com})] \label{ax:com}
Both classes \textit{Cof} and \textit{W\hspace{-0.5pt}eq} are closed under composition.

\myitem[(\textit{Cof})] \label{ax:cof}
For any object $X$,  the unique morphism $\emptyset\longhookrightarrow X$ is a cofibration.
Furthermore, if  $X\longrightarrow\widetilde{X}$ is any morphism and $X\cof Y$ is a cofibration, then so is the pushout \mbox{$\widetilde{X}\cof \widetilde{X}\sqcup_{X}Y$}.

\myitem[(\textit{W\hspace{-0.5pt}eq})] \label{ax:weq}
For any commutative diagram\medskip
\begin{center}
\begin{tikzpicture}[->,baseline={([yshift=5pt]current bounding box.south)}]

\node (Z) at (0,1.7) {$Z$};
\node (X) at (2,1.7) {$X$};
\node (Y) at (4,1.7) {$Y$};

\node (Zt) at (0,0) {$\widetilde{Z}$};
\node (Xt) at (2,0) {$\widetilde{X}$};
\node (Yt) at (4,0) {$\widetilde{Y}$};

\draw (X) -- (Z);
\draw (Xt) -- (Zt);
\draw[right hook->] (X) -- (Y);
\draw[right hook->] (Xt) -- (Yt);
\draw (X) -- node[font=\footnotesize,sloped,above=-1pt] {$\sim$} (Xt);
\draw (Y) -- node[font=\footnotesize,sloped,above=-1pt] {$\sim$} (Yt);
\draw (Z) -- node[font=\footnotesize,sloped,above=-1pt] {$\sim$} (Zt);

\end{tikzpicture},\medskip
\end{center}
where the vertical arrows are weak equivalences and right horizontal arrows are cofibrations, the induced map 
$Z\sqcup_{X}Y\weq\widetilde{Z}\sqcup_{\widetilde{X}}\widetilde{Y}$ is a weak equivalence.
\end{itemize}
\end{df}

We shall work with the homotopy category $\Ho (\mathcal{C}):=\mathcal{C}[\textit{W\hspace{-0.5pt}eq}^{-1}]$.
Since uncontrolled inverting of weak equivalences could lead, in principle, to unwanted collapses of the homotopy category, to prevent this from happening, we 
add the condition \ref{ax:cwW} to obtain the following.

\begin{df}\label{cw-Wdef}\leftmargini=4em
A \emph{cw-Waldhausen category} (cofibration-weakening-Waldhausen category) is an unpointed
Waldhausen category satisfying:
\begin{itemize}
\myitem[(\textit{cw-W})] \label{ax:cwW}
For every pushout diagram
\begin{center}
\begin{tikzpicture}[->,baseline={([yshift=-5pt]current bounding box.center)}]

\node (Z'') at (1.5,3) {$\,\doublewidetilde{Z}\!$};
\node (Y') at (0,1.5) {$\widetilde{Y}$};
\node (Z') at (3,1.5) {$\widetilde{Z}$};
\node (Y) at (1.5,0) {$Y$};

\path[->,font=\scriptsize,inner sep=1pt]
		(Y') edge[right hook->] node[above left] {$\widetilde{j}$} (Z'')
		(Y) edge node[above right] {$g$}  (Y')
		(Y) edge[right hook->] node[above left] {$j$} (Z')
		(Z') edge node[above right] {$\widetilde{h}$}   (Z'');
\end{tikzpicture},
\end{center}
where $j$ is a cofibration, $\widetilde{h}$ is a weak equivalence if and only if so is $g$.

\myitem[(\textit{2/3})] \label{ax:2outof3}
In the commutative diagram
\begin{center}
\begin{tikzpicture}[->]

\node (A) at (0,0) {$X$};
\node (B) at (2,0) {$Y$};
\node (C) at (4,0) {$Z$,};

\path[->] (A) edge (B) (B) edge (C) (A) edge[bend right] (C);

\end{tikzpicture}
\end{center}
if any two morphisms are weak equivalences, then so is the third one.

\end{itemize}
\end{df}

Note that $\widetilde{j}$ is automatically a cofibration because of \ref{ax:cof}.
Thanks to the condition \mbox{\ref{ax:cwW}},
we can work within $\Ho(\mathcal{C})$ with the left fraction calculus of the form \mbox{\textit{W\hspace{-0.5pt}eq}$^{-1}\circ$\textit{Cof}}.
We refer to the morphisms in the class \mbox{\textit{W\hspace{-0.5pt}eq}$^{-1}\circ$\textit{Cof}} as \emph{weak cofibrations}.
We can represent them as cospans
\begin{center}
\begin{tikzpicture}[scale=1.2,->,baseline={([yshift=5pt]current bounding box.south)}]

\node (Yt) at (1.2,1.3) {$\widetilde{Y}$};
\node (X) at (0,0) {$X$};
\node (Y) at (2.4,0) {$Y$};

\draw[right hook->] (X) -- (Yt);
\draw (Y) -- node[font=\footnotesize,sloped,above] {$\sim$} (Yt);

\end{tikzpicture} ,
\end{center}
which we can think of as generalized morphisms denoted by $X\rightarrowtail Y$.
Under the \mbox{\ref{ax:cwW}} assumption, their composition in the homotopy category $\Ho(\mathcal{C})$ can be defined as follows:
\begin{center}
\begin{tikzpicture}[scale=1.5,baseline={([yshift=5pt]current bounding box.south)}]

\node (Z') at (1,1) {$\widetilde{Z}$};
\node (Y'') at (4,1) {$\widetilde{Y}$};
\node (Y') at (0,0) {$Y$};
\node (Z) at (2,0) {$Z$};
\node (o) at (2.5, 0.5) {$\circ$};
\node (X) at (3,0) {$X$};
\node (Y) at (5,0) {$Y$};
\node (=) at (5.5, 0.5) {$=$};
\node (Z''') at (7,1) {$\doublewidetilde{Z}$};
\node (X') at (6,0) {$X$};
\node (Z'') at (8,0) {$Z$};

\path[->,font=\footnotesize]
	(Y') edge[right hook->] node[above left] {$j$} (Z')
	(Z) edge node[above right] {$h$} node[sloped,below] {$\sim$} (Z')
	(X) edge[right hook->] node[above left] {$i$} (Y'')
	(Y) edge node[above right] {$g$} node[sloped,below] {$\sim$} (Y'')
	(X') edge[right hook->] node[above left] {$\widetilde{j}\circ i$} (Z''')
	(Z'') edge node[above right]{$\widetilde{h}\circ h$} node[sloped,below] {$\sim$} (Z''');

\end{tikzpicture} .
\end{center}
Here $\widetilde{j}$ and $\widetilde{h}$ are the arrows completing the pushout square in the diagram
\begin{center}
\begin{tikzpicture}[scale=1.5,baseline={([yshift=5pt]current bounding box.south)}]

\node (Z'') at (2,3) {$\doublewidetilde{Z}$};
\node (Y') at (1,2) {$\widetilde{Y}$};
\node (Z') at (3,2) {$\widetilde{Z}$};
\node (X) at (0,1) {$X$};
\node (Y) at (2,1) {$Y$};
\node (Z) at (4,1) {$Z$};

\path[->,font=\footnotesize]
	(Y') edge [dashed, right hook->] node[above left]{$\widetilde{j}$} (Z'')
	(Z') edge[dashed] node[above right] {$\widetilde{h}$}  node[sloped,below] {$\sim$} (Z'')   
	(X) edge[right hook->] node[above left] {$i$} (Y')
	(Y) edge node[above right] {$g$} node[sloped,below] {$\sim$} (Y')
	(Y) edge[right hook->] node[above left] {$j$} (Z')
	(Z) edge node[above right] {$h$} node[sloped,below] {$\sim$} (Z');

\end{tikzpicture}.
\end{center}
Equivalently, we can write
\[
(h^{-1}\circ j)  \circ    (g^{-1}\circ i) = (\widetilde{h} \circ h)^{-1}\circ (\widetilde{j}\circ i) .
\]
Furthermore, note that the axiom \ref{ax:2outof3} is motivated by the fact that weak equivalences become isomorphisms in the homotopy category, and isomorphisms satisfy the 2-out-of-3 property.

We are interested in cw-Waldhausen categories with a symmetric monoidal structure which is, in some sense, compatible with the Waldhausen structure. 
Such a monoidal structure should be a weakening of the notion of the Cartesian product. The minimal requirements are as follows:
\begin{enumerate}[label={\arabic*)},leftmargin=1.6em]
\item The monoidal unit is the terminal object.
\item The monoidal product is compatible with cofibrations.
\item The weak equivalences are compatible with the monoidal product with objects in some distinguished subcategory.
\end{enumerate}
This leads to the next definitions.

\begin{df}\label{df:wcontr}
An object $X$ in a cw-Waldhausen category is called \emph{weakly contractible} if the (unique) map $X\to\star$ to the terminal object is a weak equivalence.
\end{df}

\pagebreak

\begin{rmk}\label{rmk:weak}
Assume that there exists a morphism $x:\star\to X$. 
Then it follows from the axiom \ref{ax:2outof3} applied to the diagram
\begin{center}
\begin{tikzpicture}[->,scale=0.8]

\node (A) at (0,0) {$\star$};
\node (B) at (2,0) {$X$};
\node (C) at (4,0) {$\star$\;,};

\path[->] (A) edge node[above] {$x$} (B) (B) edge (C) (A) edge[bend right] (C);

\end{tikzpicture}
\end{center}
and the fact that the unique morphism $\star\to\star$ is an identity, that
$X$ is weakly contractible if and only if $x$ is a weak equivalence.
\end{rmk}

\begin{df}\label{cw-sWdef}\leftmargini=5.5em
A \emph{cw-Waldhausen symmetric monoidal category} $\mathcal{C}$ is both a cw-Waldhausen and a symmetric monoidal category, with the monoidal product denoted by $\times$ and with the terminal object as a monoidal unit, together with a distinguished full subcategory~$\mathcal{C}^{\times}$, such that:
\begin{itemize}\itemsep=2pt
\myitem[(\textit{Con-Weq})] \label{ax:contr}
Any morphism $X\to Y$ between weakly contractible objects is a weak equivalence.

\myitem[(\textit{Cof-Mon})] \label{ax:cofmon}
For every cofibration $f: X\cof Y$ and every object $Z$ in $\mathcal{C}$, the induced morphisms $Z\times f: Z\times X\cof Z\times Y$ and $f\times Z : X\times Z\cof Y\times Z$ are cofibrations.

\myitem[(\textit{Sub})] \label{ax:sub}
For every pushout diagram\vspace{-5pt}
\[
\begin{tikzpicture}[scale=1.2,baseline={([yshift=-3pt]current bounding box.center)}]

\node (Z'') at (1,2) {$W$};
\node (Y') at (0,1) {$X$};
\node (Z') at (2,1) {$Y$};
\node (Y) at (1,0) {$Z$};

\path[->]
		(Y') edge[right hook->] (Z'')
		(Y) edge (Y')
		(Y) edge[right hook->] (Z')
		(Z') edge (Z'');
\end{tikzpicture},
\]
in $\mathcal{C}$, if $X,Y,Z$ are objects in $\mathcal{C}^{\times}$, then
$W$ is in $\mathcal{C}^{\times}$.

\myitem[(\textit{W\hspace{-0.5pt}eq-Mon})] \label{ax:weqmon}
The distinguished subcategory $\mathcal{C}^{\times}$ is a monoidal subcategory of $\mathcal{C}$ such that, for every weak equivalence $f:X\weq Y$ and every object $S$ in $\mathcal{C}^{\times}$, the induced morphisms
$S\times f: S\times X\weq S\times Y$ and $f\times S: X\times S\weq Y\times S$ are weak equivalences.

\end{itemize}
\end{df}

\begin{df}
Let  $\mathcal{C}$ be an unpointed Waldhausen category and
\begin{equation}\label{eq:wpo}
\begin{tikzpicture}[scale=1,baseline=(current bounding box.center)]

\node (X) at (1.5,3) {$\widetilde{W}$};
\node (Z'') at (1.5,1.5) {$W$};
\node (Y') at (0,1.5) {$X$};
\node (Z') at (3,1.5) {$Y$};
\node (Y) at (1.5,0) {$Z$};

\path[->,font=\scriptsize,inner sep=2pt]
		(Y) edge (Y')
		(Y) edge (Z')
		(Z') edge (Z'')
		(Z'') edge node[right] {$\psi$} (X)
		(Y') edge (X)
		(Z') edge (X);
\end{tikzpicture}
\end{equation}
be a commutative diagram in $\mathcal{C}$ such that
such that $\psi$ is a weak equivalence and the outer square is a pushout. The induced square
\begin{center}
\begin{tikzpicture}[scale=1.3]

\node (Z'') at (1,2) {$W$};
\node (Y') at (0,1) {$X$};
\node (Z') at (2,1) {$Y$};
\node (Y) at (1,0) {$Z$};

\path[->,font=\scriptsize,inner sep=2pt]
		(Y') edge (Z'')
		(Y) edge (Y')
		(Y) edge (Z')
		(Z') edge (Z'');
\end{tikzpicture}
\end{center}
in $\Ho(\mathcal{C})$ is called a \emph{weak pushout}.
(Note that, if the morphism $Z\to Y$ is a cofibration, then the morphism $X\to W$ is a weak cofibration
by \ref{ax:cof}.)
\end{df}

\subsection{The  cw-Waldhausen category of compact quantum spaces}
Recall that the pullback  \cite[Definition~15.3.1]{b-b98} of
two *-homorphisms
\mbox{$A_1\xrightarrow{\,\pi_1\;}A_{12}\xleftarrow{\;\pi_2\,}A_2$}
can be realized as
\begin{equation}\label{eq:canpullback}
A_1\times_{A_{12}}A_2:=\big\{(a_1,a_2)\in A_1\times A_2\;\big|\;\pi_1(a_1)=\pi_2(a_2)\big\}
\end{equation}
with *-homomorphisms $A_1\xleftarrow{}A_1\times_{A_{12}}A_2\xrightarrow{}A_2$ given by the projections on the two factors. 
Pullbacks in the category of unital C*-algebras dualize to pushouts in the opposite category \CQS.
The category of unital C*-algebras has a terminal object given by the zero C*-algebra, which defines, by duality, the initial object $\emptyset$ of \CQS.

As a symmetric monoidal structure on C*-algebras we consider the maximal tensor product \cite[p.~76]{BN72} that gives, by duality, the monoidal product \eqref{eq:monprod} in \CQS. The initial object $\C$ in the category of C*-algebras dualizes to the terminal object $\star$ in \CQS\ (the one-point space), which is a monoidal unit. Note that in \CQS\ the categorical product is dual to the free product of C*-algebras, so it is different from the above monoidal product above. However, in the subcategory of compact Hausdorff spaces these two notions coincide.

A unital *-homomorphism between unital C*-algebras that induces an isomorphism in K-theory will be called a \emph{K-equivalence}.
In the category \CQS, we define \emph{cofibrations} as the opposites of surjective unital \mbox{*-homomorphisms} and \emph{weak equivalences} as the opposites of K-equivalences. We shall use the symbol $\weq$ to denote both a K-equivalence $C(\mathsf Y)\weq C(\mathsf X)$ and its opposite morphism $\mathsf X\weq\mathsf Y$.

Now, we turn our attention to the K{\"u}nneth formula, which is always satisfied by commutative unital C*-algebras (compact Hausdorff spaces).
For the maximal tensor product of C*-algebras, the relevant class $\mathcal{N}_{\max}$ is defined in \cite{Uuy11} as follows. A C*-algebra $A$ belongs to 
$\mathcal{N}_{\max}$ if and only if $K_*(A\otimes B)=0$ for every C*-algebra $B$ with $K_*(B)=0$.
A C*-algebra satisfies the K{\"u}nneth theorem with respect to the maximal tensor product if and only if it belongs to $\mathcal{N}_{\max}$ \cite[Theorem 4.1]{Uuy11}.
Therefore, we will call $\mathcal{N}_{\max}$ the \emph{K{\"u}nneth class}.

In the context of \CQS, we consider the full subcategory whose objects are dual to unital C*-algebras in the K{\"u}nneth class $\mathcal{N}_{\max}$. This subcategory, denoted by $\CQS^{\times}$ in what follows, will play the role of the distinguished subcategory $\mathcal{C}^{\times}$ in Definition \ref{cw-sWdef}.

\begin{thm}\label{hWald}
With the structure and the distinguished subcategory $\CQS^{\times}$ introduced above, \CQS\ is a cw-Waldhausen symmetric monoidal category.
\end{thm}

\begin{proof}
Both axioms \ref{ax:iso} and \ref{ax:com} are clearly satisfied.

\ref{ax:cof} A *-homomorphism into the zero algebra is clearly surjective.
Moreover, if $B\twoheadrightarrow A$ is a surjective  unital *-homomorphism and $\widetilde{A}\to A$  any *-homomorphism then the pullback  *-homomorphism  $\widetilde{A}\times_{A} B\twoheadrightarrow \widetilde{A}$ is unital surjective.

\ref{ax:weq}
Reversing all the arrows, \cite[Theorem 2.4]{d2025milnor} becomes exactly the axiom \ref{ax:weq}. 

\ref{ax:cwW}
Assume that in the following pullback diagram of C*-algebras
\begin{center}
\begin{tikzpicture}[scale=1.3]

\node (BB) at (1,3) {$C$};
\node (A) at (0,2) {$A$};
\node (B) at (2,2) {$B$};
\node (D) at (1,1) {$D$};

\path[-To,font=\scriptsize]
		(BB) edge node[above left] {$\pi$} (A)
		(BB) edge node[above right] {$\beta$}  (B)
		(A) edge node[below left] {$\delta$} (D)
		(B) edge[->>] node[below right] {$\tau$} (D);
\end{tikzpicture}
\end{center}
the *-homomorphism $\tau$ is surjective. Then so is $\pi$, and $\beta$ is a K-equivalence if and only if $\delta$ is a K-equivalence.
Indeed, since surjective *-homomorphisms are regular epimorphisms, they are stable under all pullbacks, which proves the surjectivity of $\pi$. Next, thanks to the surjectivity of $\tau$, we have the Mayer--Vietoris six-term exact sequence \cite{Hilgert1986}
\begin{center}
\begin{tikzpicture}[-To,baseline={([yshift=5pt]current bounding box.south)}]

\node (a1) at (0,1.7) {$K_0(C)$};
\node (a2) at (4,1.7) {$K_0(A)\oplus K_0(B)$};
\node (a3) at (8,1.7) {$K_0(D)$};
\node (b1) at (0,0) {$K_1(D)$};
\node (b2) at (4,0) {$K_1(A)\oplus K_1(B)$};
\node (b3) at (8,0) {$K_1(C)$};

\begin{scope}[font=\scriptsize]
\draw (a1) -- node[above,scale=1.2] {$\binom{\pi_{*}}{\beta_{*}}$} (a2);
\draw (a2) -- node[above] {$(\delta_{*}, -\tau_{*})$} (a3);
\draw (a3) -- (b3);
\draw (b3) -- node[above,scale=1.2] {$\binom{\pi_{*}}{\beta_{*}}$} (b2);
\draw (b2) -- node[above] {$(\delta_{*}, -\tau_{*})$} (b1);
\draw (b1) -- (a1);
\end{scope}

\end{tikzpicture},
\end{center}
which we use as follows. 
First, if $\delta$ is a K-equivalence, we have a section of the map $(\delta_{*},-\tau_{*})$, provided by composing the standard embedding \mbox{$i_{A} : K_*(A)\rightarrow K_*(A)\oplus K_*(B)$} with $\delta_{*}^{-1}$. This cuts the six-term exact sequence into split short exact sequences
\begin{center}
\begin{tikzpicture}[-To]

\node (a1) at (0.5,0) {$0$};
\node (a2) at (3,0) {$K_*(C)$};
\node (a3) at (7,0) {$K_*(A)\oplus K_*(B)$};
\node (a4) at (11,0) {$K_*(D)$};
\node (a5) at (13.5,0) {$0$\ .};

\begin{scope}[font=\scriptsize]
\draw (a1) -- (a2);
\draw (a2) -- node[above,scale=1.2] {$\binom{\pi_{*}}{\beta_{*}}$} (a3);
\draw (a3) -- node[above] {$(\delta_{*}, -\tau_{*})$} (a4);
\draw (a4) -- (a5);
\draw[dashed] (a4) edge[bend left=30] node[below] {$i_{A} \circ \delta_{*}^{-1}$} (a3);
\end{scope}

\end{tikzpicture}
\end{center}
The above splittings produce idempotent endomorphisms of $K_*(A)\oplus K_*(B)$ 
\[
p:=i_{A} \circ \delta_{*}^{-1}\circ (\delta_{*},-\tau_{*})=p\circ p,\quad
p^{\perp}:= \id - p ,
\]
such that 
\[
\ker (p) = \ker (\delta_{*},-\tau_{*}) = \ima \binom{\pi_*}{\beta_*}, \quad \coima(p^{\perp}) = \coima (p_{B}) .
\]
Here, $p_{B}$ is the canonical projection from $K_*(A)\oplus K_*(B)$ onto $K_*(B)$.
We thus arrived at the commutative diagram
\begin{center}
\begin{tikzpicture}[scale=1.5,baseline={([yshift=5pt]current bounding box.south)}]

\node (P) at (0,0) {$\ker (p)$};
\node (C) at (0,2) {$K_{*}(C)$};
\node (AB) at (3,1) {$K_{*}(A)\oplus K_{*}(B)$};
\node (PP) at (6,0) {$\coima (p^{\perp})$};
\node (B) at (6,2) {$K_{*}(B)$};

\path[-To]
     (C) edge  node[left] 
     [scale=0.75]{$\cong$} (P)
     (C) edge  node[below left][scale=1.00]{$\binom{\pi_{*}}{\beta_{*}}$} (AB)   
		(C) edge node[above] [scale=0.75]{$\beta_{*}$} (B)
		(P) edge (AB)
		(P) edge node[above] [scale=0.75]{$\cong$} (PP)
		(AB) edge  node[below right][scale=0.75]{$p_{B}$} (B)
		(AB) edge (PP)
		(PP) edge  node[right] 
     [scale=0.75]{$\cong$} (B);
		
		\end{tikzpicture}
\end{center}
which proves that
$\beta_{*}$ is an isomorphism, i.e.~$\beta$ is a  K-equivalence.

Now, to prove the opposite implication, assume in turn that $\beta$ is a  K-equivalence. This implies that we have a map, given by the composition of $\beta_{*}^{-1}$ with the standard projection $p_{B} : K_*(A)\oplus K_*(B)\rightarrow K_*(B)$,  which is a retraction of the map $\binom{\pi_{*}}{\beta_{*}}$ that cuts the six-term exact sequence into split short exact sequences
\begin{center}
\begin{tikzpicture}[-To]

\node (a1) at (0.5,0) {$0$};
\node (a2) at (3,0) {$K_*(C)$};
\node (a3) at (7,0) {$K_*(A)\oplus K_*(B)$};
\node (a4) at (11,0) {$K_*(D)$};
\node (a5) at (13.5,0) {$0$\ .};

\begin{scope}[font=\scriptsize]
\draw (a1) -- (a2);
\draw (a2) -- node[above,scale=1.2] {$\binom{\pi_{*}}{\beta_{*}}$} (a3);
\draw (a3) -- node[above] {$(\delta_{*}, -\tau_{*})$} (a4);
\draw (a4) -- (a5);
\draw[dashed] (a3) edge[bend left=30] node[below] {$\beta_{*}^{-1}\circ p_B$} (a2);
\end{scope}

\end{tikzpicture}
\end{center}
The above retraction yields idempotent endomorphisms of $K_*(A)\oplus K_*(B)$
$$
q:=\binom{\pi_{*}}{\beta_{*}} \circ \beta_{*}^{-1}\circ p_{B}=q\circ q,\quad q^{\perp}:=\id - q ,
$$
such that 
\[
\ima q^{\perp} = \ima i_{A},\quad  \coker q =\coker \binom{\pi_{*}}{\beta_{*}}=\coima \left(\delta_{*},-\tau_{*}\right)  .
\]
We thus arrived at the commutative diagram
\begin{center}
\begin{tikzpicture}[scale=1.5]

\node (QQ) at (0,0) {$\ima q^{\perp}$};
\node (A) at (0,2) {$K_{*}(A)$};
\node (AB) at (3,1) {$K_{*}(A)\oplus K_{*}(B)$};
\node (Q) at (6,0) {$\coker q$};
\node (D) at (6,2) {$K_{*}(D)$};

\path[-To]
     (A) edge  node[left] 
     [scale=0.75]{$\cong$} (P)
     (A) edge  node[below left][scale=0.75]{$i_{A}$} (AB)   
		(A) edge node[above] [scale=0.75]{$\delta_{*}$} (D)
		(QQ) edge (AB)
		(QQ) edge node[above] [scale=0.75]{$\cong$} (Q)
		(AB) edge  node[below right][scale=0.75]{$\left(\delta_{*},-\tau_{*}\right)$} (D)
		(AB) edge (Q)
		(Q) edge  node[right] 
     [scale=0.75]{$\cong$} (D);
		
		\end{tikzpicture}
\end{center}
which proves that $\delta_{*}$ is an isomorphism,  i.e.~ $\delta$ is a K-equivalence.

\ref{ax:2outof3}
Consider the commutative diagram of unital *-homomorphisms of unital C*-algebras
\begin{center}
\begin{tikzpicture}

\node (A) at (0,0) {$A$};
\node (B) at (3,0) {$B$};
\node (C) at (6,0) {$C$};

\path[->] (A) edge node[above] {$f$} (B) (B) edge node[above] {$g$} (C) (A) edge[bend right] node[below] {$g\circ f$} (C);

\end{tikzpicture}
\end{center}
and the induced diagram in K-theory
\begin{center}
\begin{tikzpicture}

\node (A) at (0,0) {$K_*(A)$};
\node (B) at (3,0) {$K_*(B)$};
\node (C) at (6,0) {$K_*(C)$\;.};

\path[->] (A) edge node[above] {$K(f)$} (B) (B) edge node[above] {$K(g)$} (C) (A) edge[bend right] node[below] {$K(g)\circ K(f)$} (C);

\end{tikzpicture}
\end{center}
Since isomorphisms of graded abelian groups satisfy the 2-out-of-3 property, if two out of three maps $f$, $g$, and $g\circ f$ are K-equivalences, then so is the third one.

\ref{ax:contr}
Suppose that $\mathsf X$ is K-contractible. This means that the morphism \mbox{$\mathsf X\to\star$} is a weak equivalence, i.e.~the dual unital *-homomorphism $\C\to A:=C(\mathsf X)$ is a K-equivalence. 
Next, let $B:=C(\mathsf Y)$, where $\mathsf Y$ is also K-contractible. Then, if $\varphi:B\to A$ is any unital *-homomorphism, $\varphi$ is clearly a K-equivalence because $K_*(A)$ and $K_*(B)$ are generated by the class of the unit.

\ref{ax:cofmon}
From the right exactness of the maximal tensor product, we deduce that, if $\phi:B\twoheadrightarrow A$ is a surjective unital *-homomorphism and $C$ is any unital C*-algebra, then $C\otimes \phi:C\otimes B\to C\otimes A$ 
and $\phi\otimes C:B\otimes C\to A\otimes C$
are surjective as well. Setting $A:=C(\mathsf X)$, $B:=C(\mathsf Y)$ and $C:=C(\mathsf Z)$, and passing to the opposite category, we get the axiom \ref{ax:cofmon}.

\pagebreak

\ref{ax:sub}
Consider a pullback diagram of unital C*-algebras and unital *-homomorphisms
\begin{center}
\begin{tikzpicture}[scale=1.2]

\node (BB) at (1,3) {$P$};
\node (A) at (0,2) {$A$};
\node (B) at (2,2) {$B$};
\node (D) at (1,1) {$C$};

\path[-To,font=\scriptsize]
		(BB) edge (A)
		(BB) edge (B)
		(A) edge (D)
		(B) edge[->>] (D);
\end{tikzpicture}
\end{center}
with $A,B,C$ in $\mathcal{N}_{\max}$, and take $D$ to be any C*-algebra such that $K_*(D)=0$. Then
\[
K_*(A\otimes D)=K_*(B\otimes D)=K_*(C\otimes D)=0
\]
by the definition of $\mathcal{N}_{\max}$.
Now, from \cite[Remark 3.10]{Ped}
it follows that
\begin{center}
\begin{tikzpicture}[scale=1.3]

\node (BB) at (1,3) {$P\otimes D$};
\node (A) at (0,2) {$A\otimes D$};
\node (B) at (2,2) {$B\otimes D$};
\node (D) at (1,1) {$C\otimes D$};

\path[-To,font=\scriptsize]
		(BB) edge (A)
		(BB) edge (B)
		(A) edge (D)
		(B) edge[->>] (D);
\end{tikzpicture}
\end{center}
is also a pullback diagram. Finally, from the associated Mayer--Vietoris six-term exact sequence, we deduce that $K_*(P\otimes D)=0$, which proves that $P\in\mathcal{N}_{\max}$. By passing to the opposite category, we get the axiom \ref{ax:sub}.

\ref{ax:weqmon}
Since $\C\otimes B\cong B$ for all C*-algebras $B$, we infer that $\C\in\mathcal{N}_{\max}$, and so $\CQS^{\times}$ contains the monoidal unit. Moreover, $A\otimes B\in\mathcal{N}_{\max}$ whenever $A,B\in\mathcal{N}_{\max}$ (see \cite[Lemma~4.4(3)]{Uuy11}), so $\CQS^{\times}$ is a monoidal subcategory of \CQS.
Next, given a K-equivalence $\varphi: B\rightarrow A$, consider the following diagram comparing K\"{u}nneth short exact sequences that hold for every $C$ in $\mathcal{N}_{\max}$
\cite[Theorem~4.1]{Uuy11}:
\begin{center}
\begin{tikzpicture}[-To,baseline={([yshift=5pt]current bounding box.south)}]

\node (a1) at (0,2) {$0$};
\node (a2) at (3,2) {$K_*(C) \otimes K_*(A)$};
\node (a3) at (7,2) {$K_*(C\otimes A)$};
\node (a4) at (11,2) {${Tor}(K_*(C), K_*(A))$};
\node (a5) at (14,2) {$0$};

\node (b1) at (0,0) {$0$};
\node (b2) at (3,0) {$K_*(C) \otimes K_*(B)$};
\node (b3) at (7,0) {$K_*(C\otimes B)$};
\node (b4) at (11,0) {${Tor}(K_*(C), K_*(B))$};
\node (b5) at (14,0) {$0$};

\draw (a1) -- (a2);
\draw (a2) -- (a3);
\draw (a3) -- (a4);
\draw (a4) -- (a5);
\draw (b1) -- (b2);
\draw (b2) -- (b3);
\draw (b3) -- (b4);
\draw (b4) -- (b5);
\draw (b2) -- node[left,font=\scriptsize] {$K_*(C) \otimes K_*(\varphi)$} (a2);
\draw (b3) -- node[fill=white,font=\scriptsize] {$K_*(C\otimes\varphi)$}  (a3);
\draw (b4) -- node[right,font=\scriptsize] {${Tor}(K_*(C), K_*(\varphi))$} (a4);

\end{tikzpicture}.
\end{center}
The first and third vertical arrows are isomorphisms because we are applying two functors to the isomorphism $K_*(\varphi)$. Now it follows from the exactness of the rows and the Five Lemma that the middle vertical arrow is also an isomorphism, which means that the \mbox{*-homomorphism} $C\otimes \varphi: C\otimes B\rightarrow C\otimes A$ is a K-equivalence.
\end{proof}

From now on, we will consider \CQS\ as equipped with the above cw-Waldhausen symmetric monoidal structure. 
In order to avoid confusion with classical weak contractibility as weak homotopy equivalence to a point, in the context of compact quantum spaces, we will refer 
to weak contractibility (in the sense of Definition~\ref{df:wcontr}) as \emph{K-contractibility}.

A widely used class of C*-algebras in the context of the K{\"u}nneth theorem is the \emph{bootstrap class} $\mathcal{N}$ \cite[Definition 22.3.4]{b-b98}. 
Since C*-algebras in $\mathcal{N}$ satisfy the K{\"u}nneth theorem \cite[Theorem 23.1.3]{b-b98}, and the minimal and maximal tensor products coincide 
for them, it follows from \cite[Lemma~4.4]{Uuy11} that $\mathcal{N}\subseteq\mathcal{N}_{\max}$.
Furthermore, any separable commutative C*-algebra belongs to $\mathcal{N}$ (see 22.3.5(d) in \cite{b-b98}), and hence to $\mathcal{N}_{\max}$. 
Therefore, restricting $\mathcal{N}$ to commutative unital C*-algebras yields the class $\mathcal{N}^c$ of all separable commutative unital C*-algebras
(compact metrizable spaces),
which is contained in the K{\"u}nneth class $\mathcal{N}_{\max}$ restricted to commutative C*-algebras. However, as the distinguished full subcategory 
$\mathcal{C}^{\times}$ of $\CHS$, we choose that internally defined  in $\CHS$ much in the same way as we defined 
$\mathcal{N}_{\max}$.
 Here, the objects of $\CHS^{\times}$ are all compact Hausdorff spaces $X$ with the property 
 that $K^*(Y)=0$ for a locally compact Hausdorff space $Y$
implies that $K^*(X\times Y)=0$. Clearly, all the duals of the aforementioned intersection of $\mathcal{N}_{\max}$  with commutative unital 
C*-algebras belong to~$\CHS^{\times}$.

Recall that the monoidal product in $\CHS$ is given by the Cartesian product with product topology,
the initial object is the empty space,
the terminal object is the one-point space,
cofibrations are closed embeddings,
and we choose weak equivalences to be K-equivalences.
Now, regarding $\CHS$ as a full monoidal subcategory of $\CQS$, we obtain the following comparison theorem by restricting the reasoning in Theorem \ref{hWald}
to commutative C*-algebras and then reversing arrows.

\begin{prop}\label{cor:cHcwW}
With the structure and the distinguished subcategory $\CHS^{\times}$ introduced above, \CHS\ is a cw-Waldhausen symmetric monoidal category.
Furthermore, the Gelfand--Naimark full embedding $\varphi: \CHS\hookrightarrow \CQS$
\begin{enumerate}[label={\arabic*)},leftmargin=1.6em]\item is strongly monoidal, 
\item preserves cofibrations and weak equivalences, 
\item satisfies the condition
$\varphi^{-1}(\CQS^{\times})\subseteq \CHS^{\times}$.
\end{enumerate}
\end{prop}

\begin{rmk}\label{rem:trivialK}
Let $\mathsf X$ be a compact quantum space and $A:=C(\mathsf X)$.
With respect to the above cw-Waldhausen symmetric monoidal structure,
$\mathsf X$ is K-contractible if and only if $K_0(A)$ is generated by the class $[1_A]$ and $K_1(A)=0$.
On the other hand, recall that $\mathsf X$ is called \emph{contractible} if there exists a character $x:A\to\C$
and a *-homomorphism
$h:A\to C([0,1],A)$
to the cylinder of $A$ such that ${ev}_0\circ h=\eta\circ x$ and  ${ev}_1\circ h=\id_A$, where $\eta:\C\to A$ is the unit map. This is an obvious generalization of the notion of contractibility of a compact Hausdorff space. It is also not difficult to check that this is equivalent to \cite[Definition~2.6]{dhn}. Finally,
by the homotopy invariance of K-theory, if $\mathsf X$ is contractible, then it is also K-contractible (cf.~the comment after \cite[Definition~2.6]{dhn}).
\end{rmk}

We end this section with a series of examples illustrating our cw-Waldhausen formalism at work. We begin with an example concerning contractibility, and continue with  examples and non-examples of weak equivalences in \CQS.

\begin{example}\label{ex:trivialK}
The quantum disk $\mathsf D_q$ dual to the Toeplitz C*-algebra $\mathcal{T}$ is not contractible (see below Corollary~2.7 in \cite{dhn}). The same argument as in \cite{dhn} proves that the polydisk $\mathsf D_q^{\times n}$ is also not contractible for any $n\geq 1$.
On the other hand, it follows from the K{\"u}nneth formula that $\mathsf D_q^{\times n}$
is K-contractible for all $n\geq 1$.
Better still, quantum balls $\mathsf B_q^{2n}$ are also K-contractible for all $n\geq 1$ \cite[Example~6.5]{HS02}. They are K-equivalent to quantum polydisks (see Section~\ref{sec:avb2}) but they are not isomorphic to them (see Proposition~\ref{prop:noniso}).
\end{example}

\begin{example}\label{ex:toeplitz}
Let $A$ be a unital C*-algebra with an action $\alpha$ of $\Z$. The crossed product $A\rtimes_\alpha \Z$ is the universal C*-algebra generated by $A$ and a unitary $u$ satisfying $uau^*=\alpha_1(a)$ for all $a\in A$. Next, let $\mathcal{T}$ be the Toeplitz algebra, whose generating unilateral shift we denote by $s$.
The unital C*-subalgebra of $(A\rtimes_\alpha \Z)\otimes\mathcal{T}$ generated by $A\otimes 1$ and $u\otimes s$ is called the \emph{Toeplitz algebra of the pair $(A,\alpha)$}, and is denoted by $\mathcal{T}_{(A,\alpha)}$ \cite{PV80bis}. The inclusion $A\to \mathcal{T}_{(A,\alpha)}$ is a K-equivalence, see \cite[Lemma~2.3]{PV80bis} and \cite[Theorem~2.12]{Yan}.
As a special case, when $A:=\C$, the action $\alpha$ is trivial and $u=1$, we get a K-equivalence $\C\to\mathcal{T}=\mathcal{T}_{(\C,{id})}$ given by the unit map, which instanciates \ref{ax:contr}.
\end{example}

\begin{example}
The collapse of a closed quantum subspace $\mathsf Z$ of a compact quantum space $\mathsf X$ to a point is described by the pushout diagram (i), which dually reads as the pullback diagram (ii):
\begin{equation}\label{eq:CQSquotient}
\text{(i)}\quad
\begin{tikzpicture}[scale=1.2,baseline={([yshift=-3pt]current bounding box.center)}]

\node (Z'') at (1,2) {$\mathsf X/\mathsf Z$};
\node (Y') at (0,1) {$\star$};
\node (Z') at (2,1) {$\mathsf X$};
\node (Y) at (1,0) {$\mathsf Z$};

\path[->]
		(Y') edge[right hook->] (Z'')
		(Y) edge (Y')
		(Y) edge[right hook->] (Z')
		(Z') edge (Z'');
\end{tikzpicture},
\hspace{2cm}
\text{(ii)}\quad
\begin{tikzpicture}[scale=1.2,baseline={([yshift=-3pt]current bounding box.center)}]

\node (BB) at (1,3) {$I^+$};
\node (A) at (0,2) {$\C$};
\node (B) at (2,2) {$A$};
\node (D) at (1,1) {$A/I$};

\path[-To]
	(BB) edge[->>] (A)
	(BB) edge (B)
	(A) edge (D)
	(B) edge[->>] (D);
\end{tikzpicture} .
\end{equation}
Here, $A:=C(\mathsf X)$, $I:=\ker\{C(\mathsf X)\twoheadrightarrow C(\mathsf Z)\}$, and $C(\mathsf X/\mathsf Z)\cong I^+$, where $I^+$ is the minimal unitization of $I$. It follows from \ref{ax:cwW} that the collapsing morphism $\mathsf X\to\mathsf X/\mathsf Z$ is a weak equivalence (i.e.~$I^+\to A$ is a K-equivalence) 
if and only if $\mathsf Z$ is K-contractible.
\end{example}

\begin{example}
Consider the following dual pair of pushout and pullback diagrams
\begin{equation}\label{eq:wedge}
\begin{tikzpicture}[scale=1.2]

\node (Z'') at (1,2) {$\mathsf X\vee\mathsf Y$};
\node (Y') at (0,1) {$\mathsf X$};
\node (Z') at (2,1) {$\mathsf Y$};
\node (Y) at (1,0) {$\star$};

\path[->]
		(Y') edge[right hook->] (Z'')
		(Y) edge (Y')
		(Y) edge[right hook->] (Z')
		(Z') edge (Z'');
\end{tikzpicture}\hspace{2cm}
\begin{tikzpicture}[scale=1.2]

\node (BB) at (1,3) {$A\times_{\C} B$};
\node (A) at (0,2) {$A$};
\node (B) at (2,2) {$B$};
\node (D) at (1,1) {$\C$};

\path[-To]
		(BB) edge[->>] (A)
		(BB) edge (B)
		(A) edge (D)
		(B) edge[->>] (D);
\end{tikzpicture}
\end{equation}
describing a wedge sum of compact quantum spaces,
where $A:=C(\mathsf X)$ and $B:=C(\mathsf Y)$. From \ref{ax:cwW}, we deduce that $\mathsf Y\to\mathsf X\vee\mathsf Y$ is a weak equivalence if and only if the morphism $\star\to\mathsf X$ is a weak equivalence,
i.e.~$\mathsf X$ is K-contractible (see Remark \ref{rmk:weak}).
\end{example}

The next examples instantiate the fact that KK-equivalence is weaker than K-equivalence.

\begin{example}
Let $X$ and $Y$ be any two compact Hausdorff spaces such that $K^*(X)$ and $K^*(Y)$ are isomorphic as graded abelian groups but not as graded rings. Then $C(X)$ and $C(Y)$ are KK-equivalent but not K-equivalent because a K-equivalence would induce an isomorphism of K-theory rings.
An example is given by the complex projective space
$\CP^n$ and the discrete space $P_{n+1}$ consisting of $n+1$ points. Their K-theory groups are isomorphic
but the rings $K^0(\CP^n)\cong\Z[x]/(x^{n+1})$ and $K^0(\C^{n+1})\cong\Z^{n+1}$ are not \cite{AT60}.
\end{example}

\begin{example}
For $n\geq 2$, the C*-algebras $\C$ and $M_n(\C)$
are Morita equivalent, hence KK-equivalent,
but there exists no K-equivalence between them in either direction.
Indeed, there exists no unital *-homomorphism $M_n(\C)\to\C$, and the unique unital *-homomorphism $\C\to M_n(\C)$ induces the multiplication by $n$ under the canonical identification of both $K_0(\C)$ and $K_0(M_n(\C))$ with $\Z$ given by the trace map.
\end{example}

\subsection{Weak and strict quantum CW-complexes}
Recall that a finite \emph{filtration by skeleta} of a topological space $X$ is a sequence of closed embeddings 
\begin{equation}\label{skelX}
X^{d_{0}}\longhookrightarrow X^{d_{1}}\longhookrightarrow\cdots
\longhookrightarrow X^{d_{n-1}}\longhookrightarrow X^{d_{n}}=:X .
\end{equation}
Here $X^{d_{0}}$ is a finite discrete space and, for each $1\leq k\leq n$, the $k$-th arrow in \eqref{skelX} is the top-left arrow of the pushout diagram
\begin{equation}\label{eq:attachingcells}
\begin{tikzpicture}[scale=1.7,baseline={([yshift=-3.5pt]current bounding box.center)}]

\node (A) at (1,2) {$X^{d_{k}}\!\!$};
\node (B) at (0,1) {$X^{d_{k-1}}$};
\node (C) at (2,1) {$\coprod_{i=1}^{m_{k}} {B}_i^{d_{k}}$};
\node (D) at (1,0) {$\coprod_{i=1}^{m_{k}} {S}_i^{d_{k}-1}$};

\path[->]
     	(B) edge[right hook->] (A)
		(C) edge (A)
		(D) edge (B)
		(D) edge[right hook->,font=\footnotesize,inner sep=2pt] node[below right]{$\partial$} (C);

\end{tikzpicture} \!.
\end{equation}
Here $\partial$ denotes the boundary map from a finite disjoint union of spheres to a finite disjoint union of balls, and  $0=d_{0}<d_{1}<\ldots  <d_{n}$ is the dimension sequence.
We refer to the pushout in \eqref{eq:attachingcells} as \emph{attaching cells} of dimension $d_k$.
A topological space for which a finite filtration by skeleta exists is called a \emph{finite CW-complex}. (We are interested only in finite CW-complexes because they are compact Hausdorff spaces.)

We now want to generalize this notion 
to compact quantum spaces and to allow the more flexible notion of weak cofibration based on K-equivalence.
The following definition is motivated by the K-theorical properties of the classical boundary map $S^{d-1}\to B^d$.

\begin{df}\label{def:Ksphereandball}
Let $d$ be a positive integer. An object $\mathsf{S}^{d-1}$ in $\CQS^\times$ is called a \emph{K-sphere} if it has the same K-theory as $S^{d-1}$, and an object $\mathsf{B}^{d}$ in $\CQS^\times$ is called a \emph{K-ball} if it has the same K-theory as $B^d$. 
A cofibration $\partial:\mathsf{S}^{d-1}\hookrightarrow\mathsf{B}^{d}$ from a K-sphere $\mathsf{S}^{d-1}$ to a K-ball 
$\mathsf{B}^{d}$ is called a \emph{K-boundary map}
if it induces the short exact sequence
\begin{equation}
\begin{tikzcd}[
ampersand replacement=\&,
row sep=1em,
column sep=1.5cm,
/tikz/la/.style={label={[yshift=-8mm,rotate=270,anchor=west,inner sep=2pt,label={right:\;\mathbb{Z}}]:\cong}}
]
0 \arrow[r] \&
|[la]| K^{0}(\mathsf{B}^{d}) \arrow[r]  \&
|[la]| K^{0}(\mathsf{S}^{d-1}) \arrow[r] \&
0 \arrow[r] \&
0
\end{tikzcd}\label{0even}
\end{equation}
for $d$ even, and the short exact sequence
\begin{equation}
\noeqref{0odd}
\begin{tikzcd}[
ampersand replacement=\&,
row sep=1em,
column sep=1.5cm,
/tikz/la/.style={label={[yshift=-8mm,rotate=270,anchor=west,inner sep=2pt,label={right:\;\mathbb{Z}}]:\cong}},
/tikz/lb/.style={label={[yshift=-8mm,rotate=270,anchor=west,inner sep=2pt,label={right:\;\mathbb{Z}\oplus\mathbb{Z}}]:\cong}}
]
0 \arrow[r] \&
|[la]| K^{0}(\mathsf{B}^{d}) \arrow[r]  \&
|[lb]| K^{0}(\mathsf{S}^{d-1}) \arrow[r] \&
\Z \arrow[r] \&
0
\end{tikzcd}\label{0odd}
\end{equation}
for $d$ odd. 
\end{df}

\begin{df}\label{def:finiteqKCWc}
A finite \emph{weak filtration by skeleta} of a compact quantum space $\mathsf X$ is a sequence of weak cofibrations in $\Ho(\CQS)$
\begin{equation}\label{wskelX}
\begin{tikzcd}[arrow style=tikz]
\mathsf X^{d_{0}}\ar[r, rightarrowtail] &
\mathsf X^{d_{1}}\ar[r, rightarrowtail] &
\cdots\ar[r, rightarrowtail] &
\mathsf X^{d_{n-1}}\ar[r, rightarrowtail] &
\mathsf X^{d_{n}}=\mathsf X
\end{tikzcd}
\end{equation}
such that $\mathsf X^{d_{0}}$ 
is a \emph{finite discrete quantum space} (i.e.~it is described by a finite dimensional C*-algebra),
for all $1\leq k\leq n$, the $k$-th arrow in 
\eqref{wskelX} is the top-left arrow of the weak pushout diagram
\begin{equation}\label{eq:qpush}
\begin{tikzpicture}[scale=1.7,baseline={([yshift=-3pt]current bounding box.center)}]

\node (A) at (1,2) {$\mathsf X^{d_{k}}$};
\node (B) at (0,1) {$\mathsf X^{d_{k-1}}$};
\node (C) at (2,1) {$\coprod_{i=1}^{m_{k}} \mathsf{B}^{d_{k}}_{i}$};
\node (D) at (1,0) {$\coprod_{i=1}^{m_{k}} \mathsf{S}^{d_{k}-1}_{i}$};

\path[->,font=\footnotesize,inner sep=2pt]
     	(B) edge[>->] (A)
		(C) edge (A)
		(D) edge node[below left]{$a_{k}$} (B)
		(D) edge[right hook->] node[below right]{$\partial$} (C);

\end{tikzpicture}.
\end{equation}
Here $\partial=\coprod_{i=1}^{m_{k}}\partial_i$ is a coproduct of K-boundary maps $\partial_i:\mathsf{S}^{d_{k}-1}_{i}\hookrightarrow\mathsf{B}^{d_{k}}_{i}$ from a K-sphere to a K-ball, and $a_k$ is a given morphism called the \emph{attaching morphism}.
A compact quantum space $\mathsf X$ for which a weak filtration by skeleta exists is called a \emph{finite weak
quantum CW-complex}. The filtration \eqref{wskelX} is called \emph{strict} if its arrows are induced by cofibrations in $\CQS$ and all pushout diagrams \eqref{eq:qpush} are induced by pushouts in $\CQS$.
\end{df}

\begin{rmk}
Due to the Bott periodicity of topological K-theory, the only sensible notion of dimension in Definition \ref{def:Ksphereandball} is modulo 2. Therefore, the classical notion of increasing sequence of dimensions of skeleta in our framework does not make sense. Consequently, our notion of a strict CW-complex admits more general examples even classically.
It is also worth noticing that our notion of K-boundary map from a K-sphere to a K-ball is stable under twisting by matrix algebras, as is illustrated in the example below.\end{rmk}

\begin{example}
The morphism dual to the quotient map $C(B^d)\otimes M_n(\C)\to C(S^{d-1})\otimes M_n(\C)$
is a K-boundary map
from a K-sphere to a K-ball in the sense of Definition \ref{def:Ksphereandball}.
Consequently, every NCCW complex in the sense of \cite[Section~2.4]{eil-lor-ped-98} is a finite strict quantum CW-complex in our sense.
Note that we do not assume that a K-ball is K-contractible. In fact, in this example the K-ball
is not K-contractible if $n\geq 2$.
\end{example}

The following theorem is a generalization of Atiyah's K{\"u}nneth theorem for finite CW-complexes \cite{ATI62}.

\begin{thm}\label{atiyahkun}
Every finite strict quantum CW-complex is an object in $\CQS^\times$.
\end{thm}

\begin{proof}
We will prove, by induction on $k$, that each skeleton $\mathsf X^{d_k}$ is an object in $\CQS^\times$.
For $k=0$, this is obvious because finite-dimensional C*-algebras are in the bootstrap class, which is contained in $\mathcal{N}_{\max}$.
Assume now that $X^{d_{k-1}}$ is in $\CQS^\times$ for some $k\geq 1$.
By definition, K-spheres and K-balls belong to $\CQS^\times$.
Also, our monoidal structure distributes over finite coproducts, so the morphisms $\partial$ in \eqref{eq:qpush} are in $\CQS^\times$. 
Therefore, applying \ref{ax:sub} to
the pushout diagram \eqref{eq:qpush}
in $\CQS$ yields that $\mathsf X^{d_k}$ is an object in $\CQS^\times$, thus proving the inductive step.
\end{proof}

Recall now that every one-surjective pullback diagram of unital C*-algebras has an associated Mayer--Vietoris six-term exact sequence in K-theory.
As a consequence, any weak pushout \eqref{eq:wpo} with 
$Z\cof Y$ a cofibration, has an associated Mayer--Vietoris six-term exact sequence given by the sequence induced by the outer pushout diagram and the isomorphism induced by the K-equivalence $\psi:W\weq \widetilde{W}$.
Applying this to \eqref{eq:qpush}, we get:
\[
\begin{tikzpicture}[-To,baseline={([yshift=5pt]current bounding box.south)}]

\node (a1) at (0,2) {$K^0(\mathsf X^{d_{k}})$};
\node (a2) at (6,2) {$K^0(\mathsf X^{d_{k-1}})\oplus \bigoplus_{i=1}^{m_{k}}K^0( \mathsf{B}^{d_{k}}_{i})$};
\node (a3) at (12,2) {$\bigoplus_{i=1}^{m_{k}}K^0(\mathsf{S}^{d_{k}-1}_{i})$};

\node (b1) at (0,0) {$\bigoplus_{i=1}^{m_{k}}K^1(\mathsf{S}^{d_{k}-1}_{i})$};
\node (b2) at (6,0) {$K^1(\mathsf X^{d_{k-1}})\oplus \bigoplus_{i=1}^{m_{k}}K^1( \mathsf{B}^{d_{k}}_{i})$};
\node (b3) at (12,0) {$K^1(\mathsf X^{d_{k}})$\ .};

\begin{scope}[font=\scriptsize]
\draw (a1) -- (a2);
\draw (a2) -- node[above] {$(a_{k}^{*},{}-\partial^{*})$} (a3);
\draw (a3) -- (b3);
\draw (b3) -- (b2);
\draw (b2) -- node[above] {$(a_{k}^{*},{}-\partial^{*})$} (b1);
\draw (b1) -- (a1);
\end{scope}

\end{tikzpicture}
\]
Since every homomorphism of $\mathbb{Z}/2\mathbb{Z}$-graded abelian groups $\partial_{i}^{*}: K^{*}(\mathsf{B}^{d_{k}}_{i})\rightarrow K^{*}(\mathsf{S}^{d_{k}-1}_{i})$ is an embedding onto a direct summand, and $K^1( \mathsf{B}^{d_{k}}_{i})=0$,
the above sequence boils down to the following six-term exact sequence:
\begin{equation}\label{eq:sixtermrel}
\begin{tikzpicture}[-To,baseline=(current bounding box.center)]

\node (a1) at (0,2) {$K^0(\mathsf X^{d_{k}})$};
\node (a2) at (5,2) {$K^0(\mathsf X^{d_{k-1}})$};
\node (a3) at (10,2) {$\bigoplus_{i=1}^{m_{k}}K^0(\mathsf{S}^{d_{k}-1}_{i})/K^0( \mathsf{B}^{d_{k}}_{i})$};

\node (b1) at (0,0) {$\bigoplus_{i=1}^{m_{k}}K^1(\mathsf{S}^{d_{k}-1}_{i})$};
\node (b2) at (5,0) {$K^1(\mathsf X^{d_{k-1}})$};
\node (b3) at (10,0) {$K^1(\mathsf X^{d_{k}})$\ .};

\begin{scope}[font=\scriptsize]
\draw (a1) -- (a2);
\draw (a2) -- (a3);
\draw (a3) -- node[right] {$\delta_{01}$} (b3);
\draw (b3) -- (b2);
\draw (b2) -- (b1);
\draw (b1) -- node[left] {$\delta_{10}$} (a1);
\end{scope}

\end{tikzpicture}
\end{equation}

Now we are ready to compute inductively the K-theory of finite quantum CW-complexes, provided the images and kernels of the $\delta$-maps are known.

\begin{thm}\label{Kthm}
For any finite weak filtration by skeleta as in \eqref{wskelX}, and for all $1\leq k\leq n$,
there exist non-canonical isomorphisms
\begin{alignat}{2}
& K^{0}(\mathsf X^{d_{k}}) \cong K^{0}(\mathsf X^{d_{k-1}})\oplus \ima \delta_{10}  && \text{if $d_k$ is even and $K^{0}(\mathsf X^{d_{k-1}})$ is free} , \label{for Milnor} \\
& K^{1}(\mathsf X^{d_{k-1}}) \cong K^{1}(\mathsf X^{d_{k}})\oplus \ker \delta_{10}  && \text{if $d_k$ is even} , \label{for vanish} \\
& K^{0}(\mathsf X^{d_{k-1}}) \cong K^{0}(\mathsf X^{d_{k}})\oplus \ker \delta_{01} \qquad && \text{if $d_k$ is odd} , \label{dkodd} \\
& K^{1}(\mathsf X^{d_{k}}) \cong K^{1}(\mathsf X^{d_{k-1}})\oplus \ima \delta_{01}  && \text{if $d_k$ is odd and $K^{1}(\mathsf X^{d_{k-1}})$ is free.} \label{k1odd}
\end{alignat}
\end{thm}

\begin{proof}
It follows directly from Definition \ref{def:Ksphereandball} that the six-term exact sequence \eqref{eq:sixtermrel} is equivalent to the pair of short exact sequences
\begin{gather}\label{ses0even}
0\longrightarrow
\ima \delta_{10}\longrightarrow
K^{0}(\mathsf X^{d_{k}})\longrightarrow
K^{0}(\mathsf X^{d_{k-1}})\longrightarrow
0,
\\
\label{ses1even}
\noeqref{ses1even}
0\longrightarrow
K^{1}(\mathsf X^{d_{k}})\longrightarrow
K^{1}(\mathsf X^{d_{k-1}})\longrightarrow
\ker \delta_{10}\longrightarrow
0,
\end{gather}
for $d_{k}$ even, and to the pair of short exact sequences
\begin{gather}\label{ses0odd}
\noeqref{ses0odd}
0\longrightarrow
K^{0}(\mathsf X^{d_{k}})\longrightarrow
K^{0}(\mathsf X^{d_{k-1}})\longrightarrow
\ker \delta_{01}\longrightarrow
0,
\\
\label{ses1odd}
0 \longrightarrow
\ima \delta_{01} \longrightarrow
K^{1}(\mathsf X^{d_{k}}) \longrightarrow
K^{1}(\mathsf X^{d_{k-1}}) \longrightarrow
0,
\end{gather}
for $d_{k}$ odd. Since a subgroup of a free abelian group is free abelian, the kernels of the \mbox{$\delta$-maps} are free abelian, and we have non-canonical isomorphisms \eqref{for vanish} and \eqref{dkodd}.
For $d_{k}$ even, if $K^{0}(\mathsf X^{d_{k-1}})$ is free, the short exact sequence \eqref{ses0even} splits and we have an additional non-canonical isomorphism \eqref{for Milnor}. Much in the same way,
for $d_{k}$ odd, if $K^{1}(\mathsf X^{d_{k-1}})$ is free, the short exact sequence \eqref{ses1odd} splits and we have an additional non-canonical isomorphism~\eqref{k1odd}.
\end{proof}

As an application of Theorem \ref{Kthm}, we compute the K-theory groups of a class of finite quantum CW-complexes including classical examples such as finite discrete spaces, wedge sums of even-dimensional spheres, complex projective spaces and more.

\begin{cor}
\label{KweakK}
Suppose we have a finite weak filtration by skeleta
\eqref{wskelX} of $\mathsf X$,
with $\mathsf X^{d_{0}}$ K-contractible.
Assume also that, in \eqref{eq:qpush}, $d_k$ is even and $m_k=1$ for all $1\leq k\leq n$.
Then, for all $0\leq k\leq n$,
\begin{equation}
K^0(\mathsf X^{d_k})\cong\Z^{k+1}\quad
\text{and} \quad K^1(\mathsf X^{d_k})\cong 0 .
\end{equation}
\end{cor}
\begin{proof}
We prove the theorem by induction on $k$.
It is trivially true for $k=0$. Assume, next, that it is true for $k-1$ with $k\geq 1$. Since $K^{0}(\mathsf X^{d_{k-1}})$ is free, the isomorphisms
\begin{gather}
K^{0}(\mathsf X^{d_{k}}) \cong K^{0}(\mathsf X^{d_{k-1}})\oplus \ima \delta_{10} \cong\Z^k \oplus \ima \delta_{10} , \label{eq:latter} \\
K^{1}(\mathsf X^{d_{k}})\oplus \ker \delta_{10} \cong K^{1}(\mathsf X^{d_{k-1}})\cong 0 , \label{eq:former}
\end{gather}
\noeqref{eq:latter}%
follow from \eqref{for Milnor} and \eqref{for vanish}, respectively.
Furthermore, from \eqref{eq:former}, we deduce that $K^{1}(\mathsf X^{d_{k}})$ and $\ker \delta_{10}$ are both $0$.
Hence,
\begin{equation}
\delta_{10}\colon \Z\cong K^1(\mathsf{S}^{d_{k}-1})\longrightarrow K^{0}(\mathsf X^{d_{k}})
\end{equation}
is injective, so $\ima \delta_{10}\cong\Z$. Finally,  from \eqref{eq:latter}, we obtain
$K^{0}(\mathsf X^{d_{k}})\cong \Z^{k}\oplus\Z\cong\Z^{k+1}$, thus proving the inductive step.
\end{proof}

In the proof of the preceding corollary, we obtain the following construction of the generators of K-groups.
For all $\mathsf X$ as in Corollary \ref{KweakK}, there exists a non-canonical isomorphism
\begin{equation}\label{eq:Ksplit}
(\kappa,\delta_{10}):K^0(\mathsf X^{d_{k-1}})\oplus K^1(\mathsf{S}^{d_{k}-1})
\longrightarrow K^0(\mathsf X^{d_{k}}) .
\end{equation}
Here $\delta_{10}$ is the Milnor connecting homomorphism from the six-term exact sequence 
\eqref{eq:sixtermrel} and $\kappa$ is a section of the split epimorphism
\begin{equation}\label{kisplit}
K^0(\mathsf X^{d_{k}})\longrightarrow K^0(\mathsf X^{d_{k-1}})
\end{equation}
induced by the weak cofibration in  \eqref{eq:qpush}.

Furthermore, the weak filtration \eqref{wskelX} induces a tower
of $\Z/2\Z$-graded K-groups
\begin{equation}\label{Ktower}
\Z=K^{*}(\mathsf X^{d_{0}})\leftarrow K^{*}(\mathsf X^{d_{1}})\leftarrow\cdots\leftarrow K^{*}(\mathsf X^{d_{n-1}})\leftarrow K^{*}(\mathsf X^{d_{n}}) ,
\end{equation}
where all arrows are surjective and split by \eqref{eq:Ksplit}. Better still, denoting by $g$ a generator of 
$K^1(\mathsf{S}^{d_{k}-1})\cong\Z$, we 
 call  $\delta_{10}(g)\in K^0(\mathsf X^{d_{k}})$ \emph{the Milnor generator}.
Now, inductively,  splittings of \eqref{Ktower} taken in all dimensions up to $k=n$ provide a basis of $K^0(\mathsf X)$ consisting of the last Milnor generator 
and liftings of the Milnor generators from previous skeleta.

We end with the following definition needed in the forthcoming sections.

\begin{df}\label{wceq}
A \emph{weak cellular equivalence}
is a commutative diagram
\[
\begin{tikzcd}[arrow style=tikz,column sep=1.3cm,row sep=1.3cm,labl/.style={anchor=south, rotate=90, inner sep=.5mm}]
\mathsf X^{d_{0}}\ar[r, rightarrowtail]\ar[d,"\sim" labl] &
\mathsf X^{d_{1}}\ar[r, rightarrowtail]\ar[d,"\sim" labl] &
\cdots\ar[r, rightarrowtail] &
\mathsf X^{d_{n-1}}\ar[r, rightarrowtail]\ar[d,"\sim" labl] &
\mathsf X^{d_{n}}=\mathsf X\ar[d,"\sim" labl] \\
\mathsf Y^{d'_{0}}\ar[r, rightarrowtail] &
\mathsf Y^{d'_{1}}\ar[r, rightarrowtail] &
\cdots\ar[r, rightarrowtail] &
\mathsf Y^{d'_{n-1}}\ar[r, rightarrowtail] &
\mathsf Y^{d'_{n}}=\mathsf Y
\end{tikzcd}
\]
 in $\Ho(\CQS)$ such that the rows are weak filtrations by skeleta,
the vertical arrows are isomorphisms in the homotopy category coming from weak equivalences in $\CQS$,
and $d_k$ and $d'_k$ have the same parity for all $0\leq k\leq n$.
(Recall that in this construction the dimension only makes sense modulo $2$.)
\end{df}

\section{Quantum CW-complex structures of quantized surfaces}

\noindent
In this section, we present some examples of finite strict and weak quantum CW complexes coming
from quantizations of surfaces. These include all Podle\'s quantum spheres \cite{pod87}, the mirror
quantum sphere~\cite{hms06}, the multipushout quantum real projective plane, its covering multipushot quantum 2-sphere~\cite{hrz21}, a new quantum sphere that we introduce herein, and Toeplitz compact
quantum surfaces~\cite{gw20}, which include the quantum real projective plane \cite{hms03,h-pm96}, and quantum
teardrops~\cite{bf12}. Let us stress here that all the filtrations by skeleta in this section are strict,
except for the three quantum spheres, equatorial Podle\'s, mirror, and multipushout, that have  weak
filtrations by skeleta. 
%However, there is a weak cellular equivalence between their filtration
%by skeleta and the strict filtration by skeleta of the standard Podle\'s quantum sphere.

\subsection{The Podle\'s and mirror quantum $2$-spheres}\label{sec:2spheres}
Probably the easiest non-trivial examples of finite weak and strict quantum CW-complexes are given by the Podle\'s \cite{pod87} 
and mirror quantum spheres \cite{hms06}. We denote by $C(\mathsf{S}_q^2)$,
$C(\mathsf{S}^{2,+}_{q,\infty})$ and $C(\mathsf{S}^{2,-}_{q,\infty})$ the C*-algebra of the standard Podle\'s quantum sphere, the equatorial Podle\'s quantum sphere and the mirror quantum sphere, respectively. The spaces of characters of $C(\mathsf{S}^2_q)$ and $C(\mathsf{S}^{2,\pm}_{q,\infty})$ are a point and a circle, respectively,
so these C*-algebras are not isomorphic. It is far more tricky to distiguish $C(\mathsf{S}^{2,+}_{q,\infty})$ from $C(\mathsf{S}^{2,-}_{q,\infty})$, but we know from
\cite[Proposition~2.3]{hms06} that they are not even Morita equivalent. However, we will see  that all these three C*-algebras are K-equivalent.

To this end, we  combine the pullback presentations of these C*-algebras (see \cite[Proposition~1.2]{sheu1991quantization} and \cite{hms06}) 
 into the following
commutative diagram of unital C*-algebras and unital *-homomorphisms:
\begin{equation}\label{eq:compareagainS2}
\begin{tikzpicture}[scale=1.7,baseline={([yshift=-5pt]current bounding box.center)}]

\node (A) at (1.1,2) { $C(\mathsf S_q^2)$ };
\node (B) at (0,1) { $\C$ };
\node (C) at (2.2,1) { $\mathcal{T}$ };
\node (D) at (1.1,0) { $C({S}^1)$ };
\node (E) at (5.6,2) { $C(\mathsf S^{2\pm}_{q,\infty})$ };
\node (F) at (4.5,1) { $\mathcal{T}$ };
\node (G) at (6.7,1) { $\mathcal{T}$ };
\node (H) at (5.6,0) { $C({S}^1)$ };

\path[-To,font=\footnotesize] (A) edge (B)
			(A) edge (C)
			(B) edge node[below left] {$\eta'$} (D)
			(C) edge node[below right,pos=0.7] {$\sigma$} (D)
			(E) edge (F)
			(E) edge (G)
			(F) edge node[below left] {$\sigma^\pm$}  (H)
			(G) edge node[below right] {$\sigma$} (H)
			(A) edge[bend left=10]  node[above] {$\alpha^\pm$} (E)
			(B) edge[bend right=20] node[below,pos=0.7] {$\eta$} (F)
			(C) edge[bend left=20] node[above,pos=0.4] {$\id$} (G)
			(D) edge[bend right=10] node[below] {$\id$} (H);
			
\end{tikzpicture}.
\end{equation}
Here $\mathcal{T}$ is the Toeplitz C*-algebra, 
i.e.~the universal C*-algebra generated by an isometry $s$,
$\sigma:\mathcal{T}\to C(S^1)$ is the \emph{symbol map} sending $s$ to the standard unitary generator $u$ of $C(S^1)$ (the restriction to $S^1$ of the identity 
function $\C\to\C$). Next, $\eta$ and $\eta'$ are the unit maps, $\sigma^+=\sigma$, and $\sigma^-$ is the *-homomorphism uniquely defined by
$\sigma^-(s):=u^*$.
The unlabelled maps are canonically defined by the pullback structure, and $\alpha^\pm$ exist by the universal property of pullbacks.
Since $\sigma$ is surjective and $\eta$ is a K-equivalence (Example \ref{ex:toeplitz}),
it follows from \ref{ax:weq} that $\alpha^\pm$ are also K-equivalences.

In the opposite category $\CQS$, the cofibration $\partial:{S}^1\to\mathsf{D}_q$ dual
to the quotient map $\sigma:\mathcal{T}=C(\mathsf{D}_q)\to C({S}^1)$ is a K-boundary map from a K-sphere to a K-ball in the sense of
 Definition~\ref{def:Ksphereandball}.
Indeed, both $C(S^1)$ and $\mathcal{T}$ are in the bootstrap class, and one can easily prove that $\sigma$ induces the sequence \eqref{0even}.
Consequently, from the left pullback in \eqref{eq:compareagainS2} we get the strict filtration by skeleta:
\begin{equation}\label{eq:CWstrPod}
\star\longhookrightarrow \mathsf{S}^2_q\,.
\end{equation}

Next, contracting the identity maps in \eqref{eq:compareagainS2}, we obtain the following diagram:
\begin{equation}\label{eq:fancysquares}
\begin{tikzpicture}[scale=1.5,baseline=(current bounding box.center),>=To]

\node (A) at (0,2) { $C(\mathsf S_q^2)$ };
\node (B) at (0,0) { $\C$ };
\node (E) at (3,2) { $C(\mathsf S^{2\pm}_{q,\infty})$ };
\node (F) at (3,0) { $\mathcal{T}$ };
\node (G) at (6,2) { $\mathcal{T}$ };
\node (H) at (6,0) { $C({S}^1)$ };

\path[-To,font=\footnotesize] (A) edge[->>] (B)
			(A) edge[bend left=20] (G)
			(B) edge[bend right=20] node[below] {$\eta'$} (H)
			(E) edge[->>] (F)
			(E) edge (G)
			(F) edge node[above] {$\sigma^\pm$}  (H)
			(G) edge[->>] node[right] {$\sigma$} (H)
			(A) edge node[below] {$\alpha^\pm$} (E)
			(B) edge node[above] {$\eta$} (F);
			
\end{tikzpicture}.
\end{equation}
Since the square on the right and the outer rectangle are pullback diagrams, the left square is also a pullback diagram \cite[p.~72, Exercise~8.b]{MacLane}, and its dual diagram
\begin{equation}\label{eq:pushoutS2S2pm}
\begin{tikzpicture}[baseline=(current bounding box.center),scale=1.3]

\node (A) at (1,2) { $\mathsf{S}^2_q$ };
\node (B) at (0,1) { $\star$ };
\node (C) at (2,1) { $\mathsf{S}^{2,\pm}_{q,\infty}$ };
\node (D) at (1,0) { $\mathsf{D}_q$ };

\path[-To,font=\footnotesize] 
	(B) edge[right hook-To] (A)
	(C) edge node[above,sloped] {$\sim$} (A)
	(D) edge node[above,sloped] {$\sim$} (B)
	(D) edge[right hook-To] (C);
			
\end{tikzpicture}
\end{equation}
is a pushout. Furthermore, since the morphism $\mathsf{D}_q\to\star$ is a weak equivalence (Remark \ref{rem:trivialK}), the morphism $\mathsf{S}^{2,\pm}_{q,\infty}\to\mathsf{S}^2_q$ is also a weak equivalence by the axiom \ref{ax:cwW}.
Finally, the roof of the diagram \eqref{eq:pushoutS2S2pm} defines a weak cofibration
\[
\begin{tikzcd}[arrow style=tikz]
\star\ar[r, rightarrowtail] & \mathsf{S}^{2,\pm}_{q,\infty}
\end{tikzcd} ,
\]
which is a weak filtration by skeleta of $\mathsf{S}^{2,\pm}_{q,\infty}$.

%The diagram \eqref{eq:pushoutS2S2pm} in the homotopy category becomes a cellular weak equivalence
%\begin{equation}\label{eq:weqskS2}
%\begin{tikzpicture}[baseline=(current bounding box.center)]
%\node (A) at (0,2.2) { $\star$ };
%\node (B) at (0,0) { $\star$ };
%\node (C) at (3,2.2) { $\mathsf{S}^{2,\pm}_{q,\infty}$ };
%\node (D) at (3,0) { $\mathsf{S}^2_q$ };
%\path[font=\footnotesize]
%		(A) edge[double equal sign distance] (B)
%		(A) edge[>-To] (C)
%		(C) edge[-To] node[above,sloped] {$\sim$} (D)
%		(B) edge[right hook-To] (D);			
%\end{tikzpicture}
%\end{equation}
%from the weak filtration of $\mathsf{S}^{2,\pm}_{q,\infty}$ to the strict filtration of $\mathsf{S}^2_q$.
%For a geometric interpretation of this weak equivalence, we can rewrite \eqref{eq:pushoutS2S2pm} in the following form
%\begin{equation}\label{eq:tubB2}
%\begin{tikzpicture}[baseline=(current bounding box.center),inner sep=2pt]
%\node (A) at (0,2.2) { $\star$ };
%\node (B) at (0,0) { $\star$ };
%\node (C) at (3,2.2) { $\mathsf{S}^{2,\pm}_{q,\infty}$ };
%\node (D) at (1.5,1.1) { $\mathsf{D}_q$ };
%\node (E) at (3,0) { $\mathsf{S}^2_q$ };
%\path[font=\footnotesize] 
%		(A) edge[double equal sign distance] (B)
%		(B) edge[right hook-To] (E)
%		(C) edge[-To] node[above,sloped] {$\sim$} (E)
%		(D) edge[-To] node[above,sloped] {$\sim$} (A)
%		(D) edge[right hook-To] (C);
%\end{tikzpicture}
%\end{equation}
%to see that the weak equivalence collapses a tubular neighbourhood of the $0$-skeleton of $\mathsf{S}^{2,\pm}_{q,\infty}$.

\subsection{Toeplitz compact quantum surfaces}\label{sec:toepsurf}
Compact quantum orientable surfaces $\mathsf{T}_{g,q}$ of genus $g\geq 0$ and compact quantum non-orientable surfaces $\mathsf{P}_{g,q}$ of demi-genus $g\geq 1$
were introduced in \cite{gw20}. In particular, 
$\mathsf{P}_{1,q}=\RPq^2_q$ is the quantum real projective plane of \cite{hms03,h-pm96}.
They are defined by the following pullback diagrams
of unital C*-algebras and unital *-homomorphism:
\[
\begin{tikzpicture}[baseline={([yshift=-4pt]current bounding box.center)},scale=1.5]

\node (A) at (1,2) {$C(\mathsf{T}_{g,q})$};
\node (B) at (0,1) {$C(\bigvee^{2g}{S}^1)$};
\node (C) at (2,1) {$\mathcal{T}$};
\node (D) at (1,0) {$C({S}^1)$};

\path[-To,font=\footnotesize]
     (A) edge[->>] (B)
		(A) edge[right hook-To] (C)
		(B) edge[right hook-To] node[below left] {$\rho^+_g$} (D)
		(C) edge[->>] node[below right] {$\sigma$} (D);

\end{tikzpicture},
\hspace{1.5cm}
\begin{tikzpicture}[baseline={([yshift=-4pt]current bounding box.center)},scale=1.5]

\node (A) at (1,2) {$C(\mathsf{P}_{g,q})$};
\node (B) at (0,1) {$C(\bigvee^{g}{S}^1)$};
\node (C) at (2,1) {$\mathcal{T}$};
\node (D) at (1,0) {$C({S}^1)$};

\path[-To,font=\footnotesize]
     (A) edge[->>] (B)
		(A) edge[right hook-To] (C)
		(B) edge[right hook-To] node[below left] {$\rho^-_g$} (D)
		(C) edge[->>] node[below right] {$\sigma$} (D);

\end{tikzpicture}.
\]
Here, $\sigma$ is the symbol map and $\rho^\pm_g$ are defined in \cite[Section 2]{gw20}.
%From a dual point of view, we have the pushout diagrams
%\[
%\begin{tikzpicture}[baseline=(current bounding box.center),scale=1.7]
%
%\node (A) at (1,2) {$\mathsf{T}_{g,q}$};
%\node (B) at (0,1) {$\bigvee^{2g}{S}^1$};
%\node (C) at (2,1) {$\mathsf{D}_q$};
%\node (D) at (1,0) {${S}^1$};
%
%\path[-To]
%     	(B) edge[right hook->] (A)
%		(C) edge (A)
%		(D) edge (B)
%		(D) edge[right hook->] node[below right,font=\footnotesize] {$\partial$} (C);
%
%\end{tikzpicture}
%\hspace{1.5cm}
%\begin{tikzpicture}[baseline=(current bounding box.center),scale=1.7]
%
%\node (A) at (1,2) {$\mathsf{P}_{g,q}$};
%\node (B) at (0,1) {$\bigvee^g{S}^1$};
%\node (C) at (2,1) {$\mathsf{D}_q$};
%\node (D) at (1,0) {${S}^1$};
%
%\path[-To]
%     	(B) edge[right hook->] (A)
%		(C) edge (A)
%		(D) edge (B)
%		(D) edge[right hook->] node[below right,font=\footnotesize] {$\partial$} (C);
%
%\end{tikzpicture}
%\]
The dual pushout diagrams yield the following filtrations by skeleta:
\[
\star\cof {\textstyle\bigvee^{2g}{S}^1}\cof\mathsf{T}_{g,q} \;,\qquad
\star\cof {\textstyle\bigvee^{g}{S}^1}\cof\mathsf{P}_{g,q}
\;.
\]
In the special case $g=1$, the second filtration becomes
$\star\hookrightarrow\RP^1\hookrightarrow\RPq_q^2$.
 Let us note here that the general  Hong--Szyma\'nski quantum real projective spaces $\RPq^{n}_q$ \cite{HS02} also enjoy having quantum CW-complex 
 structures, which can be shown using graph C*-algebras
 in the even case, and using topological-graph C*-algebras in the odd case~\cite{atulphd}.
 
\subsection{The multipushout quantum real projective plane and its covering sphere}\label{sec:multirp}
Multipushout quantum complex projective spaces prompted the development of our  weak quantum CW-complex
theory, and yield its fundamental example.  In the
previous sections, we provided a strict quantum CW-complex structure for $\RPq^{2}_q=\mathsf{P}_{1,q}$
and a weak quantum CW-complex structure its covering sphere $\mathsf{S}^{2,+}_{q,\infty}$. Much in the same way, herein we will provide a 
strict quantum CW-complex structure for 
the multipushout quantum real projective plane $\RPq^{2}_{\mathcal{T}}$ and a weak quantum CW-complex structure its covering
sphere~$\mathsf{S}^{2}_{\R,\mathcal{T}}$.

Roughly speaking, the multipushout quantum real projective plane $\RPq^{2}_{\mathcal{T}}$  
is constructed by gluing two quantum discs $\mathsf D_q$ to obtain a quantum M{\"o}bius strip,
and then gluing to it the third quantum disc $\mathsf D_q$ over its boundary circle. Its covering sphere $\mathsf{S}^{2}_{\R,\mathcal{T}}$ is
obtained by gluing three pairs of quantum discs to obtain a quantum cube surface equipped with the antipodal $\mathbb{Z}/2\mathbb{Z}$-action,
so $\RPq^{2}_{\mathcal{T}}=\mathsf{S}^{2}_{\R,\mathcal{T}}/(\mathbb{Z}/2\mathbb{Z})$ by design. Better still, a multipushout quantum disc 
$\mathsf D_\mathcal{T}$ is constructed out of three usual quantum discs $\mathsf D_q$, and a new quantum sphere $\widetilde{\mathsf S}_{\R,\mathcal{T}}^2$
is obtained by shrinking its boundary circle to a point. This new sphere plays the role of the standard Podle\'s quantum sphere in the multipushout   
setting. Moreover, we use the fact that $\mathsf{S}^{2}_{\R,\mathcal{T}}$ 
can be assembled from two 
discs $\mathsf D_\mathcal{T}$, much as $\mathsf{S}^{2,+}_{q,\infty}$ is obtained by gluing two discs $\mathsf D_q$. Again much as before,
we use the crucial fact that $\RPq^{2}_{\mathcal{T}}$ is the quotient of $\mathsf D_\mathcal{T}$ by the antipodal  $\mathbb{Z}/2\mathbb{Z}$-action on its boundary
circle. 

More precisely, recall that there are two quotient maps $\sigma_i^\circ:C(\mathsf D_{\mathcal{T}})\to C(S^1)$, $i=1,2$, related by the antipodal map on $S^1$
(see~\cite[(3.47)]{hrz21} and the line above).
From~\cite[(3.60)]{hrz21} and \ref{ax:sub} it follows that $\mathsf D_{\mathcal{T}}$ is an object in $\CQS^\times$.
Next, $K^0(\mathsf D_{\mathcal{T}})\cong\Z$ and $K^1(\mathsf D_{\mathcal{T}})\cong 0$ \cite[Lemma~3.5]{hrz21}.
Moreover, as $\sigma_1^\circ$ is a unital *-homomorphism and the $K^0$-group of $S^1$ is generated by the constant function $1$, 
it follows that the induced homomorphism
\begin{equation}\label{eq:K0sigma}
K^0(\sigma_1^\circ):K^0(\mathsf D_{\mathcal{T}})\cong\Z\longrightarrow\Z\cong K^0(S^1)
\end{equation}
is surjective, and so it is an isomorphism. Thus, \eqref{0even} is satisfied, and the dual cofibration $S^1\to\mathsf D_{\mathcal{T}}$ is a K-boundary map in the sense of Definition~\ref{def:Ksphereandball}.

Now, in analogy with \eqref{eq:compareagainS2}, we have the following commutative diagram:
\begin{equation}\label{newsphere}
\begin{tikzpicture}[scale=1.7,baseline=(current bounding box.center)]

\node (A) at (1.1,2) { $C(\widetilde{\mathsf S}_{\R,\mathcal{T}}^2)$ };
\node (B) at (0,1) { $\C$ };
\node (C) at (2.2,1) { $C(\mathsf D_{\mathcal{T}})$ };
\node (D) at (1.1,0) { $C({S}^1)$ };
\node (E) at (5.6,2) { $C(\mathsf S^2_{\R,\mathcal{T}})$ };
\node (F) at (4.5,1) { $C(\mathsf D_{\mathcal{T}})$ };
\node (G) at (6.7,1) { $C(\mathsf D_{\mathcal{T}})$ };
\node (H) at (5.6,0) { $C({S}^1)$ };

\path[-To,font=\footnotesize] (A) edge[->>] (B)
			(A) edge (C)
			(B) edge node[below left] {$\eta'$} (D)
			(C) edge[->>] node[below right,pos=0.7] {$\sigma_1^\circ$} (D)
			(E) edge[->>] node[above left] {$\pi_2$} (F)
			(E) edge node[above right] {$\pi_1$} (G)
			(F) edge node[below left] {$\sigma_2^\circ$}  (H)
			(G) edge[->>] node[below right] {$\sigma_1^\circ$} (H)
			(A) edge[bend left=10]  node[above] {$\beta$} (E)
			(B) edge[bend right=20] node[below,pos=0.7] {$\eta$} (F)
			(C) edge[bend left=20] node[above,pos=0.4] {$\id$} (G)
			(D) edge[bend right=10] node[below] {$\id$} (H);
			
\end{tikzpicture}.
\end{equation}
Here, the left square is a pullback diagram defining the C*-algebra $C(\widetilde{\mathsf S}_{\R,\mathcal{T}}^2)$
of a new quantum 2-sphere, and the right square is the pullback diagram in~\cite[(3.44)]{hrz21}.
Much as for the standard Podle\'s quantum sphere, from the left pullback in \eqref{newsphere}, we get the strict filtration by skeleta:
\begin{equation}
\star\longhookrightarrow \widetilde{\mathsf S}_{\R,\mathcal{T}}^2\,.
\end{equation}

Next, contracting the identity maps in the previous diagram, we obtain the following diagram:
\begin{equation}
\begin{tikzpicture}[scale=1.5,baseline=(current bounding box.center),>=To]

\node (A) at (0,2) { $C(\widetilde{\mathsf S}^2_{\R,\mathcal{T}})$ };
\node (B) at (0,0) { $\C$ };
\node (E) at (3,2) { $C(\widetilde{\mathsf S}^2_{\R,\mathcal{T}})$ };
\node (F) at (3,0) { $C(\mathsf D_{\mathcal{T}})$ };
\node (G) at (6,2) { $C(\mathsf D_{\mathcal{T}})$ };
\node (H) at (6,0) { $C({S}^1)$ };

\path[-To,font=\footnotesize] (A) edge[->>] (B)
			(A) edge[bend left=20] (G)
			(B) edge[bend right=20] node[below] {$\eta'$} (H)
			(E) edge[->>] (F)
			(E) edge (G)
			(F) edge node[above] {$\sigma_2^\circ$}  (H)
			(G) edge[->>] node[right] {$\sigma_1^\circ$} (H)
			(A) edge node[below] {$\beta$} (E)
			(B) edge node[above] {$\eta$} (F);
			
\end{tikzpicture}.
\end{equation}
Since the square on the right and the outer rectangle are pullback diagrams, the left square is also a pullback diagram \cite[p.~72, Exercise~8.b]{MacLane}, 
and its dual diagram
\begin{equation}
\begin{tikzpicture}[baseline=(current bounding box.center),scale=1.4]

\node (A) at (1,2) { $\widetilde{\mathsf S}^2_{\R,\mathcal{T}}$ };
\node (B) at (0,1) { $\star$ };
\node (C) at (2,1) { $\mathsf S^2_{\R,\mathcal{T}}$ };
\node (D) at (1,0) { $\mathsf{D}_{\mathcal{T}}$ };

\path[-To,font=\footnotesize] 
	(B) edge[right hook-To] (A)
	(C) edge node[above,sloped] {$\sim$} (A)
	(D) edge node[above,sloped] {$\sim$} (B)
	(D) edge[right hook-To] (C);
			
\end{tikzpicture}
\end{equation}
is a pushout. We thus get a diagram analogous to \eqref{eq:fancysquares}, from which we obtain the following weak filtration by skeleta:
\[
\begin{tikzcd}[arrow style=tikz]
\star\ar[r, rightarrowtail] & \mathsf{S}^2_{\R,\mathcal{T}}
\end{tikzcd} .
\]

Finally, we identify the subalgebra of $C(S^1)$ of elements invariant under the antipodal map with~$C(\RP^1)$.
The C*-subalgebra $C(\mathsf D_{\mathcal{T}})^+\subset C(\mathsf D_{\mathcal{T}})$ is defined in \cite[(3.50)]{hrz21} as the preimage under $\sigma_1^\circ$ of $C(\RP^1)$.
In other words, $C(\mathsf D_{\mathcal{T}})^+$ is defined by the pullback diagram
\begin{equation}\label{eq:taupull}
\begin{tikzpicture}[baseline=(current bounding box.center),scale=1.5]

\node (A) at (1,2) {$C(\mathsf D_{\mathcal{T}})^+$};
\node (B) at (0,1) {$C(\RP^1)$};
\node (C) at (2,1) {$C(\mathsf D_{\mathcal{T}})$};
\node (D) at (1,0) {$C({S}^1)$};

\path[-To,font=\footnotesize]
     (A) edge[->>] (B)
		(A) edge[right hook-To] (C)
		(B) edge[right hook-To] (D)
		(C) edge[->>] node[below right,inner sep=1pt,pos=0.45] {$\sigma_1^\circ$} (D);

\end{tikzpicture}
\end{equation}
By combining \cite[Lemma 3.4]{hrz21} and \cite[(3.18)]{hrz21}, we see that $C(\mathsf D_{\mathcal{T}})^+\cong C(\RPq^2_{\mathcal{T}})$. Hence, \eqref{eq:taupull}
gives the following strict filtration by skeleta
\[
\star\longhookrightarrow\RP^1\longhookrightarrow \RPq^2_{\mathcal{T}} .
\]

\subsection{Quantum teardrops}\label{sec:exotic}
In this section, we present an example enjoying a peculiar quantum effect: quantum teardrops $\WP^1_q(1,l)$ as defined in \cite{bf12} can  be
viewed as multiwedges $\bigvee^l\mathsf S^2_q$ of standard Podle\'s quantum spheres (see the graph \cite[(3.2)]{bs18}) despite the fact that, for $l>1$, 
the classical teardrop
$W\!P^1(1,l)\cong S^2$ is not homeomorphic, or even K-equivalent,
to  the multiwedge $\bigvee^l S^2$ of classical 2-spheres.

Let $l\geq 1$. The quantum teardrop $\WP^1_q(1,l)$ is defined in \cite{bf12} as the quotient of $\mathsf{SU}_q(2)$ by a suitable weighted action of~${U}(1)$.
(Alternatively, we can arrive at $\WP^1_q(1,l)$ by taking a $U(1)$-quotient of the quantum lens space $\mathsf L^3_q (l;1,l)$ whose graph-algebraic presentation is 
in~\cite{HS03,GZ24}.)
Observe that $\WP^1_q(1,1)=\mathsf{S}^2_q$.
For $k \geq 2$, there is a pullback diagram \cite[Proposition~4.4]{adht17} that defines the pushout:
\begin{center}
\begin{tikzpicture}[scale=1.6,baseline={([yshift=-4pt]current bounding box.center)}]

\node (A) at (1,2) {$\WP^1_q(1,k )$};
\node (B) at (0,1) {$\WP^1_q(1,k -1)$};
\node (C) at (2,1) {$\mathsf{D}_q$};
\node (D) at (1,0) {${S}^1$};

\path[->]
     	(B) edge (A)
		(C) edge (A)
		(D) edge (B)
		(D) edge[right hook-To] node[below right,inner sep=1pt,font=\footnotesize] {$\partial$} (C);

\end{tikzpicture} .
\end{center}
Next, since $\partial$ is a cofibration, it follows from the axiom \ref{ax:cof} that so is the morphism $\WP^1_q(1,k -1)\to\WP^1_q(1,k )$. Therefore, we have a strict 
filtration by skeleta
\begin{equation}\label{eq:WP}
\star\cof\mathsf{S}^2_q\cof\WP^1_q(1,2)
\cof\ldots
\cof\WP^1_q(1,l-1)
\cof\WP^1_q(1,l) \;.
\end{equation}
Now we can use our Corollary~\ref{KweakK}  to reprove the fact $K^0(\WP^1_q(1,k)))\cong\Z^{k+1}$ (see~\cite[Corollary 5.3]{bf12}).

\section{Quantum CW-complex structures of quantum complex projective spaces}
\noindent
Complex projective spaces enjoy very simple and natural CW-complex structure: 
the standard  projective linear embedding $\CP^{n-1}\hookrightarrow\CP^{n}$ can be obtained from the following pushout diagram:
\begin{equation}\label{classpull}
\begin{tikzpicture}[<-,scale=1.5,baseline={(current bounding box.center)}]

\node (A) at (1,1) {$\CP^n$};
\node (B) at (0,0) {$\CP^{n-1}$};
\node (C) at (2,0) {$B^{2n}$};
\node (D) at (1,-1) {$S^{2n-1}$};

\path
	(A) edge[<-right hook] (B)
	(A) edge (C)
	(B) edge (D)
	(C) edge[<-right hook] (D);
			
\end{tikzpicture}\,.
\end{equation}
Here the left bottom arrow is the defining Hopf fibration and the right bottom arrow is embedding the sphere into the ball as its boundary.
We thus have a simple but fundamental mechanism of attaching higher dimensional cells that yields the following filtration by skeleta:
\begin{equation}
\star=\CP^0\cof \CP^1\cof\CP^2\cof\ldots\cof\CP^n.
\end{equation}
Having this vantage point one can not only compute the K-groups of complex projective spaces, but also provide free generators of $K^0(\CP^n)$ 
 using  splittings of the canonical surjections $K^0(\CP^{k+1})\twoheadrightarrow K^0(\CP^{k})$ and  the Milnor connecting 
homomorphisms applied to generators of $K^1$ of odd spheres.

The goal of this section is to obtain quantum CW-complex structures for two different types of quantum complex projective spaces.
While the K-theory of the Vaksman--Soibelman quantum complex projective spaces $\CPq^n_q$ has been studied extensively
(e.g., see \cite{d2010bounded,AHT20,abl15}), including the Atiyah--Todd identities \cite{d2026k}, even the computation of the K-groups of the multipushout 
quantum complex projective
spaces \cite{hkz12} was not easy, and was accomplished only in~\cite{hnpsz18}. The Vaksman--Soibelman part interprets pullback results in \cite{adht17} 
as a strict quantum
CW-complex structure. Thus, its primary function is to serve as a fundamental example of our general theory. On the other hand, the weak quantum CW-complex
of the multipushout quantum complex projective spaces $\CPq^n_H$ provides us with  new tools to unravel their K-theory. Now  we can use the Milnor 
connecting homomorphisms, just as in the classical setting,
% with non-canonical injective group homomorphisms $K_0(C(\CPq^k_H))\hookrightarrow K_0(C(\CPq^{k+1}_H))$ 
not only to compute the K-groups, but also to provide free
generators of  $K_0(C(\CPq^n_H))$, thus extending 
to arbitrary $n$ the $n=2$ case study of~\cite{d2025milnor}.

All building blocks on the Vaksman--Soibelman side are graph C*-algebras, which greatly facilitates working with them. We have a $U(1)$-equivariant pullback
diagram of the C*-algebras of quantum balls and spheres, and  restricting this diagram to fixed-point subalgebras yields a pullback of graph C*-algebras that is 
a $q$-defomation of
the pullback diagram dualizing the pushout diagram~\eqref{classpull}.
Now we can follow the classical path to obtain a strict quantum CW-complex structure. However, on the multipushout
side, we lack a *-homomorphism $C(\CPq^{n}_H)\to C(\CPq^{n-1}_H)$. All we have instead is
\begin{equation}\label{eq:geom}
C(\CPq^{n}_H)\relbar\joinrel\twoheadrightarrow \big(C(\mathsf S^{2n-1}_H)\otimes\mathcal{T}\big)^{U(1)}\stackrel{\text{\raisebox{-3pt}{\scriptsize }}}{\longhookleftarrow} C(\CPq^{n-1}_H) \,,
\end{equation}
where the left arrow is a surjective *-homomorphism  and the right arrow is a K-equivalence. This precisely what led us to formally invert *-homomorphisms that
induce isomorphisms in K-theory, and prompted the idea of the cw-Waldhausen category.

We can give the following geometrical interpretation to \eqref{eq:geom}.
In the classical setting, one can choose the unit-disk tubular neighbourhood 
$\ctub{n-1}{n}$ of the linear embedding \mbox{$\CP^{n-1}\hookrightarrow\CP^n$}, identified with the unit-disk neighbourhood 
of the zero section in the normal bundle of this embedding. 
%Since this normal bundle is isomorphic to a direct sum of copies of the Hopf line bundle, each equipped with its standard Hermitian structure, we can 
%put on the normal bundle the associated max norm of their direct sum. This results in a compact neighbourhood of the zero section in the normal bundle 
%consisting of normal polydisks,
%that are closed unit balls in this norm. 
The tubular neighbourhood can be understood as the associated bundle of disks:
%:, acted on diagonally by ${U}(1)$, associated with the principal ${U}(1)$-
%bundle ${S}^{2n+1}\to\CP^n$ of unit circles in the complex Hopf line bundle, which leads to the identification
\[
\CP^n\supset\ctub{n-1}{n}\cong{S}^{2n-1}\times^{{U}(1)}{D}.
\]
Combining the above inclusion with the projection of the  disk bundle, we obtain
\begin{equation}
\CP^{n}\longhookleftarrow {S}^{2n-1}\times^{{U}(1)}{D}\relbar\joinrel\twoheadrightarrow \CP^{n-1},
\end{equation}
which is the classical version of~\eqref{eq:geom}. Better still, as the disk $D$ is contractible, \mbox{${S}^{2n-1}\times^{{U}(1)}{D}$} is homotopic to
${S}^{2n-1}/{U}(1)\cong \CP^{n-1}$, so the pullback of the right projection 
\begin{equation}
\big(C(S^{2n-1})\otimes C(D)\big)^{U(1)} \longhookleftarrow C(\CP^{n-1})
\end{equation}
is a K-equivalence. It is important to stress here that the pullback of the zero section of the disk bundle yields also a *-homomomorphism
\begin{equation}
\big(C(S^{2n-1})\otimes C(D)\big)^{U(1)} \relbar\joinrel\twoheadrightarrow C(\CP^{n-1}),
\end{equation}
but the lack of the middle point of the quantum disk $\mathsf D_q$ precludes the existence of the zero section of the quantum-disk quantum bundle 
causing the lack of such a *-homomomorphism
in the noncommutative setting. Finally, let us mention that we will need the more general *-homomorphisms
\begin{equation}\label{eq:geom2}
C(\CPq^{N}_H)\relbar\joinrel\twoheadrightarrow \big(C(\mathsf{S}^{2n-1}_H)\otimes\mathcal{T}^{N-n+1}\big)^{{U}(1)}\stackrel{\text{\raisebox{-3pt}{\scriptsize }}}{\longhookleftarrow} C(\CPq^{n-1}_H) \,,
\end{equation}
with the right arrow being a K-equivalence. Here the classical interpretation is given by a bundle of polydisks indentified with a compact tubular neighbourhood
of $\CP^{n-1}$ in $\CP^{N}$:
\begin{equation}
\CP^{N}\supset\ctub{n-1}{N}\cong {S}^{2n-1}\times^{{U}(1)}{D^{\times N-n+1}}\relbar\joinrel\twoheadrightarrow \CP^{n-1}.
\end{equation}
This justifies our definition of the quantum tubular neighbourhood~$\tub{n}{N}$:
\begin{equation}
C(\tub{n}{N}):=\big(C(\mathsf S^{2n+1}_H)\otimes\mathcal{T}^{\otimes N-n}\big)^{U(1)}.
\end{equation}

\subsection{Strict quantum CW-complex structure of \texorpdfstring{$\CPq_q^n$}{CPqn}}\label{sec:strictCW}
\noindent
In this section, we first recall the graph C*-algebra realization \cite{HS02} of Vaksman--Soibelman \cite{VS91}
quantum odd spheres  $\mathsf{S}^{2n+1}_q$ and Hong--Szyma{\'n}ski  \cite{HS08} quantum even balls~$\mathsf{B}^{2n}_q$.
 Then we will recall a pullback realization of $C(\mathsf{S}^{2n+1}_q)$, and the corresponding pullback realization of the C*-algebra $C(\CPq_q^n)$
 of the Vaksman--Soibelman quantum complex projective space
 $\CPq_q^n$, which is defined as the ${U}(1)$-invariant subalgebra of $C(\mathsf{S}^{2n+1}_q)$ for the diagonal (gauge) action.
The latter pullback is at the core of the strict quantum CW-complex structure of~$\CPq_q^n$.
The strict filtration by skeleta will come from 
the *-homomorphism
$ C(\mathsf{B}^{2n}_q)\to C(\mathsf{S}^{2n-1}_q)$ given by the one-sink extension of the graph used to obtain $C(\mathsf{S}^{2n-1}_q)$ as a graph C*-algebra.

As all the above C*-algebras are graph C*-algebras of row-finite graphs, first we need to recall some basic definition.
Let $E=(E^0,E^1,s,r)$ be a directed graph,
where $E^1$ is the set of edges, $E^0$ the set of vertices, and $s,r:E^1\to E^0$ are the source and range maps, respectively.
The graph $E$ is called \emph{row finite} if $s^{-1}(v)$ is a finite set for all $v\in E^0$.
A \emph{sink} is a vertex $v$ such that $s^{-1}(v)=\emptyset$, and a \emph{loop} is an edge $e$ such that $r(e)=s(e)$.
It is known~\cite{bates2000c} that the \emph{graph C*-algebra} $C^*(E)$ of a row-finite graph $E$ is the universal \mbox{C*-algebra} generated by mutually orthogonal projections $\{P_v\,|\,v\in E^0\}$ and partial isometries $\{S_e\,|\,e\in E^1\}$ subject to the Cuntz--Krieger relations:
\[
\forall\;e\in E^1:\;S_e^*S_e =P_{r(e)} \,, \qquad
\forall\;v\in s(E^1):\sum_{e\in s^{-1}(v)}S_eS_e^*=P_v\,.
\]
Finally, all graph C*-algebra admit the so-called \emph{gauge action}, $\gamma:U(1)\to\mathrm{Aut}(C^*(E))$, given on generators by $\gamma_\lambda(S_e):=\lambda S_e$ and $\gamma_\lambda(P_v):=P_v$, for all $\lambda\in U(1)$, $e\in E^1$ and $v\in E^0$. The C*-subalgebra of gauge-invariant elements is called the \emph{core} of $C^*(E)$,
and is denoted by $C^*(E)^{U(1)}$.

We begin with fixing notation for quantum balls and spheres. We identify the C*-algebra $C(\mathsf{S}^{2n+1}_q)$, $n\geq 0$, with the graph \mbox{C*-algebra} of the graph $\Lambda_{2n+1}$:
\begin{equation}\label{lambda}
\begin{tikzpicture}[baseline={([yshift=7pt]current bounding box.center)}]

\clip (-0.6,-2.6) rectangle (9.7,1.7);

      \node[main node] (1) {};
      \node (2) [right of=1] {};
      \node (3) [right of=2] {};
      \node (4) [right of=3] {};
      \node (5) [right of=4] {};
      \node (6) [right of=5] {};

      \filldraw (1) circle (0.05) node[below left] {\footnotesize $v_0$};
      \filldraw (2) circle (0.05);
      \filldraw (3) circle (0.05);
      \filldraw (4) circle (0.05);
      \filldraw (6) circle (0.05) node[below right] {\footnotesize $v_n$};

      \path[freccia] (1) edge[ciclo] node[above]{\footnotesize $\scriptsize e_{00}$} (1);
			\path[freccia] (2) edge[ciclo] node[above]{\footnotesize $\scriptsize e_{11}$} (2);
			\path[freccia] (3) edge[ciclo] node[above]{\footnotesize $\scriptsize e_{22}$} (3);
			\path[freccia] (4) edge[ciclo] node[above]{\footnotesize $\scriptsize e_{33}$} (4);
      \path[freccia] (6) edge[ciclo] node[above] {\footnotesize $\scriptsize e_{nn}$} (6);

      \path[freccia] (1) edge node[above]{\footnotesize $\scriptsize e_{01}$} (2)
			               (2) edge node[above]{\footnotesize $\scriptsize e_{12}$} (3)
										 (3) edge node[above]{\footnotesize $\scriptsize e_{23}$} (4);
      \path[freccia,dashed] (4) edge (6);

      \path[freccia] (1) edge[bend right=30] node[below]{\footnotesize $\scriptsize e_{02}$} (3)
			               (1) edge[bend right=50] node[below, pos=0.4]{\footnotesize $\scriptsize e_{03}$} (4)
										 (2) edge[bend right=50] node[above]{\footnotesize $\scriptsize e_{13}$} (4);
	
	\path[freccia] (1) edge[bend right=60] node[above, pos=0.4]{\footnotesize $\scriptsize e_{0n}$} (6)
	(2)  edge[bend right=60] node[above, pos=0.39]{\footnotesize $\scriptsize e_{1n}$} (6)
	(3) edge[bend right=60] node[above, pos=0.4]{\footnotesize $\scriptsize e_{2n}$} (6)
	(4) edge[bend right=60] node[above, pos=0.4]{\footnotesize $\scriptsize e_{3n}$} (6);
	\end{tikzpicture}.
\end{equation}
The above graph has vertices $v_0,\ldots,v_n$ and an edge $e_{ij}:v_i\to v_j$ for every $i\leq j$. Note that all vertices here are ranges of edges, so the graph C*-
algebra is generated by the partial isometries corresponding to the edges. Next, by removing the loop $e_{nn}$ from the graph $\Lambda_{2n+1}$ we obtain the 
graph $\Gamma_{2n}$
\begin{equation}\label{eq:ballgraph}
\begin{tikzpicture}[baseline=(current bounding box.center)]

\clip (-0.6,-2.6) rectangle (9.7,1.8);

      \node[main node] (1) {};
      \node (2) [right of=1] {};
      \node (3) [right of=2] {};
      \node (4) [right of=3] {};
      \node (5) [right of=4] {};
      \node (6) [right of=5] {};

      \filldraw (1) circle (0.05) node[below left] {\footnotesize $\tilde{v}_0$};
      \filldraw (2) circle (0.05);
      \filldraw (3) circle (0.05);
      \filldraw (4) circle (0.05);
      \filldraw (6) circle (0.05) node[below right] {\footnotesize $\tilde{v}_n$};

      \path[freccia] (1) edge[ciclo] node[above]{\footnotesize $\tilde{e}_{00}$} (1);
			\path[freccia] (2) edge[ciclo] node[above]{\footnotesize $\tilde{e}_{11}$} (2);
			\path[freccia] (3) edge[ciclo] node[above]{\footnotesize $\tilde{e}_{22}$} (3);
			\path[freccia] (4) edge[ciclo] node[above]{\footnotesize $\tilde{e}_{33}$} (4);

      \path[freccia] (1) edge node[above]{\footnotesize $\tilde{e}_{01}$} (2)
			               (2) edge node[above]{\footnotesize $\tilde{e}_{12}$} (3)
										 (3) edge node[above]{\footnotesize $\tilde{e}_{23}$} (4);
      \path[freccia,dashed] (4) edge (6);

      \path[freccia] (1) edge[bend right=30] node[below]{\footnotesize $\tilde{e}_{02}$} (3)
			               (1) edge[bend right=50] node[below, pos=0.4]{\footnotesize $\tilde{e}_{03}$} (4)
										 (2) edge[bend right=50] node[above,pos=0.4]{\footnotesize $\tilde{e}_{13}$} (4);
	
	\path[freccia] (1) edge[bend right=60] node[above, pos=0.4]{\footnotesize $\tilde{e}_{0n}$} (6)
	(2)  edge[bend right=60] node[above, pos=0.4]{\footnotesize $\tilde{e}_{1n}$} (6)
	(3) edge[bend right=60] node[above, pos=0.4]{\footnotesize $\tilde{e}_{2n}$} (6)
	(4) edge[bend right=60] node[above, pos=0.4]{\footnotesize $\tilde{e}_{3n}$} (6);
\end{tikzpicture}
\end{equation}
whose graph C*-algebra we identify with $C(\mathsf{B}_q^{2n})$. (It is well known that $C(\mathsf{B}^2_q)\cong\mathcal{T}$, e.g., see \cite{r-i05}.) 
Furthermore, if $n\geq 1$, by removing $\tilde{v}_n$ and $r^{-1}(\tilde{v}_n)$ from $\Gamma_{2n}$, we obtain the graph $\Lambda_{2n-1}$ of 
$C(\mathsf{S}^{2n-1}_q)$.
We can view $\mathsf{B}^{2n}_q$ has a half equator of $\mathsf{S}^{2n+1}_q$ and $\mathsf{S}^{2n-1}_q$ as the boundary of $\mathsf{B}^{2n}_q$.
These cofibrations are defined by the surjective *-homomorphisms
\begin{equation}\label{eq:partialn}
\begin{aligned}
&r_n:C(\mathsf{S}^{2n+1}_q)\longrightarrow C(\mathsf{B}^{2n}_q) , & \qquad\quad &
\beta_n:C(\mathsf{B}^{2n}_q)\longrightarrow C(\mathsf{S}^{2n-1}_q) \;, \\
& r_n(S_{e_{ij}}):=\begin{cases}
S_{\tilde e_{ij}} & \text{if }(i,j)\neq (n,n), \\
P_{\tilde v_n} & \text{if }(i,j)=(n,n), 
\end{cases} &&
\beta_n(S_{\tilde e_{ij}}):=\begin{cases}
S_{e_{ij}} & \text{if }j\neq n, \\
0 & \text{if } j=n.
\end{cases}
\end{aligned}
\end{equation}
Observe that $r_n$ is not gauge equivariant, while
 $\beta_n$ and the composition $\beta_n\circ r_n$ are.

Following the classical intuition concerning balls and spheres, we can assemble the following ${U}(1)$-equivariant pullback diagram \cite{adht17}: 
\begin{equation}\label{eq:VSspheresPB}
\begin{tikzpicture}[baseline={([yshift=-3]current bounding box.center)},->,inner sep=2pt]

\node (A) at (2,4) { $C(\mathsf{S}^{2n+1}_q)_\bullet$ };
\node (B) at (0,2) { $C(\mathsf{S}^{2n-1}_q)_\bullet$ };
\node (C) at (4,2) { $C(\mathsf{B}^{2n}_q)\otimes C({S}^1)_\bullet$ };
\node (D) at (2,0) { $C(\mathsf{S}^{2n-1}_q)\otimes C({S}^1)_\bullet$ };

\path[-To,font=\scriptsize] 
	(A) edge[->>] node[above left] {$\beta_n\circ r_n$} (B)
	(A) edge node[above right] {$(r_n\otimes\id)\circ\delta$} (C)
	(B) edge node[below left] {$\delta$} (D)
	(C) edge[->>] node[below right] {$\beta_n\otimes\id$} (D);

\end{tikzpicture}.
\end{equation}
Here $\delta$ is the coaction dual to the ${U}(1)$-gauge action.
Passing to the subalgebras of ${U}(1)$-invariants and induced *-homomorphisms, we get the following pullback diagram \cite{adht17}:
\begin{equation}\label{eq:CWCPnq}
\begin{tikzpicture}[baseline={([yshift=-3]current bounding box.center)},->,inner sep=2pt]

\node (A) at (2,4) { $C(\CPq_q^n)$ };
\node (B) at (0,2) { $C(\CPq_q^{n-1})$ };
\node (C) at (4,2) { $C(\mathsf{B}^{2n}_q)$ };
\node (D) at (2,0) { $C(\mathsf{S}^{2n-1}_q)$ };

\path[font=\scriptsize] 
	(A) edge[->>] node[above left] {$(\beta_n\circ r_n)^{U(1)}$\!\!} (B)
	(A) edge (C)
	(B) edge (D)
	(C) edge[->>] node[below right] {$\beta_n$} (D);
			
\end{tikzpicture}.
\end{equation}

The cofibration
$\mathsf{S}^{2n-1}_q\hookrightarrow\mathsf{B}^{2n}_q$ dual to the *-homomorphism $\beta_n$ is a K-boundary map from a K-sphere to a K-ball in the 
sense of Definition \ref{def:Ksphereandball}. To see this, observe that $\mathsf{B}^{2n}_q$ is K-contractible (Example~\ref{ex:trivialK}),
$K^0(\mathsf{S}^{2n-1}_q)\cong\Z$ is generated by the class of the unit element (proved in \cite[Section~4.3]{HL04} using index pairing) and 
$K^1(\mathsf{S}^{2n-1}_q)\cong\Z$.
The unital \mbox{*-homomorphism} $\beta_n:C(\mathsf{B}^{2n}_q)\to C(\mathsf{S}^{2n-1}_q)$ induces then an isomorphism in $K_0$ as in \eqref{0even}. Finally, 
every graph C*-algebras is in the bootstrap class \cite[Theorem 5.5]{Kumjian2000}, which proves that $\mathsf{B}^{2n}_q$ and $\mathsf{S}^{2n-1}_q$
 belong to $\CQS^\times$.
In the category \CQS, the diagram dual to \eqref{eq:CWCPnq} is the desired pushout diagram of attaching a cell.
We thus obtain the following strict filtration by skeleta:
\begin{equation}\label{VSskel}
\star=\CPq_q^0\longhookrightarrow \CPq^1_q\longhookrightarrow \CPq^2_q\longhookrightarrow \ldots\longhookrightarrow \CPq_q^n \;.
\end{equation}

\begin{rmk}
In \cite[Section~5.3]{CHT21}, there is an alternative version of the pullback diagram~\eqref{eq:CWCPnq}. All the C*-algebras are the same and the
boundary quotient  *-homomorphism $\beta_n$ is the same, 
but the bottom injective *-homomorphism (the attaching morphism) is different because its image is not the $U(1)$-invariant
subalgebra of $C(\mathsf S^{2n-1}_q)$ as it is in our case. Consequently,  the $C(\CPq^n_q)\to C(\CPq^{n-1}_q)$ surjection,
which in \cite[Section~5.3]{CHT21} comes directly from an embedding of graphs, might be different. 
We thus obtain another strict quantum CW-complex structure of~$\CPq^n_q$.
\end{rmk}

\subsection{Weak quantum CW-complex structure of \texorpdfstring{$\CPq_{H}^{n}$}{CPHn}}\label{sec:weakCW}
Inspired by the Heegaard splitting of ${S}^3$, the C*-algebra $C(\mathsf{S}^3_H)$ of the Heegard quantum 3-sphere $\mathsf{S}^3_H$ was defined in \cite{BHMS04} by a suitable pullback construction. This was generalized to the multipullback C*-algebra $C(\mathsf{S}^{2n+1}_H)$ construction defining in~\cite{hnpsz18}.
Herein, we will use a presentation of $C(\mathsf{S}^{2n+1}_H)$ as the universal C*-algebra given by generators and relations that is afforded in \cite[Theorem~3.3]{hnpsz18}.
For $n\geq 0$, we identify $C(\mathsf{S}^{2n+1}_H)$ with the universal C*-algebra generated by isometries $s_0,\ldots,s_n$ subject to the commutation relations
\[
[s_i,s_j] =0=[s_i,s_j^*] ,
\]
for all $0\leq i\neq j\leq n$, and the relation
\begin{equation}\label{eq:sphererel}
\prod_{i=0}^n \big(1-s_is_i^*)=0 \,.
\end{equation}

Now, let us introduce the notation
\[
t_i:=\underbrace{1\otimes\ldots\otimes 1}_{i\text{ times}}{}\otimes s\otimes\underbrace{1\otimes\ldots\otimes 1}_{n-i\text{ times}} \;,\qquad i=0,\ldots,n,
\]
for the generators of $\mathcal{T}^{\otimes n+1}$, where $s$ is the isometry generating $\mathcal{T}$.
Right coactions of $C(S^1)$ on $C(\mathsf{S}^{2n+1}_H)$ and $\mathcal{T}^{\otimes n+1}$ dualize the diagonal $U(1)$-actions, and are all denoted by $\delta$. They are given on generators by
$\delta(s_i)=s_i\otimes u$ and $\delta(t_i)=t_i\otimes u$, respectively, for $0\leq i\leq n$.
For all $n\geq 1$, there is a short exact sequence of $U(1)$-C*-algebras \cite[(5.5)]{hnpsz18}:
\begin{equation}\label{eq:Kbound}
0\longrightarrow\mathcal{K}^{\otimes n}\longrightarrow\mathcal{T}^{\otimes n}\stackrel{\sigma_n}{\longrightarrow}C(\mathsf{S}^{2n-1}_H)\longrightarrow 0 ,
\end{equation}
where $\sigma_n$ is defined explicitly on generators by $\sigma_n(t_i):=s_i$, for all $0\leq i<n$.

The cofibration
$\mathsf S^{2n-1}_H\to \mathsf D_q^{\times n}$ dual to the *-homomorphim $\sigma_n$ is a K-boundary map from a K-sphere to a K-ball in the sense of Definition~\ref{def:Ksphereandball}. Indeed, since $\mathcal{K}^{\otimes n},\mathcal{T}^{\otimes n}\in\mathcal{N}_{\max}$, from 
\eqref{eq:Kbound} and \cite[Lemma 4.4(2)]{Uuy11} it follows that $C(\mathsf S^{2n-1}_H)\in\mathcal{N}_{\max}$.
The K{\"u}nneth formula yields $K_0(\mathcal{T}^{\otimes n})=\Z[1]$ and $K_1(\mathcal{T}^{\otimes n})=0$, 
while \cite[Theorem 5.2]{hnpsz18} states that $K_0(C(\mathsf S^{2n-1}_H))=\Z[1]$ and $K_1(C(\mathsf S^{2n-1}_H))\cong\Z$. Finally, since $\sigma_n$ is a unital *-homomorphism, it automatically follows that it induces an isomorphism in K-theory, so \eqref{0even} is satisfied.

We are now ready to recall:

\begin{lemma}[{\protect\cite[Lemma 6.2]{hnpsz18}}]\label{lemma:3.1}
For every $n\geq 1$, there is a
${U}(1)$-equivariant pullback diagram
\begin{equation}\label{eq:Mpull}
\begin{tikzpicture}[baseline={([yshift=-4pt]current bounding box.center)}]

\node (A) at (2,3.6) { $C(\mathsf{S}^{2n+1}_H)_\bullet$ };
\node (B) at (0,1.8) { $C(\mathsf{S}^{2n-1}_H)_\bullet\otimes\mathcal{T}_\bullet$ };
\node (C) at (4,1.8) { $\mathcal{T}^{\otimes n}_\bullet\otimes C({S}^1)_\bullet$ };
\node (D) at (2,0) { $C(\mathsf{S}^{2n-1}_H)_\bullet\otimes C({S}^1)_\bullet$ };

\path[->>,font=\scriptsize] 
      (A) edge node[above left] {$p_1$} (B)
			(A) edge node[above right] {$p_2$} (C)
			(B) edge node[below left] {$\pi_1$} (D)
			(C) edge node[below right] {$\pi_2$} (D);
			
\end{tikzpicture}.
\end{equation}
Here $\pi_1:=\id\otimes\sigma_1$, $\pi_2:=\sigma_n\otimes\id$, and the remaining two maps are defined in terms of generators by:
\begin{align*}
p_1(s_i)  &:=\bigg\{\!\begin{array}{ll}
s_i\otimes 1 & \forall\;i=0,\ldots,n-1, \\
1\otimes t & \text{if }i=n,
\end{array}
\\
p_2(s_i)  &:=\bigg\{\!\begin{array}{ll}
t_i\otimes 1 &\forall\;i=0,\ldots,n-1, \\
1\otimes u &\text{if }i=n.
\end{array}
\end{align*}
\end{lemma}

Next, we need another technical lemma, which is an easy consequence of Lemma \ref{lemma:3.1}.

\begin{lemma}\label{lemma:p1k}
For all $n,k\geq 0$, we have a ${U}(1)$-equivariant pullback diagram:
\begin{equation}\label{eq:MpullT}
\begin{tikzpicture}[baseline={([yshift=-4pt]current bounding box.center)}]

\node (A) at (2.5,4) { $C(\mathsf{S}^{2n+1}_H)_\bullet\otimes\mathcal{T}^{\otimes k}_\bullet$ };
\node (B) at (0,2) { $C(\mathsf{S}^{2n-1}_H)_\bullet\otimes\mathcal{T}^{\otimes k+1}_\bullet$ };
\node (C) at (5,2) { $\mathcal{T}^{\otimes n}\otimes C({S}^1)_\bullet\otimes\mathcal{T}^{\otimes k}$ };
\node (D) at (2.5,0) { $C(\mathsf{S}^{2n-1}_H)\otimes C({S}^1)_\bullet\otimes\mathcal{T}^{\otimes k}$ };

\path[-To,font=\scriptsize] 
      (A) edge node[above left] {$p_1^k$} (B)
			(A) edge node[above right] {$p_2^k$} (C)
			(B) edge node[below left] {$\pi_1^k$} (D)
			(C) edge node[below right] {$\pi_2^k$} (D);
			
\end{tikzpicture}.
\end{equation}
Here
$p_1^k:=p_1\otimes\id_{\mathcal{T}^{\otimes k}}$,
$\pi_1^k:=\phi_1\circ (\pi_1\otimes\id_{\mathcal{T}^{\otimes k}})$,
$\pi_2^k:=\sigma_n\otimes\id_{C({S}^1)}\otimes\id_{\mathcal{T}^k}$,
$p_2^k:=\phi_2\circ (p_2\otimes\id_{\mathcal{T}^{\otimes k}})$,
where $p_1,p_2,\pi_1,\pi_2$ are the *-homomorphisms as in Lemma~\ref{lemma:3.1}, and the *-isomorphisms
of $U(1)$-C*-algebras
\begin{alignat*}{2}
\phi_1: && \; C(\mathsf{S}^{2n-1}_H)_\bullet\otimes C({S}^1)_\bullet\otimes\mathcal{T}^{\otimes k}_\bullet &\longrightarrow C(\mathsf{S}^{2n-1}_H)\otimes C({S}^1)_\bullet\otimes\mathcal{T}^{\otimes k} , \\
\phi_2: && \mathcal{T}^{\otimes n}_\bullet\otimes C({S}^1)_\bullet\otimes\mathcal{T}^{\otimes k}_\bullet &\longrightarrow\mathcal{T}^{\otimes n}\otimes C({S}^1)_\bullet\otimes\mathcal{T}^{\otimes k} ,
\end{alignat*}
are special cases of \eqref{eq:gaugetrick}.
\end{lemma}
\begin{proof}
By tensoring all algebras in \eqref{eq:Mpull} with $\mathcal{T}^{\otimes k}$ and all maps with the identity, we get a ${U}(1)$-equivariant pullback diagram
\begin{equation}\label{eq:5.4b}
\begin{tikzpicture}[baseline=(current bounding box.center)]

\node (A) at (2.5,4) { $C(\mathsf{S}^{2n+1}_H)_\bullet\otimes\mathcal{T}^{\otimes k}_\bullet$ };
\node (B) at (0,2) { $C(\mathsf{S}^{2n-1}_H)_\bullet\otimes\mathcal{T}^{\otimes k+1}_\bullet$ };
\node (C) at (5,2) { $\mathcal{T}^{\otimes n}_\bullet\otimes C({S}^1)_\bullet\otimes\mathcal{T}^{\otimes k}_\bullet$ };
\node (D) at (2.5,0) { $C(\mathsf{S}^{2n-1}_H)_\bullet\otimes C({S}^1)_\bullet\otimes\mathcal{T}^{\otimes k}_\bullet$ };

\path[-To,font=\scriptsize] 
      (A) edge node[above left] {\smash{$p_1\otimes\id_{\mathcal{T}^{\otimes k}}$}} (B)
			(A) edge node[above right] {\smash{$p_2\otimes\id_{\mathcal{T}^{\otimes k}}$}} (C)
			(B) edge node[below left] {\smash{$\pi_1\otimes\id_{\mathcal{T}^{\otimes k}}$}} (D)
			(C) edge node[below right] {\smash{$\pi_2\otimes\id_{\mathcal{T}^{\otimes k}}$}} (D);
			
\end{tikzpicture},
\end{equation}
where $p_1,p_2,\pi_1,\pi_2$ are the maps as in Lemma~\ref{lemma:3.1}.
Using the isomorphisms $\phi_1$ and $\phi_2$ of ${U}(1)$-C*-algebras,
we transform the diagram \eqref{eq:5.4b} into \eqref{eq:MpullT}. Finally, note that
\[
\pi_2^k =\phi_1\circ (\pi_2\otimes\id_{\mathcal{T}^{\otimes k}})\circ\phi_2^{-1}=\pi_2\otimes\id_{\mathcal{T}^{\otimes k}}
\]
because of the equivariance of $\pi_2$.
\end{proof}

We now establish the main technical result of this section, i.e.~the K-equivalence between tubular neighbourhoods.

\begin{lemma}[Tubular Neighbourhood Lemma]\label{lemma:3.3}
For all $n,k\in \N$, the $U(1)$-equivariant map
\[
\id\otimes 1_{\mathcal{T}}:C(\mathsf S^{2n+1}_H)_\bullet\otimes\mathcal{T}^{\otimes k}_\bullet \to C(\mathsf S^{2n+1}_H)_\bullet\otimes\mathcal{T}^{\otimes k+1}_\bullet
\]
induces a K-equivalence between the $U(1)$-fixed points subalgebras
\begin{equation}\label{eq:weqtub}
C(\tub{n}{n+k})\weq C(\tub{n}{n+k+1}).
\end{equation}
\end{lemma}

\begin{proof}
Let us consider first the following commutative diagram:
\begin{equation*}
\begin{tikzpicture}[yscale=2,xscale=1.8,baseline={([yshift=-4pt]current bounding box.center)}]

\node (A) at (1.2,2) { $C(\mathsf S^{2n+1}_H)\otimes\mathcal{T}^k$ };
\node (B) at (0,1) { $C(\mathsf S^{2n-1}_H)\otimes\mathcal{T}^{k+1}$ };
\node (C) at (2.4,1) { $\mathcal{T}^n\otimes C(S^1)\otimes\mathcal{T}^k$ };
\node (D) at (1.2,0) { $C(\mathsf S^{2n-1}_H)\otimes C(S^1)\otimes\mathcal{T}^k$ };
\node (E) at (5.7,2) { $C(\mathsf S^{2n+1}_H)\otimes\mathcal{T}^{k+1}$ };
\node (F) at (4.5,1) { $C(\mathsf S^{2n-1}_H)\otimes\mathcal{T}^{k+2}$ };
\node (G) at (6.9,1) { $\mathcal{T}^n\otimes C(S^1)\otimes\mathcal{T}^{k+1}$ };
\node (H) at (5.7,0) { $C(\mathsf S^{2n-1}_H)\otimes C(S^1) \otimes\mathcal{T}^{k+1}$ };

\path[->] (A) edge (B)
			(A) edge (C)
			(B) edge (D)
			(C) edge (D)
			(E) edge (F)
			(E) edge (G)
			(F) edge (H)
			(G) edge (H)
			(A) edge[bend left=10] (E)
			(B) edge[bend right=20] (F)
			(C) edge[bend left=20] (G)
			(D) edge[bend right=10] (H);
			
\end{tikzpicture}\!.
\end{equation*}
The left square is \eqref{eq:MpullT}, the right square is \eqref{eq:MpullT} with $k$ replaced by $k+1$,
and the horizontal arrows are all given by $\id\otimes 1_{\mathcal{T}}$.
Passing to the fixed-point subalgebras we get the commutative diagram
\begin{equation}\label{eq:induction}
\begin{tikzpicture}[yscale=2,xscale=1.8,baseline={([yshift=-15pt]current bounding box.center)}]

\node (A) at (1.2,2) { $\big(C(\mathsf S^{2n+1}_H)\otimes\mathcal{T}^k\big)^{U(1)}$ };
\node (B) at (0,1) { $\big(C(\mathsf S^{2n-1}_H)\otimes\mathcal{T}^{k+1}\big)^{U(1)}$ };
\node (C) at (2.4,1) { $\mathcal{T}^n\otimes\mathcal{T}^k$ };
\node (D) at (1.2,0) { $C(\mathsf S^{2n-1}_H)\otimes \mathcal{T}^k$ };
\node (E) at (5.7,2) { $\big(C(\mathsf S^{2n+1}_H)\otimes\mathcal{T}^{k+1}\big)^{U(1)}$ };
\node (F) at (4.5,1) { $\big(C(\mathsf S^{2n-1}_H)\otimes\mathcal{T}^{k+2}\big)^{U(1)}$ };
\node (G) at (6.9,1) { $\mathcal{T}^n\otimes \mathcal{T}^{k+1}$ };
\node (H) at (5.7,0) { $C(\mathsf S^{2n-1}_H)\otimes \mathcal{T}^{k+1}$ };

\path[->,font=\footnotesize] (A) edge (B)
			(A) edge (C)
			(B) edge (D)
			(C) edge (D)
			(E) edge (F)
			(E) edge (G)
			(F) edge (H)
			(G) edge (H)
			(A) edge[bend left=10] node[above] {$\phi_{n,k}$} (E)
			(B) edge[bend right=20] node[below,pos=0.6] {$\phi_{n-1,k+1}$} (F)
			(C) edge[bend left=20] (G)
			(D) edge[bend right=10] (H);
			
\end{tikzpicture}\!,
\end{equation}
where the two squares are still pullback diagrams.
We now prove by induction on $n\geq 0$ that, for all $k\geq 0$, the map $\phi_{n,k}$ is a K-equivalence.

Let us start with $n=0$, and consider the $U(1)$-equivariant commutative diagram
\begin{center}
\begin{tikzpicture}[scale=2]

\node (A) at (0,1) {$C(S^1)_\bullet\otimes\mathcal{T}^k_\bullet$};
\node (B) at (2.5,1) {$C(S^1)_\bullet\otimes\mathcal{T}^{k+1}_\bullet\phantom{\,,}$};
\node (C) at (0,0) {$C(S^1)_\bullet\otimes\mathcal{T}^k$};
\node (D) at (2.5,0) {$C(S^1)_\bullet\otimes\mathcal{T}^{k+1}$\,,};

\path[->]
		(A) edge node[above]{\scriptsize $\id\otimes 1_{\mathcal{T}}$} (B)
		(A) edge (C)
		(B) edge (D)
		(C) edge node[above]{\scriptsize $\id\otimes 1_{\mathcal{T}}$} (D);

\end{tikzpicture}
\end{center}
where the vertical arrows are the isomorphisms $a\otimes b\mapsto ab_{(1)}\otimes b_{(0)}$ (cf.~\eqref{eq:gaugetrick}). Taking the $U(1)$-invariant subalgebras gives:
\begin{equation}\label{eq:nequalzero}
\begin{tikzpicture}[scale=2,baseline=(current bounding box.center)]

\node (A) at (0,1) {$\big(C(S^1)_\bullet\otimes\mathcal{T}^k_\bullet\big)^{U(1)}$};
\node (B) at (2.5,1) {$\big(C(S^1)_\bullet\otimes\mathcal{T}^{k+1}_\bullet\big)^{U(1)}$};
\node (C) at (0,0) {$\mathcal{T}^k$};
\node (D) at (2.5,0) {$\mathcal{T}^{k+1}$\,.\!\!};

\path[->]
		(A) edge node[above]{\scriptsize $\phi_{0,k}$} (B)
		(A) edge (C)
		(B) edge (D)
		(C) edge node[above]{\scriptsize $\id\otimes 1_{\mathcal{T}}$} (D);

\end{tikzpicture}
\end{equation}
The bottom arrow is a morphism between K-contractible quantum spaces, so it is a K-equivalence by
\ref{ax:contr}. Furthermore, since the vertical arrows are isomorphisms, $\phi_{0,k}$ is also a K-equivalence.

Let us assume now that $\phi_{n-1,k}$ is a K-equivalence for some $n\geq 1$ and for all $k\geq 0$, and look again at the diagram \eqref{eq:induction}. The  *-homomorphisms $C(\mathsf S^{2n-1}_H)\otimes\mathcal{T}^k\to C(\mathsf S^{2n-1}_H)\otimes\mathcal{T}^{k+1}$ and $\mathcal{T}^n\otimes\mathcal{T}^k\to\mathcal{T}^n\otimes\mathcal{T}^{k+1}$ are K-equivalences by \ref{ax:weqmon}
and \ref{ax:contr}, respectively. Finally, since $\phi_{n-1,k+1}$ is a K-equivalence by assumption, we infer from \ref{ax:weq} that the top arrow $\phi_{n,k}$ is also a K-equivalence, as needed.
\end{proof}

As the special case $k=0$, we conclude that:

\begin{cor}\label{cor:tub}
The *-homomorphism
\begin{equation}\label{eq:tubweq}
(\id\otimes 1_{\mathcal{T}})^{U(1)}:C(\CPq^{n-1}_H)\weq
C(\tub{n-1}{n})
\end{equation}
is a K-equivalence.
\end{cor}

We thus arrive at the pinnacle of this section:

\begin{thm}\label{weakCWmult}
The multipushout quantum complex projective space $\CPq^N_{H}$ has a weak filtration by skeleta
\begin{equation}\label{skel}
\CPq_{H}^{0}\rightarrowtail\CPq_{H}^{1}\rightarrowtail\cdots\rightarrowtail 
\CPq_{H}^{N-1}\rightarrowtail \CPq_{H}^{N}\,.
\end{equation}
Here, at every step, a single quantum cell is weakly attached:
\begin{equation}\label{eq:weakatt}
\begin{tikzpicture}[scale=1.5,baseline={([yshift=-3pt]current bounding box.center)}]

\node (A) at (1,2) {$\CPq_{H}^{n}$};
\node (B) at (0,1) {$\CPq_{H}^{n-1}$};
\node (C) at (2,1) {$\mathsf{D}_q^{\times n}$};
\node (D) at (1,0) {$\mathsf{S}_{H}^{2n-1}$};

\path[->]
     	(B) edge[>->](A)
	(C) edge (A)
	(D) edge  (B)
	(D) edge[right hook->] (C);

\end{tikzpicture}.
\end{equation}
\end{thm}
\begin{proof}
Passing from the pullback diagram \eqref{eq:MpullT} to the pullback diagram of the $U(1)$-invariant subalgebras, for $k=0$, we obtain the pullback diagram 
defining the bottom pushout square of the following commutative diagram:
\begin{equation}\label{twosquares}
\begin{tikzpicture}[scale=2,baseline=(current bounding box.center)]

\node (A) at (1,3) {$\widetilde{\CPq_{H}^{n}}$};
\node (B) at (0,2) {$\CPq_{H}^{n-1}$};
\node (C) at (2,2) {$\CPq_{H}^{n}$};
\node (D) at (1,1) {$\tub{n-1}{n}$};
\node (E) at (2,0) {$\mathsf{S}_{H}^{2n-1}$};
\node (F) at (3,1) {$\mathsf{D}_q^{\times n}$\,.\!};

\path[->,font=\small]
		(B) edge[right hook->,dashed] (A)
		(C) edge[dashed] node[above,sloped]{$\sim$} (A)
		(D) edge node[above,sloped] {$\sim$} (B)
		(E) edge[left hook->] (D)
		(F) edge[left hook->] (C)
		(D) edge[right hook->] (C)
		(E) edge[right hook->] (F);

\end{tikzpicture}
\end{equation}
Here, the K-equivalence $\tub{n-1}{n}\weq\CPq_H^{n-1}$ is defined by \eqref{eq:tubweq}, $\widetilde{\CPq_{H}^{n}}$ is defined as the pushout of the span $\CPq_{H}^{n-1}\stackrel{\text{\raisebox{-3pt}{\scriptsize $\sim$}}}{\longleftarrow}\tub{n-1}{n}\cof\CPq_{H}^{n}$, and we used \ref{ax:cof} and \ref{ax:cwW} to deduce that the two morphisms to $\widetilde{\CPq_{H}^{n}}$ are a cofibration and a K-equivalence, respectively. The outer rectangle is a pushout diagram because both the squares are pushout diagrams. The top (dashed) roof defines the weak cofibration $\CPq_{H}^{n-1}\rightarrowtail \CPq_{H}^{n}$
in \eqref{skel}, and fits into the weak pushout
\eqref{eq:weakatt}
weakly attaching a quantum polydisk to $\CPq^{n-1}_H$.
\end{proof}

Combining Theorem \ref{weakCWmult} with Corollary \ref{KweakK}, we obtain the following result.

\begin{cor}
$K^{0}(\CPq_{H}^{n}) \cong \Z^{n+1}$ and $K^{1}(\CPq_{H}^{n}) \cong 0$.
\end{cor}

We end this section by comparing our computation of K-groups of the multipushout quantum complex projective spaces with the original computation 
performed in \cite[Theorem~5.6]{hnpsz18}. Therein, 
the starting point is the decomposition of the odd spheres into disjoint components
(see \cite[p.~847]{hnpsz18}):
\begin{equation}\label{nest}
{S}^{2n+1}\;\;\cong\;\;{B}^{2n}\!\times{S}^1\ \coprod\ {S}^{2n-1}\!\times{D}_0\,,
\end{equation}
where $D_0$ is an open disc. The $U(1)$ quotient of this decomposition yields the decomposition of quotient spaces
\begin{equation}\label{nest2}
\CP^{n}\;\;\cong\;\;{B}^{2n}\ \coprod\ {S}^{2n-1}\!\times^{U(1)}{D}_0\,
\end{equation}
captured in the following short exact sequence of commutative C*-algebras:
\begin{equation}\label{nest2}
0\longrightarrow \big(C({S}^{2n-1})\otimes C_0({D}_0)\big)^{U(1)}\longrightarrow C(\CP^{n}) \longrightarrow C({B}^{2n}) \longrightarrow 0.
\end{equation}
The Heegaard quantum incarnation of the above short exact sequence is the short exact sequence
\begin{equation}\label{nest3}
0\longrightarrow \big(C(\mathsf S^{2n-1}_H)\otimes \mathcal{K}\big)^{U(1)}\longrightarrow C(\CPq^{n}_H) \longrightarrow \mathcal{T}^{\otimes n} \longrightarrow 0,
\end{equation}
which is a key idea in the proof of \cite[Theorem~5.6]{hnpsz18}.

A key difference between this proof and our proof is that our starting point is the pushout decomposition of the odd spheres:
\begin{equation}\label{our}
{S}^{2n+1}\;\;\cong\;\;{B}^{2n}\!\times{S}^1\ \coprod_{{S}^{2n-1}\times{S}^1}\ {S}^{2n-1}\!\times{D}.
\end{equation}
Thus, instead of decomposing $S^{2n+1}$ into a closed and an open subset that are disjoint, we decompose it into two closed subsets that intersect non-
trivialy. Now, the $U(1)$ quotient of this decomposition yields the pushout decomposition of quotient spaces
\begin{equation}\label{our2}
\CP^{n}\;\;\cong\;\;{B}^{2n}\ \coprod_{{S}^{2n-1}}\ {S}^{2n-1}\!\times^{U(1)}{D}\,,
\end{equation}
which can be recast as the pushout diagram
\begin{equation}\label{classpull2}
\begin{tikzpicture}[<-,scale=1.6,baseline={(current bounding box.center)}]

\node (A) at (1,1) {$\CP^n$};
\node (B) at (0,0) {${S}^{2n-1}\!\times^{U(1)}{D}$};
\node (C) at (2,0) {$B^{2n}$};
\node (D) at (1,-1) {$S^{2n-1}$};

\path
	(A) edge[<-right hook] (B)
	(A) edge (C)
	(B) edge (D)
	(C) edge[<-right hook] (D);
			
\end{tikzpicture}\,.
\end{equation}
The Heegaard quantum incarnation of the above pushout diagram is the pullback diagram
\begin{equation}
\begin{tikzpicture}[baseline={([yshift=-4pt]current bounding box.center)},scale=1.6]

\node (A) at (1,2) { $C(\CPq^{n}_H)$ };
\node (B) at (0,1) { $\big(C(\mathsf{S}^{2n-1}_H)\otimes\mathcal{T}\big)^{U(1)}$ };
\node (C) at (2,1) { $\mathcal{T}^{\otimes n}$ };
\node (D) at (1,0) { $C(\mathsf{S}^{2n-1}_H)$ };

\path[->,font=\scriptsize] 
      (A) edge node[above left] {} (B)
			(A) edge node[above right] {} (C)
			(B) edge node[below left] {} (D)
			(C) edge node[below right] {} (D);
			
\end{tikzpicture},
\end{equation}
which is a key ingredient in our computation of K-groups.

\section{K-equivalences of Vaksman--Soibelman and multipushout quantum spaces}\label{sec:Ktheorytype}
\noindent
The goal of this section is to prove K-equivalences between the Vaksman--Soibelman and the multipushout sides, and then harvest resulting benefits.
As the Vaksman--Soibelman side is much more computable due to its graph-algebraic presentation, we gain a lot even when an independent proof exists
on the multipushout side. First, we construct a $U(1)$-equivariant K-equivalence of quantum balls and spheres. The first benefit is an explicit formula for
a generator of $K_1(C(\mathsf S^{2n+1}_H))$, which we do not know how to compute otherwise. 

Then we prove the   weak cellular equivalence
of quantum projective spaces. On the other hand, we combine  the aforementioned $U(1)$-equivariant map between quantum spheres with the
 theorem in \cite{HM16},  stating that the association of finitely generated projective modules commutes 
with their pushforwards, to transport the associated line bundles over the Vaksman--Soibelman quantum projective spaces to 
the associated line bundles over the multipushout quantum projective spaces. Now we can claim the desired result that associated line bundles
provide free generators of $K_0(C(\CPq^n_H))$, and thus extend the $n=2$ case proven in~\cite{d2025milnor}. 

Better still, we 
use identities obtained for line bundles over the Vaksman--Soibelman quantum projective spaces using graph-algebraic techniques \cite{d2026k} to construct
the isomorphism of abelian groups $\varphi\colon K^0(\CP^n)\to K^0(\CPq^n_q)$ mapping the class of the line bundle with  the winding number $k$
to the class of the line bundle with the same winding number~$k$. This allows us to transport the ring structure from $K^0(\CP^n)$ to~$K^0(\CPq^n_q)$. 
Subsequently, we use the K-equivalence of $\CPq^n_q$ and  $\CPq^n_H$ to transport the ring structure to~$K^0(\CPq^n_H)$. Much as in the classical
case, we can infer from the ring structure $K_0$-identities expressing the classes of line bundles in terms of free generators of the $K_0$-group.
At least some of these identities can also be derived  (in a more complicated way) using groupoids for $\CPq^n_q$ (see~\cite{s-aj2019}), 
and projections in the C*-algebra~$\mathcal{T}^{\otimes n}$ for $\CPq^n_H$  (see~\cite{s-aj}).

\subsection{The K-equivalence of 
\texorpdfstring{$\mathsf{S}^{2n+1}_q$}{Sq}
and \texorpdfstring{$\mathsf{S}^{2n+1}_H$}{SH}}\label{sec:avb2}

The main focus of this section will be the following diagram ($n\geq 1$):
\begin{equation}\label{eq:ballpullback}
\begin{tikzpicture}[baseline=(current bounding box.center),yscale=2.2,xscale=3.5]

\node (a) at (0,1) {$C(\mathsf{B}^{2n}_q)$};
\node (b) at (1,1) {$\mathcal{T}^{\otimes n}$};
\node (c) at (0,0) {$C(\mathsf{S}^{2n-1}_q)$};
\node (d) at (1,0) {$C(\mathsf{S}^{2n-1}_H)$\,.\!\!};

\begin{scope}[font=\footnotesize]
\path
	(a) edge[right hook->] node[above] {$\rho_n$} (b)
	(c) edge[right hook->] node[above] {$\omega_{n-1}$} (d)
	(a) edge[->>] node[left] {$\beta_n$} (c)
	(b) edge[->>] node[right] {$\sigma_n$} (d);
\end{scope}

\end{tikzpicture}
\end{equation}
Here, the vertical *-homomorphisms $\beta_n$ and $\sigma_n$ are defined in \eqref{eq:partialn} and \eqref{eq:Kbound}, respectively.
We start by defining the *-homomorphisms $\rho_n$ and $\omega_{n-1}$, and then we show that \eqref{eq:ballpullback} is a pullback diagram. Finally, we conclude that $\rho_n$ and $\omega_{n-1}$ are K-equivalences.

\begin{lemma}\label{prop:ball2ball}
The following formulas
\begin{gather*}
\rho_n(S_{\tilde{e}_{ij}}) := t_it_jt_j^*\prod_{k=0}^{j-1}(1-t_kt_k^*), \quad
\rho_n(P_{\tilde{v}_j}) :=t_jt_j^*\prod_{k=0}^{j-1}(1-t_kt_k^*) ,
\quad \forall\;0\leq i\leq j< n, \\
\rho_n(S_{\tilde{e}_{in}}) := t_i \prod_{k=0}^{n-1}(1-t_kt_k^*), \quad
\rho_n(P_{\tilde{v}_n}) :=\prod_{k=0}^{n-1}(1-t_kt_k^*) , 
\quad \forall\;0\leq i<n, \\
\omega_{n-1}(S_{e_{ij}}) :=s_is_js_j^*\prod_{k=0}^{j-1}(1-s_ks_k^*), \quad
\omega_{n-1}(P_{v_j}) :=s_js_j^*\prod_{k=0}^{j-1}(1-s_ks_k^*) ,
\quad \forall\;0\leq i\leq j<n,
\end{gather*}
define two injective $U(1)$-equivariant unital *-homomorphisms
\[
\rho_n:C(\mathsf{B}^{2n}_q)\longrightarrow\mathcal{T}^{\otimes n}
\qquad\text{and}\qquad
\omega_{n-1}:C(\mathsf{S}^{2n-1}_q)\longrightarrow C(\mathsf{S}^{2n-1}_H) .
\]
Moreover, with $\rho_n$ and $\omega_{n-1}$ defined as above, \eqref{eq:ballpullback}  is a commutative diagram.
\end{lemma}

\begin{proof}
We begin by noting that, since \mbox{$(1-t_kt_k^*)t_kt_k^*=0$}, the projections $\rho_n(P_{\tilde{v}_j})$ are mutually orthogonal.
Next, we check that $\rho_n$ preserves the Cuntz--Krieger relations.
The first relation
$\rho_n(S_{\tilde{e}_{ij}})^*\rho_n(S_{\tilde{e}_{ij}})=\rho_n(P_{\tilde{v}_j})$ easily follows from $t_i$'s being commuting isometries.
First, to verify the second Cuntz--Krieger relation,
we check by induction on $i$ from $n$ to lower values
that
\begin{equation}\label{eq:inductionrhoP}
\sum_{j=i}^n\rho_n(P_{\tilde{v}_j})=\prod_{k=0}^{i-1}(1-t_kt_k^*) .
\end{equation}
Combining \eqref{eq:inductionrhoP} with $\rho_n(S_{\tilde{e}_{ij}})\rho_n(S_{\tilde{e}_{ij}})^*=t_i \rho_n(P_{\tilde{v}_j})  t_i^*$, we obtained the desired equality
$$
\sum_{j=i}^n\rho_n(S_{\tilde{e}_{ij}})\rho_n(S_{\tilde{e}_{ij}})^*=
t_i\bigg(\sum_{j=i}^n\rho_n(P_{\tilde{v}_j})\bigg)t_i^*=t_i\prod_{k=0}^{i-1}(1-t_kt_k^*)t_i^*=t_it_i^*\prod_{k=0}^{i-1}(1-t_kt_k^*)=\rho_n(P_{\tilde{v}_i})
$$
for all $i\neq n$. Therefore, the *-homomorphism $\rho_n$ is well defined. Similarly, one proves that $\omega_{n-1}$ is well defined.

Furthermore, one explicitly checks on generators that the diagram \eqref{eq:ballpullback} is commutative. In particular,
$\sigma_n\big(\rho_n(S_{\tilde{e}_{in}})\big)=0$ because $\prod_{k=0}^{n-1}(1-s_ks_k^*)=0$ in $C(\mathsf{S}^{2n-1}_H)$.
Also, $\rho_n$ and $\omega_{n-1}$ are clearly $U(1)$-equivariant.
Their injectivity follows from the gauge-invariant uniqueness theorem \cite[Theorem~2.2]{r-i05}. Indeed,
$\rho_n(P_{\tilde v_i})\neq 0$ for all $i$, because every $t_i$ can be identified with the unilateral shift on $\ell^2(\N)$. Finally,  $\rho_n(P_{\tilde v_i})\notin\mathcal{K}^{\otimes n}$ because $t_it_i^*\notin\mathcal{K}$,
so from \eqref{eq:Kbound} we deduce that $\omega_{n-1}(P_{v_i})=\sigma_n(\rho_n(P_{\tilde v_i}))\neq 0$ for all $i$.
\end{proof}

To prove the pullback property of \eqref{eq:ballpullback}  we need the following technical lemma.

\begin{lemma}\label{lemma:kernels}
For all $n\geq 1$, $\ker(\sigma_n)\subseteq\mathrm{im}(\rho_n)$.
\end{lemma}

\begin{proof}
For starters, recall that $\ker(\sigma_n)=\mathcal{K}^{\otimes n}$ due to \eqref{eq:Kbound}.
For the sake of the present proof, we identify 
$\mathcal{T}$ with its image under the representation on $\ell^2(\N)$ sending the generating isometry to the unilateral shift.
Under this identification, $\mathcal{T}^{\otimes n}$ becomes a subalgebra of $\mathcal{B}(\ell^2(\N)^{\otimes n})$, and
$\mathcal{K}^{\otimes n}$ becomes the C*-algebra $\mathcal{K}(\ell^2(\N)^{\otimes n})$.
Furthermore, since
\[
\rho_n(P_{\tilde v_n})=\prod_{k=0}^{n-1}(1-t_kt_k^*)\in\mathcal{K}(\ell^2(\N)^{\otimes n})\setminus\{0\} ,
\]
to show the desired inclusion,
by \cite[Corollary I.10.4]{d-kr96}, it suffices to show that the commutant $\mathrm{im}(\rho_n)'$ of $\mathrm{im}(\rho_n)$ in $\mathcal{B}(\ell^2(\N)^{\otimes n})$ is trivial.

To this end, let us use \eqref{eq:inductionrhoP} and define
\[
x_i:=\sum_{j=i}^n\rho_n(S_{\tilde{e}_{ij}})=t_i\sum_{j=i}^n\rho_n(P_{\tilde{v}_j})=t_i\prod_{k=0}^{i-1}(1-t_kt_k^*)
=(1-ss^*)^{\otimes i}\otimes s \otimes 1^{\otimes (n-i-1)},
\]
for all $0\leq i<n$.
Next, take $a:=1^{\otimes k}\otimes b$ with $b\in\mathcal{B}(\ell^2(\N)^{\otimes (n-k)})$ and $0\leq k<n$, and compute
\[
[a,x_k]=\Big(\prod_{i=0}^{k-1}(1-t_it_i^*)\Big)[a,t_k] =
(1- ss^*)^{\otimes k}\otimes [b,s\otimes 1^{\otimes n-k-1}] \;.
\]
A similar formula holds for $[a,x_k^*]$.
On the other hand, a simple computation proves that $b\in\mathcal{B}(\ell^2(\N)^{\otimes (n-k)})$ commutes with 
$s\otimes 1^{\otimes n-k-1}$ and its adjoint if and only if $b=1\otimes c$ for some $c\in\mathcal{B}(\ell^2(\N)^{\otimes (n-k-1)})$.
Hence, we see that $a$ commutes with $x_k$ and $x_k^*$ if and only if $b$ commutes with $s\otimes 1^{\otimes n-k-1}$ and its adjoint, if and only if $a=1^{\otimes k+1}\otimes c$ for some $c\in\mathcal{B}(\ell^2(\N)^{\otimes (n-k-1)})$.
Now, to use induction we define 
\[
A_k:=1^{\otimes k}\otimes\mathcal{B}(\ell^2(\N)^{\otimes (n-k)})\subseteq\mathcal{B}(\ell^2(\N)^{\otimes n}) ,\qquad\forall\;0\leq k<n,
\qquad A_n:=\C 1.
\]
Clearly, $\mathrm{im}(\rho_n)'\subseteq A_0=\mathcal{B}(\ell^2(\N)^{\otimes n})$.
For the inductive step, suppose that $\mathrm{im}(\rho_n)'\subseteq A_k$ for some $0\leq k<n$, and take
$a\in\mathrm{im}(\rho_n)'$.
Finally, since $a$ commutes with $x_k$ and $x_k^*$, by the argument above, $a\in A_{k+1}$. Consequently, $\mathrm{im}(\rho_n)'\subseteq A_n=\C 1$ as desired.
\end{proof}

We are ready to claim the following key result.

\begin{thm}\label{prop:finalstatementball}
The diagram \eqref{eq:ballpullback} is a pullback diagram.
\end{thm}

\begin{proof}
Since $\rho_n$ and $\omega_{n-1}$ are injective,
%it is enough to prove that 
by \cite[Proposition~3.1]{Ped}, it suffices to prove that
\[
\mathrm{im}(\rho_n)=\sigma_n^{-1}\omega_{n-1}\big(C(\mathsf{S}^{2n-1}_q)\big) .
\]
As the inclusion 
$\mathrm{im}(\rho_n)\subseteq \sigma_n^{-1}\omega_{n-1}\big(C(\mathsf{S}^{2n-1}_q)\big)$ follows from the 
commutativity of the diagram \eqref{eq:ballpullback}, it only remains to prove the opposite inclusion.
To this end, take any $x\in \sigma_n^{-1}\omega_{n-1}\big(C(\mathsf{S}^{2n-1}_q)\big)$. 
Then $\sigma_n(x)=\omega_{n-1}(y)$ for some $y\in C(\mathsf{S}^{2n-1}_q)$, and
$y=\beta_n(z)$ for some $z\in C(\mathsf{B}^{2n}_q)$ because $\beta_n$ is surjective. Now, from
$$
\sigma_n(\rho_n(z))=\omega_{n-1}(\beta_n(z)))=\omega_{n-1}(y)=\sigma_n(x) ,
$$
we deduce that $x-\rho_n(z)\in\ker(\sigma_n)$, so
$x-\rho_n(z)\in\mathrm{im}(\rho_n)$
by Lemma \ref{lemma:kernels}.
Therefore, $x-\rho_n(z)=\rho_n(t)$ for some $t\in C(\mathsf{B}^{2n}_q)$, so $x=\rho_n(z+t)\in\mathrm{im}(\rho_n)$, as desired.
\end{proof}

Note that the morphism dual to $\rho_n$ is a K-equivalence because of \ref{ax:contr}.
By virtue of Theorem \ref{prop:finalstatementball},
we can now apply \ref{ax:cwW} to the pushout dual to \eqref{eq:ballpullback} to obtain the following corollary.

\begin{cor}\label{thm:maint}
The maps $\rho_n:C(\mathsf B^{2n}_q)\to\mathcal{T}^{\otimes n}$ and $\omega_{n-1}:C(\mathsf{S}^{2n-1}_q)\to C(\mathsf{S}^{2n-1}_H)$ defined in Lemma~\ref{prop:ball2ball} are K-equivalences.
\end{cor}

As an application of the above result, we obtain:

\begin{cor}\label{cor:generator}
$K_1(C(\mathsf{S}^{2n+1}_H))\cong\Z$ is generated by the unitary
\begin{equation}\label{eq:expr}
U:=1+(s_n-1)\prod_{k=0}^{n-1}(1-s_ks_k^*) .
\end{equation}
\end{cor}

\begin{proof}
We know \cite[Theorem 2.3]{adht17} that $K_1(C(\mathsf{S}^{2n+1}_q))\cong\Z$ is generated by the class of the unitary $V:=S_{e_{nn}}+(1-S_{e_{nn}}S_{e_{nn}}^*)$,
and $\omega_n(V)$ generates $K_1(C(\mathsf{S}^{2n+1}_H))$ because
$\omega_n$ is a K-equivalence
Now, using \eqref{eq:sphererel}, we compute
\[
\omega_n(S_{e_{nn}})=
s_ns_ns_n^*\prod_{k=0}^{n-1}(1-s_ks_k^*)
=s_n\prod_{k=0}^{n-1}(1-s_ks_k^*) .
\]
Applying again \eqref{eq:sphererel}, we see that
$\omega_n(V)=U$.
\end{proof}

Pasting the pullback diagrams \eqref{eq:CWCPnq} and \eqref{eq:ballpullback}, we get a new pullback diagram that dualizes to the pushout diagram
\begin{center}
\begin{tikzpicture}[scale=1.5,baseline={([yshift=-3pt]current bounding box.center)}]

\node (A) at (1,2) {$\CPq_{q}^{n}$};
\node (B) at (0,1) {$\CPq_{q}^{n-1}$};
\node (C) at (2,1) {$\mathsf{D}_q^{\times n}$};
\node (D) at (1,0) {$\mathsf{S}_{H}^{2n-1}$};

\path[->]
     	(B) edge[right hook->](A)
	(C) edge (A)
	(D) edge (B)
	(D) edge[right hook->] (C);

\end{tikzpicture}.
\end{center}
As a byproduct, we thus get another strict CW-structure of $\CP_q^n$, where the cofibrations are the same as in \eqref{VSskel}, but at each step instead of a 
quantum ball we attach a quantum polydisk.

We end this section by noticing that the Hong--Szyma{\'n}ski quantum balls and quantum polydisks, although K-equivalent, are not isomorphic.

\begin{prop}\label{prop:noniso}
If $n\geq 1$ and $k\geq 2$, then
the C*-algebras $C(\mathsf{B}^{2n}_q)$ and $C(\mathsf{D}^{\times k}_q)=\mathcal{T}^{\otimes k}$ are not isomorphic.
\end{prop}

\begin{proof}
Suppose that $C(\mathsf{B}^{2n}_q)\cong\mathcal{T}^{\otimes k}$.
Then, such an isomorphism composed with the symbol map gives a surjective *-homomorphism
$C(\mathsf{B}^{2n}_q)\to C(S^1)^{\otimes k}$ that factors through the canonical quotient map $\pi:C(\mathsf{B}^{2n}_q)\to C(\mathsf{B}^{2n}_q)/I$, where $I$ is the closed 
ideal generated by commutators. Therefore, it suffices to show that $C(\mathsf{B}^{2n}_q)/I\cong C(S^1)$, because there is no closed embedding of $(S^1)^{\times k}
$ in $S^1$.

To this end, recall that $C(\mathsf{B}^{2n}_q)$ is the graph C*-algebra of the graph \eqref{eq:ballgraph}.
Let $0\leq j<n$.
It follows from
\[
\sum_{k=j}^nS_{\tilde e_{jk}}S_{\tilde e_{jk}}^*=P_{\tilde v_j}=S_{\tilde e_{jj}}^*S_{\tilde e_{jj}}
\]
that $\sum_{k=j+1}^nS_{\tilde e_{jk}}S_{\tilde e_{jk}}^*=[S_{\tilde e_{jj}}^*,S_{\tilde e_{jj}}]$, so
\[
\pi\Big(\sum_{k=j+1}^nS_{\tilde e_{jk}}S_{\tilde e_{jk}}^*\Big)=
\sum_{k=j+1}^n\pi(S_{\tilde e_{jk}})\pi(S_{\tilde e_{jk}})^*=0.
\]
Now, as all summands are positive, we infer that
$\pi(S_{\tilde e_{jk}})\pi(S_{\tilde e_{jk}})^*=0$, and consequently $\pi(S_{\tilde e_{jk}})=0$, for all $0\leq j<k\leq n$. Hence,
\[
\pi(S_{\tilde e_{kk}})^*\pi(S_{\tilde e_{kk}})=\pi(P_{\tilde v_k})=\pi(S_{\tilde e_{jk}})^*\pi(S_{\tilde e_{jk}})=0
\]
for all $0\leq j<k\leq n$, so $\pi(S_{\tilde e_{kk}})=0$ for all $1\leq k\leq n$. It follows that the only generating partial isometry not in the kernel of $\pi$ is $S_{\tilde e_{00}}$, thus proving that $C(\mathsf{B}^{2n}_q)/I\cong C(S^1)$.
\end{proof}

\subsection{The weak cellular equivalence of 
\texorpdfstring{$\CPq_{q}^{n}$ and \texorpdfstring{$\CPq_{H}^{n}$}{CPHn}}{CPqn}}\label{sec:avb3}
Now we are ready to prove a theorem relating the K-theory of $\CPq_{q}^{n}$ with the K-theory of $\CPq_{H}^{n}$. We will adopt the same notation as in the 
previous section. In particular, $\rho_n$ and $\omega_{n-1}$ are the *-homomorphisms defined in Lemma~\ref{prop:ball2ball}.

\begin{thm}\label{compar}
The *-homomorphism
\[
(\omega_n)^{U(1)}:C(\CPq_q^n)\longrightarrow C(\CPq_{H}^{n})
\]
is a K-equivalence.
\end{thm}

\begin{proof}
Consider the cube diagram
\begin{equation}\label{eq:inductionomega}
\begin{tikzpicture}[scale=1.8,baseline=(current bounding box.center)]

\node (A) at (1.1,2) { $C(\mathsf{S}^{2n+1}_q)_\bullet$ };
\node (B) at (0,1) { $C(\mathsf{S}^{2n-1}_q)_\bullet$ };
\node (C) at (2.2,1) { $C(\mathsf{B}^{2n}_q)\otimes C(S^1)_\bullet$ };
\node (D) at (1.1,0) { $C(\mathsf{S}^{2n-1}_q)\otimes C(S^1)_\bullet$ };
\node (E) at (5.6,2) { $C(\mathsf{S}^{2n+1}_H)_\bullet$ };
\node (F) at (4.5,1) { $C(\mathsf{S}^{2n-1}_H)_\bullet\otimes\mathcal{T}_\bullet$ };
\node (G) at (6.7,1) { $\mathcal{T}^{\otimes n}\otimes C({S}^1)_\bullet$ };
\node (H) at (5.6,0) { $C(\mathsf{S}^{2n-1}_H)\otimes C({S}^1)_\bullet$ };

\path[-To,font=\footnotesize] 
      (A) edge (B)
			(A) edge (C)
			(B) edge (D)
			(C) edge (D)
			(E) edge (F)
			(E) edge (G)
			(F) edge (H)
			(G) edge (H)
			(A) edge[bend left=10] node[above] {$\omega_n$} (E)
			(B) edge[bend right=20] node[below,pos=0.6] {$\omega_{n-1}\otimes 1_{\mathcal{T}}$} (F)
			(C) edge[bend left=20] node[above,pos=0.4] {$\rho_n\otimes\id$} (G)
			(D) edge[bend right=10] node[below] {$\omega_{n-1}\otimes\id$} (H);
			
\end{tikzpicture},
\end{equation}
where the left square is the pullback diagram \eqref{eq:VSspheresPB} and the right square is the pullback diagram \eqref{eq:MpullT} for $k=0$. We now check that \eqref{eq:inductionomega} is a commutative diagram. The four faces to be checked are:
\begin{center}
\setlength{\tabcolsep}{0pt}%
\begin{tabular}{cc}
\begin{tikzpicture}[baseline=(current bounding box.center),xscale=0.85]

\node (A) at (0,2) { $C(\mathsf{S}^{2n+1}_q)$ };
\node (B) at (5,2) { $C(\mathsf{S}^{2n+1}_H)$ };
\node (C) at (0,0) { $C(\mathsf{B}^{2n}_q)\otimes C({S}^1)$ };
\node (D) at (5,0) { $\mathcal{T}^{\otimes n}\otimes C({S}^1)$ };

\path[-To,font=\scriptsize] 
      (A) edge node[above] {$\omega_n$} (B)
			(A) edge node[left] {$(r_n\otimes\id)\circ\delta$}  (C)
			(B) edge node[right] {$\phi\circ p_2$} (D)
			(C) edge node[below] {$\rho_n\otimes\id$} (D);
			
\end{tikzpicture} &
\begin{tikzpicture}[baseline=(current bounding box.center)]

\node (A) at (0,2) { $C(\mathsf{B}^{2n}_q)\otimes C(S^1) $ };
\node (B) at (5,2) { $\mathcal{T}^{\otimes n}\otimes C({S}^1)$ };
\node (C) at (0,0) { $C(\mathsf{S}^{2n-1}_q)\otimes C({S}^1)$ };
\node (D) at (5,0) { $C(\mathsf{S}^{2n-1}_H)\otimes C({S}^1)$ };

\path[-To,font=\scriptsize] 
      (A) edge node[above] {$\rho_n\otimes\id$} (B)
			(A) edge node[left] {$\beta_n\otimes\id$} (C)
			(B) edge node[right] {$\sigma_n\otimes\id$} (D)
			(C) edge node[below] {$\omega_{n-1}\otimes\id$} (D);
			
\end{tikzpicture} 
 \\
~\\[-5pt]
\begin{tikzpicture}[baseline=(current bounding box.center),xscale=0.9]

\node (A) at (0,2) { $C(\mathsf{S}^{2n+1}_q)$ };
\node (B) at (5,2) { $C(\mathsf{S}^{2n+1}_H)$ };
\node (C) at (0,0) { $C(\mathsf{S}^{2n-1}_q)$ };
\node (D) at (5,0) { $C(\mathsf{S}^{2n-1}_H)\otimes\mathcal{T}$ };

\path[-To,font=\scriptsize] 
      (A) edge node[above] {$\omega_n$} (B)
			(A) edge node[left] {$\beta_n\circ r_n$}  (C)
			(B) edge node[right] {$p_1$} (D)
			(C) edge node[below] {$\omega_{n-1}\otimes 1_{\mathcal{T}}$} (D);
			
\end{tikzpicture} &
\begin{tikzpicture}[baseline=(current bounding box.center)]

\node (A) at (0,2) { $C(\mathsf{S}^{2n-1}_q)$ };
\node (B) at (5,2) { $C(\mathsf{S}^{2n-1}_H)\otimes \mathcal{T}$ };
\node (C) at (0,0) { $C(\mathsf{S}^{2n-1}_q)\otimes C({S}^1)$ };
\node (D) at (5,0) { $C(\mathsf{S}^{2n-1}_H)\otimes C({S}^1)$ };

\path[-To,font=\scriptsize] 
      (A) edge node[above] {$\omega_{n-1}\otimes 1_{\mathcal{T}}$} (B)
			(A) edge node[left] {$\delta$} (C)
			(B) edge node[right] {$\phi\circ(\id\otimes\sigma_1)$} (D)
			(C) edge node[below] {$\omega_{n-1}\otimes\id$} (D);
			
\end{tikzpicture}
\end{tabular}
\end{center}
Here $p_1$ and $p_2$ are the *-homomorphisms in Lemma \ref{lemma:3.1}, and the *-automorphism $\phi$ is given by $a\otimes f\mapsto a_{(0)}\otimes a_{(1)}f$. The second diagram is simply \eqref{eq:ballpullback} tensored everywhere by $C({S}^1)$. Commutativity of the fourth diagram is also obvious. Indeed, $\phi$ applied to the image of $\omega_{n-1}\otimes 1_{\mathcal{T}}$ coincides with the coaction on the first leg, so the commutativity of the diagram is equivalent to $\delta\circ\omega_{n-1}=(\omega_{n-1}\otimes\id)\circ\delta$, which is tantamount to the $U(1)$-equivariance of~$\omega_{n-1}$.
Commutativity of the remaining two diagrams can be explicitly checked on generators. Concerning the first diagram,
$$
(\phi\circ p_2)(s_i)  =\bigg\{\!\begin{array}{ll}
t_i\otimes u &\forall\;i=0,\ldots,n-1, \\
1\otimes u &\text{if }i=n.
\end{array}
$$
Then, using the explicit formulas in Lemma \ref{prop:ball2ball}, we obtain:
$$
(\phi\circ p_2\circ\omega_n)(S_{e_{ij}})=\bigg\{\!\begin{array}{ll}
\rho_n(S_{e_{ij}})\otimes u &\text{if }i\neq n, \\
\rho_n(P_{v_n})\otimes u &\text{if }i=n.
\end{array}
$$
On the other hand,
$$
\big((r_n\otimes\id)\circ\delta\big)(S_{e_{ij}})=\bigg\{\!\begin{array}{ll}
S_{e_{ij}}\otimes u &\text{if }i\neq n, \\
P_{v_n}\otimes u &\text{if }i=n,
\end{array}
$$
so $\phi\circ p_2\circ\omega_n=(\rho_n\otimes\id)\circ(r_n\otimes\id)\circ\delta$.
Finally, to show the commutativity of the third diagram, using \eqref{eq:sphererel}, we compute
$$
(p_1\circ\omega_n)(S_{e_{ij}})=\bigg\{\!\begin{array}{ll}
\rho_n(S_{e_{ij}})\otimes 1 & \forall\;0\leq i\leq j<n, \\
0 &\text{if }j=n.
\end{array}
$$
On the other hand,
$$
(\beta_n\circ r_n)(S_{e_{ij}})=\bigg\{\!\begin{array}{ll}
S_{e_{ij}} &\forall\;0\leq i\leq j<n, \\
0 &\text{if }j=n,
\end{array}
$$
so $p_1\circ\omega_n=\omega_{n-1}\circ\beta_n\circ r_n\otimes 1_{\mathcal{T}}$.

Furthermore, since all the *-homomorphisms in the commutative diagram \eqref{eq:inductionomega} are $U(1)$-equivariant, we get the following commutative diagram between fixed-point subalgebras:
\begin{equation}\label{eq:inductioncpn}
\begin{tikzpicture}[scale=1.8,baseline={([yshift=-5pt]current bounding box.center)}]

\node (A) at (1,2) { $C(\CPq_q^n)$ };
\node (B) at (0,1) { $C(\CPq_q^{n-1})$ };
\node (C) at (2,1) { $C(\mathsf{B}^{2n}_q)$ };
\node (D) at (1,0) { $C(\mathsf{S}^{2n-1}_q)$ };
\node (E) at (5.5,2) { $C(\CPq^n_{H})$ };
\node (F) at (4.5,1) { $\smash{\big(}C(\mathsf{S}^{2n-1}_H)\otimes\mathcal{T}\smash{\big)^{{U}(1)}}$ };
\node (G) at (6.5,1) { $\mathcal{T}^{\otimes n}$ };
\node (H) at (5.5,0) { $C(\mathsf{S}^{2n-1}_H)$ };

\path[-To,font=\footnotesize] 
      (A) edge (B)
			(A) edge (C)
			(B) edge (D)
			(C) edge (D)
			(E) edge (F)
			(E) edge (G)
			(F) edge (H)
			(G) edge (H)
			(A) edge[bend left=10] node[above] {$(\omega_n)^{U(1)}$} (E)
			(B) edge[bend right=20] node[below,pos=0.6] {$(\omega_{n-1}\otimes 1_{\mathcal{T}})^{U(1)}$} (F)
			(C) edge[bend left=20] node[above,pos=0.4] {$\rho_n$} (G)
			(D) edge[bend right=10] node[below] {$\omega_{n-1}$} (H);
			
\end{tikzpicture}.
\end{equation}
Using this diagram, by induction on $n\geq 0$, we prove that $(\omega_n)^{U(1)}$ is a K-equivalence.
The first step is obvious as $(\omega_0)^{U(1)}$ is an isomorphism ($C(\CPq_q^0)\cong \C\cong C(\CPq_H^0)$). Assume now that $(\omega_{n-1})^{U(1)}$ is a K-equivalence for some $n\geq 1$.
Then the *-homomorphism $(\omega_{n-1}\otimes 1_{\mathcal{T}})^{U(1)}$ is also a K-equivalence because it is the composition of 
$(\omega_{n-1})^{U(1)}$ with the K-equivalence~\eqref{eq:tubweq}. 
Finally, $\rho_n$ and $\omega_{n-1}$ are K-equivalences by Theorem~\ref{thm:maint},
so \ref{ax:weq} applied to the diagram dual to \eqref{eq:inductioncpn} implies that $(\omega_n)^{U(1)}$ is a K-equivalence.
\end{proof}

\begin{cor}
Dualizing the morphisms in Theorem~\ref{compar},
we obtain the following weak cellular equivalence
\begin{equation}\label{weqsk}
\begin{tikzcd}
\CPq^0_H \arrow[r, rightarrowtail]
\arrow[d, rightarrow, "\sim" {xshift=2mm, anchor=center, rotate=90}] 
 & \CPq^1_H \arrow[r, rightarrowtail] 
 \arrow[d, rightarrow, "\sim" {xshift=2mm, anchor=center, rotate=90}] 
 & \ldots \arrow[r, rightarrowtail] &
 \CPq^{n-1}_H \arrow[r, rightarrowtail] 
 \arrow[d, rightarrow, "\sim" {xshift=2mm, anchor=center, rotate=90}]
 & \CPq^{n}_H
 \arrow[d, rightarrow, "\sim" {xshift=2mm, anchor=center, rotate=90}]\\
\CPq^0_q \arrow[r, hookrightarrow] & \CPq^1_q \arrow[r, hookrightarrow] & \ldots 
\arrow[r, hookrightarrow] &\CPq_q^{n-1} \arrow[r, hookrightarrow] & \CPq^{n}_q
\end{tikzcd}
\end{equation}
between the weak filtration by skeleta \eqref{skel} and the strict filtration by skeleta \eqref{VSskel}.
\end{cor}

\pagebreak

\begin{proof}
From \eqref{eq:inductioncpn} we extract the bottom-right triangle in the following diagram:
\begin{center}
\begin{tikzpicture}[baseline=(current bounding box.center),scale=1.7]

\node (A) at (4,0) { $C(\CPq^n_q)$\,.\!\! };
\node (B) at (4,2) { $C(\CPq^n_H)$ };
\node (C) at (0,0) { $C(\CPq^{n-1}_q)$ };
\node (D) at (2,1) { $\big(C(\mathsf{S}^{2n-1}_H)\otimes\mathcal{T}\big)^{U(1)}$ };
\node (E) at (0,2) { $C(\CPq^{n-1}_H)$};

\path[->,font=\scriptsize] 
	(A) edge node[below,sloped] {$\sim$} (B)
	(A) edge[->>] (C)
	(B) edge[->>] (D)
	(C) edge node[above,sloped] {$\sim$} node[below right,inner sep=0pt] {$(\omega_{n-1}\otimes 1_{\mathcal{T}})^{U(1)}$} (D)
	(C) edge node[below,sloped] {$\sim$} node[left] {$(\omega_{n-1})^{U(1)}$} (E)
	(E) edge node[below,sloped] {$\sim$} node[above right,inner sep=0pt] {$(\id\otimes 1_{\mathcal{T}})^{U(1)}$} (D);
			
\end{tikzpicture}
\end{center}
Here the *-homomorphism $(\id\otimes 1_{\mathcal{T}})^{U(1)}$ is a K-equivalence by Corollary~\ref{cor:tub}, and the left triangle is clearly commutative. Passing to the opposite category, we get
\begin{center}
\begin{tikzpicture}[baseline=(current bounding box.center),scale=1.5]

\node (A) at (4,0) {$\CPq^n_q$\,.\!\!};
\node (B) at (4,2) {$\CPq^n_H$};
\node (C) at (0,0) {$\CPq^{n-1}_q$};
\node (D) at (2,1) {$\tub{n-1}{n}$};
\node (E) at (0,2) {$\CPq^{n-1}_H$};

\path[<-,font=\scriptsize] 
	(A) edge node[below,sloped] {$\sim$} (B)
	(A) edge[<-right hook] (C)
	(B) edge[<-right hook] (D)
%	(C) edge node[above,sloped] {$\sim$} (D)
	(C) edge node[below,sloped] {$\sim$} (E)
	(E) edge node[below,sloped] {$\sim$} (D);
			
\end{tikzpicture}
\end{center}
Here one arrow is omitted because it is not relevant. In $\Ho(\CQS)$, the pushout of the top span defines the weak cofibration as in \eqref{twosquares}, 
so we get the commutative squares in~\eqref{weqsk}.
\end{proof}

\subsection{Importing the ring structure from \texorpdfstring{$K^0(\CP^n)$}{K0CPn} to \texorpdfstring{$K^0(\CPq_q^n)$}{K0CPqn} and \texorpdfstring{$K^0(\CPq_H^n)$}{K0CPHn}}\label{sec:ATH}
Every unital $U(1)$-C*-algebra $A$ has a decomposition
$A=\overline{\bigoplus_{k\in\Z}A_k}$ into spectral subspaces
\[
A_k:=\big\{a\in A\,\big|\,\gamma_\lambda(a)=\lambda^ka,\;\forall\;\lambda\in U(1)\big\} .
\]
Here $A_0$ is the C*-subalgebra of $U(1)$-invariant elements, and the $A_k$'s are $A_0$-bimodules.
%The grading is called \emph{strong} if .
When the $U(1)$-action is free \cite{ellwood2000new,baum2017free}, 
they are finitely generated projective left modules~\cite[Theorem~1.2]{ky13}, and we interpret them as the 
modules of sections of  line bundles over the quantum space described by the C*-algebra $A_0$. 

The
\emph{Peter--Weyl subalgebra} 
$\mathcal{P}_{U(1)}(A):=\bigoplus_{k\in\Z}A_k$ of~$A$
is equipped with a right $\mathcal{O}(U(1))$-coaction dual to the $U(1)$-action. The $U(1)$-action on $A$ is free if and only if the dual coaction is
Hopf--Galois~\cite{baum2017free}, which here is tantamount to the \emph{strong grading} property $A_j\cdot A_k=A_{j+k}$ for 
all $j,k\in\Z$~\cite[Theorem~8.1.7]{montgomery1993hopf}. 
If $A$ and $B$ are unital C*-algebras with free $U(1)$-actions, and $\omega:A\to B$ is a $U(1)$-equivariant (i.e.~grading preserving) *-homomorphism,
then the induced map $(\omega^{U(1)})_*:K_0(A_0)\to K_0(B_0)$ satisfies
\begin{equation}\label{eq:indmap}
(\omega^{U(1)})_*([A_k])=[B_k] .
\end{equation}
(This is a special case of \cite[Theorem~1.1]{HM16}.)

The diagonal $U(1)$-action on the odd sphere $S^{2n+1}$ yields a free action on the commutative C*-algebra~$C(S^{2n+1})$. Moreover, 
the gauge action on $C(\mathsf S^{2n+1}_q)$ is free because the graph $\Lambda_{2n+1}$ is finite and withouth sinks 
\cite[Theorem~3.15]{hazrat2013graded}. The $U(1)$-action on $C(\mathsf S^{2n+1}_H)$ is also free, as proved in 
\cite[Section~3.3.1]{hnpsz18}. To abbreviate notation, let us define
$\mathrm{L}_k:=C( S^{2n+1})_k$, $L_k:=C(\mathsf S^{2n+1}_q)_k$  and $\widetilde{L}_k:=C(\mathsf S^{2n+1}_H)_k$ for all $k\in\Z$.
Note that $\mathrm{L}_k$ is the  module of continuous sections of the following line bundle on $\CP^n$:
\[
\mathcal{L}_k:=\begin{cases}
\mathcal{L}_1^{\otimes k} & \text{if }k\geq 0, \\
\mathcal{L}_{-1}^{\otimes |k|} & \text{if }k<0 .
\end{cases}
\]
Here $\mathcal{L}_{-1}$ is the tautological line bundle, $\mathcal{L}_1$ is its dual, and $\mathcal{L}_1^{\otimes 0}$ is the trivial line bundle.

Free generators of the abelian group $K^0(\CP^n)$ were implicitly provided by Atiyah and Todd in \cite{AT60}.
Then Adams \cite{a-jf62} explicitly presented the $K^0$-ring as
\begin{equation}\label{eq:11}
K^{0}(\mathbb{C}P^n)=\mathbb{Z}[x]/(x^{n+1}) .
\end{equation}
Following \cite[Theorem 2.5, p.~214]{karoubi2009k},  we adopt the choice  $x:=[1]-[\mathcal{L}_1]$.
On the other hand, free generators of $K_0(C(\CPq^n_q))$ are given by the classes of the vertex projections $P_{v_0},\ldots,P_{v_n}$ of the graph of 
$\Lambda_{2n+1}$ in~\eqref{lambda} (see~\cite[Proposition~6.10]{d2026k}), which allows us to define the following homomorphism
of abelian groups:
\begin{equation}\label{eq:varphi}
\forall\;0\leq j\leq n\colon K_0(C(\CPq^n_q))\ni [P_{v_j}]\stackrel{\varphi}{\longmapsto} x^j-x^{j+1}  \in K^0(\CP^n) .
\end{equation}
We are now ready to claim the  key technical result of this section:
\begin{thm}\label{thm:varphi}
The homomorphism $\varphi$ of abelian groups \eqref{eq:varphi}
is an isomorphism satisfying:
\begin{equation}\label{desired}
\varphi([L_k])=[\mathrm{L}_k] \qquad\forall\;k\in\Z .
\end{equation}
\end{thm}
\begin{proof}
For starters, $\{1-x,x-x^2,\ldots,x^{n-1}-x^n,x^n\}$ is clearly a set of free generators of $K^0(\C P^n)$ because it is obtained from the
standard set  $\{1,x,\ldots,x^n\}$ of free generators by an 
invertible integer matrix. Hence, $\varphi$ is an isomorphism of abelian groups.

To prove \eqref{desired},  we first recall the formula for the generalized binomial coefficients:
\[
\binom{k}{j}:=\frac{1}{j!}\,\prod_{i=0}^{j-1}(k-i),  \quad j\geq 0,\;k\in\Z.
\]
Here, by convention, $0!:=1$ and the empty product is~$1$. 
Note also that $\binom{k}{j}$ is always an integer because of \cite[(5.14)]{graham1994concrete}.
The expansions in \cite[Section 7.2]{d2026k} for $k\geq 1$ and for $k < 1$, using  \cite[(5.14)]{graham1994concrete}, can be reformulated
into the single formula:
\begin{equation}\label{eq:above}
[L_k] =\sum_{j=0}^{n}(-1)^j\binom{k-1}{j}[P_{v_j}] \qquad
\forall\;k\in\Z.
\end{equation}
On the other hand, the generalized binomial theorem \cite[(5.13)]{graham1994concrete} gives
\[
(1-x)^{k-1}=\sum_{j=0}^n
(-1)^j\binom{k-1}{j}x^j\mod x^{n+1}
\]
for all $k\in\Z$. Hence, for all $k\in\Z$,
\begin{align}
[\mathrm{L}_k] &=[\mathrm{L}_{1}]^k \notag\\
&=(1-x)^{k-1}(1-x) \notag \\
&=\sum_{j=0}^n
(-1)^j\binom{k-1}{j}(x^j-x^{j+1}) \notag \\
&=\sum_{j=0}^n
(-1)^j\binom{k-1}{j}\varphi([P_{v_j}]) \notag \\
&=\varphi\Big(\sum_{j=0}^n
(-1)^j\binom{k-1}{j}[P_{v_j}] \Big) \notag \\
&=\varphi([L_k]) .\label{last}
\end{align}
Here, in the last equality, we used \eqref{eq:above}. 
\end{proof}

Theorem~\ref{thm:varphi} allows us to transport the ring structure of $K^0(\CP^n)$ to $K_0(\CPq^n_q)$ using the isomorphism $\varphi$ of abelian groups,
which yields:
\begin{cor}\label{Theorem5.10}
There exists a ring structure on $K_{0}(C(\CPq_q^{n}))$ whose multiplication satisfies
\begin{equation}\label{concA}
[L_j][L_k]=[L_{j+k}] \qquad \forall\; j,k\in\Z .
\end{equation}
\end{cor}

Finally, employing Theorem~\ref{compar}, we can claim:
\begin{cor}\label{multmult}
The set  $\{[\widetilde{L}_0],[\widetilde{L}_1],\ldots,[\widetilde{L}_n]\}$  is a set of free generators of~$K_{0}(C(\CPq_H^{n}))$. Moreover,
 there exists a ring structure on 
$K_{0}(C(\CPq_H^{n}))$ whose multiplication satisfies
\begin{equation}\label{concB}
[\widetilde{L}_j][\widetilde{L}_k]=[\widetilde{L}_{j+k}]
\qquad \forall\; j,k\in\Z .
\end{equation}
\end{cor}
\begin{proof}
As a special case of \eqref{eq:indmap}, the map $(\omega_n)^{U(1)}$ in Theorem~\ref{compar} satisfies
\begin{equation}\label{ltol}
(\omega_n)^{U(1)}_*:K_0(C(\CPq_q^n))\ni [L_k]\longmapsto [\widetilde L_k]\in K_0(C(\CPq_{H}^{n})) 
\end{equation}
for all $k\in\Z$.
On the other hand, we already know that
\[
K_{0}(C(\CPq_q^{n}))  =\bigoplus_{k=0}^{n}\Z[L_{k}] .
\]
(This was first proved  in 
\cite{d2010bounded} by means of the index pairing,
and is also a special case of a more general result proved in the context of graph C*-algebras \cite[Section~7.2]{d2026k}.)
Hence,
\[
K_{0}(C(\CPq_H^{n}))  =\bigoplus_{k=0}^{n}\Z[\widetilde{L}_{k}] .
\]
Finally, we use the isomorphism of abelian groups $(\omega_n)^{U(1)}_*$ to transport the ring structure of 
$K_0(C(\CPq^n_q))$ (see Corollary~\ref{Theorem5.10}) to $K_0(C(\CPq^n_H))$, and to infer \eqref{concB} from \eqref{concA}.
\end{proof}

\addtocontents{toc}{\SkipTocEntry}
\section*{Acknowledgements}
\noindent
This research is part of the EU Staff Exchange project 101086394 {\it Graph Algebras}.
The project is co-financed by the Polish Ministry of Education and Science under the program PMW, grant agreement 5305/HE/2023/2 (Piotr M.~Hajac)
and grant agreement 5448/HE/2023/2 (Tomasz Maszczyk).
We also acknowledge the EU grant H2020-MSCA-RISE-2015-691246 {\it Quantum Dynamics} for supporting the initial phase of this work, which 
was also co-financed by the Polish Ministry of Science and Higher Education from the funds allocated for science in the years 2016-2019  
through the grants W2/H2020/2016/317281 and  328941/PnH/2016 (Piotr M.\ Hajac), W43/H2020/2016/319577  and 329915/PnH/2016 
 (Tomasz Maszczyk),
W30/H2020/2016/319460  and 329390/PnH/2016
(Bartosz Zieli\'nski). 
Finally, a visit to the Institute of Mathematics of the Polish Academy of Sciences (IMPAN) of Francesco D'Andrea and Bartosz Zieli\'nski was  
supported by 
the Simons Foundation Award no.\ 663281 granted to IMPAN
 for the years 2021--2023. Last but not least, we are very grateful to Albert Sheu for his contribution to earlier versions of this paper involving
 groupoid presentations of studied C*-algebras.

\end{document}